\documentclass[letterpaper,11pt,reqno]{amsart}
\usepackage[portrait,margin=3.3cm]{geometry}
\usepackage{mathtools}
\usepackage{microtype}
\usepackage{comment}
\usepackage{mathtools} 
\usepackage{mathrsfs,xfrac} 
\usepackage[pagebackref, colorlinks=true,linkcolor=blue,citecolor=blue,urlcolor=blue]{hyperref} 
\usepackage[mathlines]{lineno} 
\usepackage{amsmath,amssymb,amsthm,amsfonts,amsbsy,latexsym,dsfont,color} 
\usepackage[numeric,initials,nobysame]{amsrefs} 
\usepackage[foot]{amsaddr}
\usepackage{upref,halloweenmath}

\usepackage[utf8]{inputenc}
\usepackage[english]{babel}

\usepackage[textsize=tiny]{todonotes}
\usepackage{regexpatch}
\makeatletter
\xpatchcmd{\@todo}{\setkeys{todonotes}{#1}}{\setkeys{todonotes}{inline,#1}}{}{}
\makeatother

\usepackage{etoolbox}          
\newcommand*\linenomathpatch[1]{%
  \cspreto{#1}{\linenomath}%
  \cspreto{#1*}{\linenomath}%
  \csappto{end#1}{\endlinenomath}%
  \csappto{end#1*}{\endlinenomath}%
}
\newcommand*\linenomathpatchAMS[1]{%
  \cspreto{#1}{\linenomathAMS}%
  \cspreto{#1*}{\linenomathAMS}%
  \csappto{end#1}{\endlinenomath}%
  \csappto{end#1*}{\endlinenomath}%
}

\expandafter\ifx\linenomath\linenomathWithnumbers
  \let\linenomathAMS\linenomathWithnumbers
  \patchcmd\linenomathAMS{\advance\postdisplaypenalty\linenopenalty}{}{}{}
\else
  \let\linenomathAMS\linenomathNonumbers
\fi

\linenomathpatch{equation}
\linenomathpatchAMS{gather}
\linenomathpatchAMS{multline}
\linenomathpatchAMS{align}
\linenomathpatchAMS{alignat}
\linenomathpatchAMS{flalign}

\makeatletter
\patchcmd{\mmeasure@}{\measuring@true}{
  \measuring@true
  \ifnum-\linenopenaltypar>\interdisplaylinepenalty
    \advance\interdisplaylinepenalty-\linenopenalty
  \fi
  }{}{}
\makeatother

\usepackage{enumerate}

\newtheorem{thm}{Theorem}[section]
\newtheorem{lem}[thm]{Lemma}
\newtheorem{cor}[thm]{Corollary}
\newtheorem{prop}[thm]{Proposition}
\newtheorem{defn}[thm]{Definition}
\newtheorem{rem}[thm]{Remark}

\renewcommand{\le}{\leqslant} 
\renewcommand{\ge}{\geqslant} 
\renewcommand{\leq}{\leqslant} 
\renewcommand{\geq}{\geqslant}

\newcommand{\eps}{\varepsilon}
\newcommand{\To}{\longrightarrow}
\newcommand{\norm}[1]{\left\Vert#1\right\Vert}
\newcommand{\abs}[1]{\left\vert#1\right\vert}

\newcommand{\equald}{\stackrel{\mathrm{d}}{=}}

 \let\gb=\beta \let\gc=\gamma  
                \let\gs=\sigma \let\gt=\tau 
  
 \let\gD=\Delta  \let\gL=\Lambda 
         \let\gS=\Sigma  
                                
\newcommand{\mc}[1]{\mathcal{#1}}
\newcommand{\cA}{\mathcal{A}}

\newcommand{\cH}{\mathcal{H}}
\newcommand{\cK}{\mathcal{K}}
\newcommand{\cN}{\mathcal{N}}
\newcommand{\cP}{\mathcal{P}}

\newcommand{\fh}{\mathfrak{h}}

\newcommand{\bN}{\mathbb{N}}
\newcommand{\bR}{\mathbb{R}}

\DeclareMathOperator{\E}{\mathds{E}}

\newtheorem{theorem}[thm]{Theorem}

\newtheorem{lemma}[thm]{Lemma}

\theoremstyle{remark}

\theoremstyle{plain}
\theoremstyle{definition}

\DeclareMathOperator{\Cov}{Cov}
\DeclareMathOperator{\Var}{Var}

\newcommand{\R}{\mathbb R}
\newcommand{\one}{\mathbf 1}

\newcommand{\ip}[1]{\left \langle #1 \right \rangle}

\renewcommand{\P}{\mathbb{P}}
\newcommand{\diff}{\mathop{}\!\mathrm{d}}
\DeclareMathOperator{\erfc}{erfc}

\newcommand{\defeq}{\vcentcolon=}
\newcommand{\eqdef}{=\vcentcolon}

\makeatletter
\def\l@subsection{\@tocline{2}{0pt}{4pc}{5pc}{}}
\makeatother

\begin{document}

\title[Complexity of Ising spin glasses]{A simple approach to the complexity of mean field Ising spin glasses}

\author[]{Benjamin McKenna$^\star$ \and Qiang Wu$^{\dagger,\ddagger}$}
\address{$^\star$School of Mathematics, Georgia Institute of Technology, Atlanta, GA, United States 
\newline 
\indent 
$^\dagger$School of Mathematics, University of Minnesota, Minneapolis, MN, United States
\newline \indent 
$^\ddagger$Department of Mathematics, University of Oregon, Eugene, OR, United States}
\email{mckenna@math.gatech.edu$^\star$ wuq@umn.edu$^\dagger$ qiangwu@uoregon.edu$^\ddagger$}
\date{November 5, 2026}

\begin{abstract}
The geometry of spherical spin glasses is relatively well understood, thanks in large part to the Kac–Rice formula. However, this formula does not directly apply to discrete spin systems, and much less is known mathematically about the geometry of Ising spin glasses. In this paper, we present a simple method for counting local optima and gapped states in a broad class of discrete spin glasses. To illustrate its scope, we compute the annealed complexities (``configurational entropies,'' in the physics language) of mixed Ising spin glasses with potentials, such as external fields, and of the Ghatak--Sherrington model. These computations recover several existing results in the mathematics and physics literature as special cases. We can also compute energy-constrained complexity, finding a critical energy threshold analogous to the spherical case, but for gapped states of pure Ising spin glasses. Finally, we use our method to establish a form of disorder persistence for gapped states, identifying a transition between $3$-spin and $4$-spin glasses.
\end{abstract}

\maketitle

\setcounter{tocdepth}{1}\tableofcontents


\section{Introduction and main results}


\subsection{Overview}
Spin glasses, as random functions in high dimensions, have found broad applications in the modern theory of machine learning and statistics. Understanding their energy landscapes has been a fundamental question, due to its connection to the relaxation of the high dimensional random system and to the behavior of algorithms searching for low energy configurations. Local minima are natural trapping points for energy decreasing dynamics, and sufficiently stable local minima can impede exploration of the configuration space. A basic question is therefore how many local minima a spin glass Hamiltonian has and how the number depends on the parameters of the model. One can further ask how stable these local minima are under perturbations of the disorder.
In this paper, we present a simple, general approach to study these questions for discrete spin glasses, which are much less understood than their continuous spherical counterparts.

\subsection{Continuous vs. discrete spin glass complexity} 

The geometry of \emph{continuous} Gaussian processes is well-established, being for example the subject of several textbooks (\cite{AdlTay2007}, \cite{AzaWsc2009}, and the very recent \cite{ArmAzaDalDeCDelLeoMor2026}). The paradigmatic example is a real-valued Gaussian process on a smooth manifold (usually the unit sphere $\mathbb{S}^{N-1}$ or $\R^N$), whose realizations are almost-surely smooth. One main tool in this area is the Kac--Rice formula, a probabilistic analogue of Federer's co-area formula, which reduces the problem to the study of random matrices.

This philosophy is quite common in the particular case when the Gaussian process is the Hamiltonian of a spin glass. To fix ideas, let $H_N : \mathbb{S}^{N-1} \to \R$ be a spin glass Hamiltonian (see the exact form \eqref{eq:mixed-H} below), and let $\mc{N}_N$ be the random variable storing the number of (Riemannian) critical points of $H_N$. Then one can use the Kac--Rice formula to compute the so-called \emph{annealed complexity}
\[
    \Sigma = \lim_{N \to \infty} \frac{1}{N} \log \E[\mc{N}_N].
\]
One can also count finer geometric properties, such as specifically local minima (rather than all critical points). The complexity $\Sigma$ depends on underlying model parameters; as these vary, one sees a geometric transition between when $\Sigma$ is positive (meaning $H_N$ has many critical points) and when it is zero (meaning $H_N$ has few critical points), and these geometric phase transitions often correspond to algorithmic and physical phase transitions. Seminal papers in this area include \cite{Fyo2004, AufBenCer2013, AufBen2013}; for recent surveys, including many more references, and more on the interplay between geometry, algorithms, and physics, we direct readers to \cite{RosFyo2023, AufMonSub2023}. 

However, much less is understood about the geometry of \emph{Ising} spin glasses, meaning for the moment real-valued Gaussian processes $H_N$ defined on $\{-1,+1\}^N$. It is not even immediately clear what to count; there are no derivatives, hence no critical points in the calculus sense, hence no Kac--Rice formula. The physics literature proposes two solutions to this problem, the first generally applicable to random functions on graphs, and the second specific to spin glasses.

The first: for any (deterministic) real-valued function $f$ on the vertices $V$ of a graph $G$, we say that a vertex $x$ is a local minimum if $f(y) - f(x) \geq 0$ for all $y$ adjacent to $x$, or more generally that $x$ is $\gamma$-gapped if $f(y) - f(x) \geq \gamma$ for some $\gamma \geq 0$ and for all $y$ adjacent to $x$. (Sometimes the term \emph{gapped state}, or \emph{stable local minimum}, is used for configurations that are $\gamma$-gapped for some $\gamma > 0$.) By thinking of $\{-1,+1\}^N$ as the vertex set of the hypercube, we can then count the random number $\mc{N}_N^{\textup{lm}}$ of local minima of $f = H_N$, and then $\Sigma = \lim_{N \to \infty} \frac{1}{N} \log \E[\mc{N}_N^{\textup{lm}}]$ is a reasonable notion of complexity. In physics, if $f$ is a Hamiltonian, then local minima are often called \emph{metastable states} and $\Sigma$ is often called the \emph{configurational entropy}. 

Second, in the special case of spin glasses, for physical reasons the Ising Hamiltonian $H_N : \{-1,+1\}^N$ is associated with a real-valued random function $f_{\textup{TAP}}$ defined on the whole unit cube $[-1,1]^N$ (named after Thouless, Anderson, and Palmer \cite{ThoAndPal1977}). Since $f_{\textup{TAP}}$ is a nice process on a continuous space, one can count its critical points using the Kac--Rice formula. Strictly speaking these are not critical points of the original model, but for some purposes they play a similar role. For recent progress in this direction, see \cite{BelConGen2024,Bou2026}. Since the methods of the second solution are essentially those of continuous spin glasses, we focus on the first. 


\subsection{Main contributions}
We summarize our main contributions as follows. 
\begin{itemize}
    \item \textbf{A general approach to discrete complexity.}
    We establish a variational formula for the annealed complexity
    of configurations satisfying prescribed local stability
    conditions. The method applies to discrete Gaussian models
    whose local fields can be represented using independent
    Gaussian variables and a fixed number of shared Gaussian
    variables. Theorem~\ref{thm:general} gives the general result
    under explicit convergence conditions, while
    Corollary~\ref{thm:main} provides a convenient sufficient
    condition for its application.

    \item \textbf{Application to mixed Ising spin glasses.}
    We compute the annealed complexities of local minima and
    gapped states for mixed Ising spin glasses with
    magnetization-dependent potentials, including external
    fields and Curie--Weiss interactions
    (Theorem~\ref{thm:complexity-p-spin}). For external fields,
    the variational formula reduces to one dimension. Its
    analysis shows that the annealed complexity of local minima
    remains positive at every finite external field and
    identifies a critical gap $\gamma_c>0$, depending on the
    mixture and the field, such that the expected number of
    $\gamma$-gapped states grows exponentially for 
    $0\le\gamma<\gamma_c$, whereas
    for $\gamma>\gamma_c$, no such states exist with high
    probability (Corollary~\ref{cor:p-spin-analysis}).

    \item \textbf{Application to disorder persistence.}
    We derive a variational formula for the annealed complexity
    of configurations that satisfy gap conditions for two
    correlated realizations of a mixed Ising Hamiltonian
    (Theorem~\ref{thm:disorder-persistence}). In particular, let $H=H_N$
    be a pure $p$-spin Hamiltonian in zero external field, let
    $H'$ be an independent copy, and consider
    $H_t=e^{-t}H+\sqrt{1-e^{-2t}}\,H'$, $t\ge0$.
    We find a transition between $p=3$ and $p=4$. For $p=2,3$,
    there is a finite threshold $t_c(p)>0$ such that, at each
    fixed $t>t_c(p)$, with high probability there is no
    configuration that is a local minimum of both $H$ and $H_t$.
    For $p\ge4$, the expected number of common local minima grows
    exponentially at every fixed finite time $t$
    (Corollary~\ref{cor:disorder-persistence-analysis}).

    \item \textbf{Application to Ghatak-Sherrington (GS) model.}
    We compute the annealed complexities of local minima and
    gapped states for the GS model, allowing arbitrary crystal
    and external fields (Theorem~\ref{thm:complexity-GS}).
    In zero external field, we prove an annealed topological
    trivialization transition: the complexity of local minima
    is positive for $D<D_c$ and zero for $D\ge D_c$, for some
    finite crystal field $D_c>0$
    (Proposition~\ref{prop:GS-topo-trivial}). This contrasts with
    the persistence of positive annealed complexity at every
    finite external field in the mixed Ising spin-glass models.

    \item \textbf{Application to energy-constrained complexity.}
    Although our general approach is written to count all local minima (or $\gamma$-gapped states), we show how to extend it in the case of Ising pure $p$-spin glasses to count the number of such states at which the Hamiltonian lies below certain energy levels (Theorem~\ref{thm:energy-constrained}). We find that, for Ising $p$-spin, there exists a critical energy $E_\infty(p,\gamma) \in (0,\infty)$ such that the annealed complexity of $\gamma$-gapped states below level $-E$ is flat for $E < E_\infty$, then strictly decreasing for $E > E_\infty$, analogous to the well-known case of spherical spin glasses. As special cases of our result (see Remarks \ref{rem:BraMoo1981shortranged}, \ref{rem:deOFon1997}, and \ref{rem:AddarioBerryEtAl2019} respectively), we recover formulas of Bray and Moore in the physics literature \cite{BraMoo1981shortranged}, de Oliveira and Fontanari in the physics literature \cite{deOFon1997}, and Addario-Berry, Devroye, Lugosi, and Oliveira in the mathematics literature \cite{AddarioBerryEtAl2019}.
\end{itemize}

\subsection{Sketch of the general approach} We briefly sketch the idea of the general approach as follows. For each spin configuration, we present a new representation for the Gaussian local fields, which become independent conditional on the common mode. We illustrate this representation in the case of Sherrington-Kirkpatrick (SK) model with zero external field, whose Hamiltonian is
\[
    H_N(\sigma)
    =
    \frac{1}{\sqrt N}\sum_{i<j}J_{ij}\sigma_i\sigma_j,
    \qquad
    \sigma\in\{-1,+1\}^N,
\]
where the $J_{ij}$ are independent standard Gaussian variables.
Write $\sigma^{(i)}$ for $\sigma$ with its $i$th spin reversed,
and define the signed local field
$G_i(\sigma)=\tfrac12[H_N(\sigma^{(i)})-H_N(\sigma)]$.
Thus, for $\gamma\ge0$, the configuration $\sigma$ is
$\gamma$-gapped precisely when $G_i(\sigma)\ge\gamma/2$ for
every $i$. We note that local fields also play a central role in pseudolikelihood estimation for spin glasses~\cite{CSW24}.

For each fixed $\sigma$, the fields $G_i(\sigma)$ are centered
Gaussian variables with variance $1-1/N$ and pairwise covariance
$1/N$. Consequently, for $N\ge3$, they admit the joint distributional
representation
\[
    \bigl(G_i(\sigma)\bigr)_{i=1}^N
    \stackrel{d}{=}
    \left(
        \sqrt{1-\frac{2}{N}}\,X_i+\frac{Z}{\sqrt N}
    \right)_{i=1}^N,
\]
where $X_1,\ldots,X_N,Z$ are independent standard Gaussian
variables. This is the SK specialization of the covariance
calculation in Section~\ref{subsec:p-spin-establishing}.

Let $\mathcal N_N$ denote the number of $\gamma$-gapped states,
and let $\Phi$ be the standard Gaussian distribution function.
The joint law of the local fields does not depend on $\sigma$.
Conditioning on $Z$ and integrating over $X_1,\ldots,X_N$
therefore gives
\[
    \mathbb E\mathcal N_N
    =
    2^N\mathbb P\!\left(
        G_i(\sigma)\ge\frac{\gamma}{2}
        \text{ for every }i
    \right)
    =
    2^N\mathbb E_Z\!\left[
        \Phi\!\left(
            \frac{Z/\sqrt N-\gamma/2}{\sqrt{1-2/N}}
        \right)^N
    \right].
\]
The correlated $N$-dimensional probability has thus become a
one-dimensional Gaussian integral.
Lemma~\ref{lem:general-nn-expansion} formalizes this conditioning
step in the general setting.

The scaled Gaussian $Z/\sqrt N$ satisfies a large deviation
principle at speed $N$ with rate function $u^2/2$. Using Varadhan's lemma (see
Section~\ref{sec:general}) then yields
\[
    \lim_{N\to\infty}\frac{1}{N}\log\mathbb E\mathcal N_N
    =
    \log 2+
    \sup_{u\in\mathbb R}
    \left\{
        \log\Phi\!\left(u-\frac{\gamma}{2}\right)
        -\frac{u^2}{2}
    \right\}.
\]
This is the zero-field SK case of
Corollary~\ref{cor:p-spin-analysis}; see also Remark~\ref{rem:recover-SK}.

In the general setting, the local-field distributions and stability conditions may also depend on the empirical profile of the configuration, namely the proportion of spins taking each allowed value. In Theorem~\ref{thm:general}, we give a general result that allows for this kind of dependence; the resulting variational problem optimizes over both the empirical profile and the shared Gaussian variables, and includes an entropy term coming from grouping configurations by their profiles. This result is used as an input for most of the rest of the results in this paper.


\subsection{Related work}
We discuss some previous works related to the discrete complexity problems and its applications. 

\subsubsection{Configurational entropy in the physics literature}

In physics, one usual motivation for configurational entropy is simply that highly nonconvex energy landscapes are supposed to be characteristic of spin glasses; thus, if a given Hamiltonian has positive complexity, this is seen as an \emph{a posteriori} justification that it is a good model for a physical spin glass.

To keep this history at a manageable length, we restrict ourselves to configurational entropy for mean-field Ising spin glasses. To the best of our knowledge, the first result in this area is that the SK model has configurational entropy approximately $0.1992$; this was obtained independently at essentially the same time by Bray and Moore \cite{BrayMoore1980,BraMoo1981}, Tanaka and Edwards \cite{TanEdw1980}, and de Dominicis \emph{et al} \cite{deDGabGarOrl1980}. These computations were later generalized in various ways, for example by adding an external field (both an annealed computation \cite{Dea1994} and a quenched computation under a replica-symmetric Ansatz \cite{Rob1981}); by considering the Ising $p$-spin models both without \cite{GroMez1984,Rie1992} and with \cite{deOFon1997} external field; and by looking for what we would call gapped states \cite{TreAmi1988}. Many of these works were motivated by the observation that metastable states of the SK model can be seen as solutions of the zero-temperature TAP equations (and many of these works were interested in the more general problem of counting solutions to the TAP equations at any temperature). 

In the physics literature, configurational entropy is often computed by using Fourier representations of the Dirac delta function, Hubbard-Stratonovich transformations to linearize quadratic terms, and steepest descent. 

\subsubsection{Configurational entropy in the mathematics literature}

The mathematics literature on configurational entropy is much sparser. To the best of our knowledge, the first mathematical computation of the configurational entropy for the SK model was only carried out recently, in work of Addario--Berry, Devroye, Lugosi, and Oliveira \cite{AddarioBerryEtAl2019}. Dandi, Gamarnik, and Zdeborov{\'a} \cite{DandiGamarnikZdeborova2023} studied spin glasses with even non-Gaussian or sparse couplings, but counted a slightly different object, which one can informally describe as configurations $\sigma \in \{-1,+1\}^N$ with magnetization $\frac{1}{N} \sum_{i=1}^N \sigma_i$ close to $0$, and which are $\gamma$-gapped in most directions rather than all of them (see also Remark \ref{rem:gamma-c}).

\subsubsection{Gapped states and connections to algorithms}

We previously mentioned physical motivations for counting local minima, but there are also algorithmic motivations; for example, local minima are natural ``wells'' for Glauber dynamics. This is especially true if the wells are deep, meaning if the configurations are gapped; thus gapped states have appeared in more recent developments of algorithmic spin glass theory. For example, Sellke \cite{Sel2025} recently gave a very short proof that, in the SK model, the existence of gapped states (as proved by \cite{DandiGamarnikZdeborova2023}) implies exponentially slow mixing of low-temperature Glauber dynamics from a worst-case initialization. 

Recently, one line of research has developed the philosophy that ``local algorithms on spin glasses can find local minima, but cannot find gapped states''; see \cite{BehArpKivZde2022} and \cite{MinSahSaw2024} for conjectures to this effect, and \cite{HuaSel2025} for rigorous results in the framework of low-degree polynomials (for both Ising and spherical spin glasses). This philosophy interacts nicely with the curiously frequent appearance of ``marginal'' local minima in spin glasses (i.e., local minima which are \emph{not} gapped): if one expects algorithms to find local minima or gapped states, and one can prove that they cannot find gapped states, then marginal local minima are all that remains. 

\subsubsection{Disorder chaos and disorder persistence}

Many spin glass models are either known or expected to be ``chaotic,'' in the sense that the system changes substantially after a small change in some parameter (like the external field, or the temperature, or the disorder). ``Disorder chaos'' refers to slight perturbations of the random coefficients of the Hamiltonian in the sense of $H$ and $H_t$ above, and usually asks for a result of the form ``If $\sigma$ is the ground state of $H$, and $\sigma_t$ is the ground state of $H_t$, then for \emph{every} $t > 0$ we have that $\sigma$ and $\sigma_t$ are almost orthogonal.'' This a kind of ``zero-temperature'' disorder chaos; the positive temperature version asks for $\sigma$ and $\sigma_t$ to be drawn from the corresponding Gibbs measures. Questions of this type have attracted longtime interest in the physics \cite{FisHus1986, BraMoo1987, KrzBou2005,Riz2009} and mathematics \cite{Cha2009, Che2013, Cha2014, CheHsiHwaShe2015, CheSen2017, ChePan2018, ArgNewSte2019, ArgHan2020, Cha2023, CheKimSen2024} literatures.

Roughly speaking, one expects such behavior because spin glasses are expected to have exponentially many almost-orthogonal approximate ground states. When the couplings are slightly perturbed, one of these points slightly outcompetes the original ground state to become the new ground state. 

In this paper, we will be interested in a slightly different question: If we track a single configuration $\sigma$ which is a local minimum at time zero, how long can it remain a local minimum, even if it gets outcompeted in terms of energy? It is not obvious that this behavior would be chaotic, since for example all that matters is the relative energies of $\sigma$ and its $N$ neighbors, not the energies of faraway configurations. Indeed, our results (Theorem \ref{thm:disorder-persistence} and Corollary \ref{cor:disorder-persistence-analysis}) show that, in an averaged sense, points can persist as local minima for order-one times. Thus we call this phenomenon ``disorder persistence'' (although it may be unlikely for a \emph{given} configuration to persist this long; see Remark \ref{rem:persistence-ldp}). 

This question seems quite natural to us, but we are not aware of its prior appearance in either the physics or mathematics literatures. However, Dean \cite{Dea1994} studied the related question of counting configurations $\sigma$ which are simultaneously local minima for the SK model with zero external field, and for the same realization of the Hamiltonian with positive external field. One could call this ``external field persistence'' (analogously to ``external field chaos'' as in \cite{PanTal2007}). 

\subsubsection{The Ghatak-Sherrington model} 
The Ghatak-Sherrington (GS) model was introduced in \cite{GhatakSherrington1977} as a generalization of the
SK model incorporating crystal field effects. The spins take values in $\{-1,0,+1\}$, and the
crystal field $D$ controls the energetic cost of the nonzero states. Unlike the SK model, the GS model therefore allows the proportion
of interacting spins to vary. In our sign convention, increasing
$D$ favors the zero state, which does not contribute to the random
interaction. This additional degree of freedom gives rise to a richer phase diagram~\cite{MottishawSherrington1985}. One particularly interesting feature is the phenomenon of
\emph{inverse freezing}. The
calculations in physics~\cite{CrisantiLeuzzi2005}
predict a re-entrant transition line: for suitable fixed crystal
fields, lowering the temperature can drive the system from a
spin glass phase back into a paramagnetic phase. 
These predictions motivate the question of how the competition
between the crystal field and the random interactions is reflected in the geometry of the energy landscape.

More directly related to complexity, da Costa and de
Ara\'ujo~\cite{daCostaDeAraujo2000} studied the zero-temperature TAP
equations of the GS model in zero external field. Using the
Bray-Moore approach, they obtained a saddle point expression for
the exponential growth rate of the disorder-averaged number of
solutions. Their numerical evaluation predicts a non-monotone
dependence on the crystal field: the complexity initially
increases from its SK value, reaches a maximum, and then decreases
to zero at a finite crystal field.

Our contribution concerns the direct counting of discrete local
minima and gapped states. Theorem~\ref{thm:complexity-GS} gives a
rigorous variational formula for their annealed complexities,
allowing both an external field and a prescribed gap. In zero
external field, Proposition~\ref{prop:GS-topo-trivial} establishes
a finite crystal field threshold $D_c>0$ such that the annealed
complexity of local minima is positive for $D<D_c$ and zero for
$D\ge D_c$. This contrasts with the mixed Ising models considered
here, whose annealed complexity of local minima remains positive
at every finite external field; see
Corollary~\ref{cor:p-spin-analysis}. The GS result establishes a rigorous
annealed topological trivialization transition as the crystal
field varies (see also Remark \ref{rem:da-costa-de-araujo} for the relationship between our work and \cite{daCostaDeAraujo2000}). This geometric transition threshold seems different from the temperature dependent thermodynamic transition associated with
inverse freezing. We also mention \cite{Panchenko2005,AuffingerChen2021,ShengWu2024} for other mathematical results on the GS model. 


\subsection{Organization of the paper}

Our general method is stated and proved in Section \ref{sec:general}. Applications to Ising spin glasses, to disorder chaos, to the Ghatak--Sherrington model, and to energy-constrained complexity are proved in Sections \ref{sec:p-spin}, \ref{sec:disorder-persistence}, \ref{sec:GS}, and \ref{sec:energy-constrained}, respectively (although Section \ref{sec:energy-constrained} includes only a proof sketch). Finally, Appendix \ref{sec:comparison} gives a detailed comparison between our method and previous methods in the literature, and Appendix \ref{appx:energy-constrained} gives the proof details for the claims in Section \ref{sec:energy-constrained}.


\subsection*{Acknowledgements}
Part of this work was conducted while the authors were attending the Random Theory 2023 workshop. B.M. was supported in part by NSF grant DMS-1760471 (to Shing-Tung Yau and Horng-Tzer Yau). Most of this work was done while Q.W. was affiliated with the University of Minnesota. Part of this work was carried out while Q.W. was visiting the Georgia Institute of Technology in March 2026, and Q.W. thanks the School of Mathematics for the hospitality. This research was supported in part by grants from the NSF (DMS-2235451) and Simons Foundation (MPS-NITMB-00005320) to the NSF-Simons National Institute for Theory and Mathematics in Biology (NITMB).


\subsection*{Disclosure of AI usage} 
This project started in summer 2023. All the complexity computations were carried out by the authors without the use of artificial intelligence. While preparing the manuscript, we used GPT-5.6 and GPT-6 to check our proofs for errors (it found some minor repairable issues in various lemmas throughout the paper), and to refine some of the results. Specifically, part of the proofs in Section \ref{subsubsec:addario-berry-energy-constrained} came from GPT; additionally, in Proposition \ref{prop:GS-topo-trivial}, the authors originally had zero complexity for large enough $D$, and GPT-6 proposed a refinement that went all the way down to $D_c$. The authors have carefully checked and evaluated the AI-generated proposals.


\section{The general method}
\label{sec:general}


\subsection{The general setting}\label{subsec:framework}
Let $\cA$ be a fixed finite alphabet of size $k$, whose elements we will call \emph{types}. Usually in practice we will take either $\cA = \{-1,1\}$ or $\cA = \{-1,0,1\}$. For $\sigma\in\cA^N$ with $N\ge 1$, define the empirical profile
\begin{align}\label{eq:profile}
\pi_a(\sigma)=\frac{|\{i\in [N]:\sigma_i=a\}|}{N},
\qquad a\in\cA,
\end{align}
which is a point of the probability simplex $\Delta_k=\{\pi\in[0,1]^k:\sum_{a}\pi_a=1\}$. The profiles realized at size $N$ form the grid $\cP_N=\{\pi\in\Delta_k:N\pi_a\in\bN\text{ for all } a \in \cA\}$.  

Let 
\begin{equation}
\label{eqn:entropy}
\cH(\pi)=-\sum_{a \in \cA}\pi_a\log\pi_a
\end{equation}
be the entropy of the profile, with convention $0\log0=0$.
The key objects of the framework are, for fixed dimensions $d\ge1$ and $r\ge0$, local Gaussian vectors $Y_i(\sigma)\in\R^d$, $1\le i\le N$, attached to every configuration, and Borel stability sets $\cK_{a,N}(\pi)\subseteq\R^d$ for every type $a$ and profile $\pi$.  The object of interest is the number of configurations such that the associated local fields belong a specified set:
\begin{align}\label{eq:count}
\cN_N=\big|\bigl\{\sigma\in\cA^N:\,Y_i(\sigma)\in\cK_{\sigma_i,N}(\pi(\sigma))\text{ for all }i\bigr\}\big|.
\end{align}
For a broad class of spin glass models, their associated Gaussian local fields $\{Y_i(\sigma)\}_{i \in [N], \sigma \in \mc{A}^N}$ have (for each $N$, each $\pi \in \mathcal{P}_N$, and each $\sigma$ with $\pi(\sigma) = \pi$) the distribution
\begin{align}\label{eq:CF}
\bigl\{Y_i(\sigma)\bigr\}_{i=1}^N
\equald
\Bigl\{\mu_{\sigma_i,N}(\pi)+A_{\sigma_i,N}(\pi)\,X_i+\tfrac1{\sqrt N}B_{\sigma_i,N}(\pi)\,Z\Bigr\}_{i=1}^N,
\end{align}
where $X_1,\dots,X_N$ are independent $\cN(0, I_d)$ vectors, $Z\sim\cN(0, I_r)$ is independent of them, and $\mu_{a,N}(\pi)\in\R^d$, $A_{a,N}(\pi)\in\R^{d\times d}$, $B_{a,N}(\pi)\in\R^{d\times r}$ may be singular. 
For $X\sim\cN(0, I_d)$, define the one-site conditional probabilities
\begin{align}\label{eq:p-def}
p_{a,N}(\pi,z)
=\P\bigl(\mu_{a,N}(\pi)+A_{a,N}(\pi)X+B_{a,N}(\pi)z\in\cK_{a,N}(\pi)\bigr),
\qquad z\in\R^r,
\end{align}
and set
\begin{align}\label{eq:LambdaN}
\Lambda_N(\pi,z)=\sum_{a \in \mc{A}}\pi_a\log p_{a,N}(\pi,z),
\end{align}
with the conventions $0\log0=0$ and $c\log0=-\infty$ for $c>0$.

Our general theorem gives the exponential asymptotics for  $\cN_N$ in \eqref{eq:count} with Gaussian fields structured as in~\eqref{eq:CF} under some natural conditions.

\begin{thm}\label{thm:general}
    For each $\gs \in \cA^N$, let $\{Y_i(\gs)\}_{i=1}^N$ be Gaussian random fields satisfying $\eqref{eq:CF}$. Let $\gL:\gD_k \times \bR^r \to [-\infty, 0]$ satisfy the following conditions:
    \begin{enumerate}
        \item[(i)] For every sequence $\pi_N \in \cP_N$ and $z_N \in \bR^r$ with $(\pi_N, z_N) \to (\pi, z)$, 
        \begin{align}\label{cond:gene-upper}
        \limsup_{N} \gL_N(\pi_N, z_N) \le \gL(\pi, z).
        \end{align}
        \item[(ii)] For every $(\pi,z)$ with $\Lambda(\pi,z)>-\infty$ and every
$\eps>0$, there are $\delta>0$ and profiles
$\pi_N\in\cP_N$ with $\pi_N\to\pi$ such that, for all sufficiently
large $N$, 
\begin{align}\label{cond:gene-lowers}
    \inf_{\|w-z\|_2\le\delta}\gL_N(\pi_N,w)
 \ge\gL(\pi,z)-\eps.
\end{align}
    \end{enumerate}
Then we have 
\[
\lim_{N\to\infty}
\frac{1}{N}\log\mathbb{E}\mathcal{N}_N
=
\sup_{\pi\in\Delta_k}
\sup_{z\in\mathbb{R}^r}
\left\{
\Lambda(\pi,z)
-\frac{\|z\|_2^2}{2}
+\mathcal{H}(\pi)
\right\}.
\]
\end{thm}

The above Theorem~\ref{thm:general} allows vanishing conditional probabilities,
which is needed for the GS model in Section~\ref{sec:GS}, where we verify
\eqref{cond:gene-upper} and \eqref{cond:gene-lowers} directly.
However, in many other applications in Sections~\ref{sec:p-spin}
and~\ref{sec:disorder-persistence}, the following corollary gives a more convenient sufficient condition to check.

\begin{cor}
\label{thm:main}
For each $\gs \in \cA^N$, let $\{Y_i(\gs)\}_{i=1}^N$ be Gaussian random fields satisfying $\eqref{eq:CF}$. Suppose that there exists a finite continuous function
\(
\Lambda:\Delta_k\times\mathbb{R}^r\longrightarrow\mathbb{R}
\)
such that, for every $M>0$, $\Lambda_N$ is finite on
$\mathcal{P}_N\times\{z:\|z\|_2\le M\}$ for all sufficiently large
$N$, and
\begin{equation}
\lim_{N\to\infty}
\sup_{\substack{\pi\in\mathcal{P}_N\\ \|z\|_2\le M}}
\left|\Lambda_N(\pi,z)-\Lambda(\pi,z)\right|
=0.
\label{eqn:main-assumption}
\end{equation}
Then
\[
\lim_{N\to\infty}
\frac{1}{N}\log\mathbb{E}\mathcal{N}_N
=
\sup_{\pi\in\Delta_k}
\sup_{z\in\mathbb{R}^r}
\left\{
\Lambda(\pi,z)
-\frac{\|z\|_2^2}{2}
+\mathcal{H}(\pi)
\right\}.
\]
\end{cor}

\begin{rem}
In all applications in this paper, $\Lambda$ will have the form
\begin{align}
\label{eqn:lambda-as-ps}
    \Lambda(\pi,z)
    =
    \sum_{a\in\mc{A}}\pi_a\log p_a(\pi,z)
\end{align}
for some limiting probabilities $p_a$, which may vanish.
For $M>0$, whenever $p_{a,N}(\pi,z)$ and $p_a(\pi,z)$ are
positive for all $a\in\mc{A}$, $\pi\in\mathcal{P}_N$, and
$\|z\|_2\le M$, we have
\begin{align}
\label{eqn:lambda-vs-p}
    \sup_{\substack{\pi\in\mathcal{P}_N\\ \|z\|_2\le M}}
    \abs{\Lambda_N(\pi,z)-\Lambda(\pi,z)}
    \le
    \sup_{\substack{\pi\in\mathcal{P}_N\\ \|z\|_2\le M}}
    \sup_{a\in\mc{A}}
    \abs{\log p_{a,N}(\pi,z)-\log p_a(\pi,z)}.
\end{align}
In Sections~\ref{sec:p-spin} and~\ref{sec:disorder-persistence},
these probabilities are positive for all sufficiently large $N$,
and we verify \eqref{eqn:main-assumption} by showing that the
right-hand side of \eqref{eqn:lambda-vs-p} tends to zero for every
$M>0$. Condition~\eqref{eqn:main-assumption} itself is more general,
since it requires convergence only of the weighted sums defining
$\Lambda_N$ and $\Lambda$.
\end{rem}

\subsection{Proof of Theorem~\ref{thm:general}}
We first introduce the following lemma, which reduces the problem of analyzing exponential asymptotics of Gaussians.
\begin{lem}
\label{lem:general-nn-expansion}
For every $N\ge 1$, we have 
\begin{equation}
\label{eqn:general-nn-expansion}
        \E\cN_N
        =\sum_{\pi\in\cP_N}\frac{N!}{\prod_{a \in \mc{A}}(N\pi_a)!}\,
        \E_Z\exp\Bigl\{N\Lambda_N\Bigl(\pi,\tfrac Z{\sqrt N}\Bigr)\Bigr\}.
\end{equation}
\end{lem}

\begin{proof}
For fixed $\sigma \in \mc{A}^N$, using that the $X_i$'s are independent of $Z$ (and of each other), we find that
\begin{equation}
\label{eqn:basic-factoring}
\begin{split}
    &\P(Y_i(\sigma) \in \mc{K}_{\sigma_i,N}(\pi(\sigma)) \text{ for all i}) \\
    &= \E_Z\left[\E_{X_1,\ldots,X_N}\left[\prod_{i=1}^N \mathds{1}_{\mu_{\sigma_i,N}(\pi(\sigma)) + A_{\sigma_i,N}(\pi(\sigma))X_i + \frac{1}{\sqrt{N}} B_{\sigma_i,N}(\pi(\sigma))Z \in \mc{K}_{\sigma_i,N}(\pi(\sigma))} \right] \right] \\
    &= \E_Z\left[ \prod_{i=1}^N p_{\sigma_i,N}\left(\pi(\sigma),\frac{Z}{\sqrt{N}} \right) \right] \\
    &= \E_Z\left[ \prod_{a \in \mc{A}} p_{a,N}\left(\pi(\sigma),\frac{Z}{\sqrt{N}} \right)^{N\pi_a(\sigma)} \right] \\
    &= \E_Z\exp\left\{N\Lambda_N\left(\pi(\sigma),\frac{Z}{\sqrt{N}}\right)\right\}.
\end{split}
\end{equation}
Thus
\[
    \E \mc{N}_N = \sum_{\sigma \in \mc{A}^N} \E_Z \exp\left\{ N \Lambda_N\left(\pi(\sigma),\frac{Z}{\sqrt{N}} \right) \right\}.
\]
Since the summands only depend on $\sigma$ through its profile $\pi(\sigma) \in \mc{P}_N$, we can write this as a sum over $\pi \in \mc{P}_N$, using that the number of $\sigma$'s in $\mc{A}^N$ such that $\pi(\sigma)$ is some given $\pi$ is $\frac{N!}{\prod_{a \in \mc{A}} (N\pi a)!}$, which completes the proof.
\end{proof}

\begin{proof}[Proof of Theorem \ref{thm:general}]

Write $Z_N=Z/\sqrt N\sim \mc{N}(0,N^{-1} I_r)$. Let $M_1,M_2,\ldots$ be i.i.d.\ uniform on
$\mathcal A$, independently of $Z_N$, and let
$\widetilde\Pi_N=N^{-1}\sum_{i=1}^N\delta_{M_i}$ be their empirical
measure. Let
\(
\Pi_N=(\widetilde\Pi_N(\{a\}))_{a\in\mathcal A}\in\mathcal P_N,
\)
since \(
\mathbb P(\Pi_N=\pi)=k^{-N}N!/\prod_{a\in\mathcal A}(N\pi_a)!,
\)
Lemma~\ref{lem:general-nn-expansion} gives
\begin{align}
\mathbb E\mathcal N_N
=k^N\mathbb E\exp\{N\Lambda_N(\Pi_N,Z_N)\}.
\label{eqn:general-joint-representation}
\end{align}

By Sanov's theorem \cite[Theorem~2.1.10]{DemZei1998}, the empirical
profile $\Pi_N$ satisfies a large-deviation principle (LDP) at speed
$N$ with good rate function $\log k-\mathcal H(\pi)$. The scaled
Gaussian $Z_N=Z/\sqrt N$ satisfies an LDP at the same speed with good
rate function $\|z\|_2^2/2$. Since $\Pi_N$ and $Z_N$ are independent, by \cite[Exercise~4.2.7]{DemZei1998}, 
their joint law satisfies an LDP with good rate function
\[
I(\pi,z)=\log k-\mathcal H(\pi)+\frac{\|z\|_2^2}{2}.
\]
 The laws of $(\Pi_N,Z_N)$ are
exponentially tight, since $\Delta_k$ is compact and $Z_N$ has
Gaussian tails. Set
\[S:=\log k+\sup_{(\pi,z)\in\Delta_k\times\mathbb R^r}
\{\Lambda(\pi,z)-I(\pi,z)\}.
\]
For the upper bound, extend $\Lambda_N$   outside
$\mathcal P_N\times\mathbb R^r$ as $-\infty$. Then $-\Lambda_N$ are nonnegative
Borel functions, and condition~\eqref{cond:gene-upper} implies
$\liminf_N[-\Lambda_N(\pi_N,z_N)]\ge-\Lambda(\pi,z)$ whenever
$(\pi_N,z_N)\to(\pi,z)$ in $\Delta_k\times\mathbb R^r$.
The generalized Varadhan upper bound
\cite[Theorem~3.5]{Mariani2018}, applied to $-\Lambda_N$ and
the joint law of $(\Pi_N,Z_N)$, gives
\begin{align}\label{eqn:general-varadhans}
\limsup_{N\to\infty}\frac1N\log\mathbb E\mathcal N_N
\le \log k+\sup_{\pi,z}\{\Lambda(\pi,z)-I(\pi,z)\}=S.
\end{align}

For the lower bound, we first deal with the case $S>-\infty$. Fix $(\pi,z)$ with
$\Lambda(\pi,z)>-\infty$, and let $\varepsilon>0$. Choose $\delta>0$
and profiles $\pi_N\to\pi$ as in \eqref{cond:gene-lowers}, and let
\[B_\delta(z)=\{w\in\mathbb R^r:\|w-z\|_2<\delta\}.\]
For large $N$, restricting \eqref{eqn:general-joint-representation}
to $\Pi_N=\pi_N$ and $Z_N\in B_\delta(z)$ gives
\[
\mathbb E\mathcal N_N
\ge k^N\mathbb P(\Pi_N=\pi_N)
\,e^{N(\Lambda(\pi,z)-\varepsilon)}
\mathbb P(Z_N\in B_\delta(z)).
\]
By the estimate in \cite[Lemma~2.1.9]{DemZei1998}, we have 
$N^{-1}\log\mathbb P(\Pi_N=\pi_N)\to\mathcal H(\pi)-\log k$.
The Gaussian LDP lower bound gives
$\liminf_N N^{-1}\log\mathbb P(Z_N\in B_\delta(z))
\ge-\|z\|_2^2/2$, since $z\in B_\delta(z)$. Consequently,
\[
\liminf_{N\to\infty}\frac1N\log\mathbb E\mathcal N_N
\ge\mathcal H(\pi)+\Lambda(\pi,z)-\frac{\|z\|_2^2}{2}-\varepsilon.
\]
Letting $\varepsilon\downarrow0$ and taking the supremum over all
$(\pi,z)$ with $\Lambda(\pi,z)>-\infty$ yields
\begin{equation}
\liminf_{N\to\infty}\frac1N\log\mathbb E\mathcal N_N\ge S.
\label{eqn:general-exponential-equivalence}
\end{equation}
If $S=-\infty$, the lower bound is automatic. Combining \eqref{eqn:general-varadhans} and
\eqref{eqn:general-exponential-equivalence}, and substituting the
formula for $I$, completes the proof.
\end{proof}


\subsection{Proof of Corollary~\ref{thm:main}}

For $M>0$, let 
\[
\eps_N(M) := \sup_{\pi \in \cP_N, \norm{w}\le M}\abs{\gL_N(\pi,w) - \gL(\pi,w)}
\]
By the assumption~\eqref{eqn:main-assumption}, we have $\eps_N(M) \to 0$. Let $\pi_N \in \cP_N, z_N \in \bR^r$ be such that $(\pi_N, z_N) \to (\pi,z)$. Now set $M> \norm{z}+1$. Then we have 
\begin{align*}
\abs{\gL_N(\pi_N,. z_N) - \gL(\pi,z)} & \le \abs{\gL_N(\pi_N, z_N) - \gL(\pi_N,z_N)} +\abs{\gL(\pi_N, z_N) - \gL(\pi,z)} \\
& \le \eps_N(M) +\abs{\gL(\pi_N, z_N) - \gL(\pi,z)},
\end{align*}
and the last term tends to zero as $N\to \infty$ by the continuity of $\gL$ and the assumption~\eqref{eqn:main-assumption}. As a result, $\lim_N \gL_N(\pi_N, z_N) = \gL(\pi,z)$, which verifies the condition~\eqref{cond:gene-upper} in Theorem~\ref{thm:general}. Note that by definition $\gL_N \le 0$, and thus for $\pi_N \to \pi$, we have 
\[
\gL(\pi, z) = \lim_N \gL_N(\pi_N, z) \le 0.
\]

It remains to verify \eqref{cond:gene-lowers}. To do this, we need a
sequence of empirical profiles approaching $\pi$ and a fixed neighborhood
of $z$ on which the required lower bound holds. Fix $(\pi,z)$ and
$\eps>0$, and choose any $\pi_N\in\cP_N$ with $\pi_N\to\pi$.
By continuity of $\Lambda$ at $(\pi,z)$, we can choose $0<\delta\le1$ such that
\[
 \sup_{\substack{\pi'\in\Delta_{k},\ \|\pi'-\pi\|_2\le\delta\\
                  \|w-z\|_2\le\delta}}
 |\Lambda(\pi',w)-\Lambda(\pi,z)|\le\frac{\eps}{2}.
\]
Set $M=\|z\|_2+1$. Since $\pi_N\to\pi$ and $\eps_N(M)\to0$,
for all sufficiently large $N$,
\[
 \|\pi_N-\pi\|_2\le\delta,
 \qquad \eps_N(M)\le\frac{\eps}{2}.
\]
Combining $\|w-z\|_2\le\delta\Rightarrow\|w\|_2\le M$, we have
\begin{align*}
 \inf_{\|w-z\|_2\le\delta}\Lambda_N(\pi_N,w)
 &\ge \inf_{\|w-z\|_2\le\delta}\Lambda(\pi_N,w)-\eps_N(M)\\
 &\ge \Lambda(\pi,z)-\frac{\eps}{2}-\eps_N(M)\\
 &\ge \Lambda(\pi,z)-\eps.
\end{align*}
This verifies~\eqref{cond:gene-lowers}. The conclusion follows
from Theorem~\ref{thm:general}.


\section{Application I: Mixed \texorpdfstring{$p$}{p}-spin model}
\label{sec:p-spin}

In this section, we apply the general results in Section~\ref{sec:general} to the mixed $p$-spin models with potentials such as external fields. In particular, our main results about this model are the following. Theorem \ref{thm:complexity-p-spin} gives a two-dimensional variational formula for the gapped complexity for general $V$. We further show that in the special case where $V$ is an external field, the variational formula reduces to one dimension and becomes more tractable to analyze. Corollary \ref{cor:p-spin-analysis} includes some of those implications based on the simplified formula. The proofs are included in Section \ref{subsec:p-spin-establishing} and \ref{subsec:p-spin-analyzing}, respectively. 


\subsection{Setting and main results}
For $\gs \in \{-1,+1\}^N$ and a potential function $V:[-1,1] \to\bR$, we define the mixed $p$-spin Hamiltonian 
\begin{align}\label{eq:mixed-H}
H_N(\sigma)=H_N^{\mathrm G}(\sigma)+NV(m_N(\sigma)),
\qquad
H_N^{\mathrm G}(\sigma)=\sum_{p=2}^\infty\frac{\beta_p}{N^{(p-1)/2}}
\sum_{I\in\binom{[N]}p}J_I^{(p)}\sigma_I,
\end{align}
where $\sigma_I\defeq\prod_{i\in I}\sigma_i$, $(\gb_p)_{p\ge 2}$ are real coefficients which are not all zero, and $m_N$ is the magnetization $m_N(\gs) = \frac1N \sum_{i\in[N]} \gs_i$. The disorders $J_I^{(p)}$ are i.i.d standard Gaussian. The potential function can be, for example, an external field $V(m)=hm$ or a Curie-Weiss perturbation $V(m)=hm-\lambda m^2/2$. Note that 
\[
\E [H_N^{\mathrm G}(\gs)H_N^{\mathrm G}(\gt)]= N \xi(R(\gs,\gt)) + O(1), \qquad R(\gs, \gt) = \frac1N \sum_{i\in[N]} \gs_i \gt_i,
\]
where
\begin{equation}\label{eq:mixed-xi}
\xi(q)=\sum_{p=2}^\infty\frac{\beta_p^2}{p!}q^p
\end{equation}
is known as the structure function of the mixed $p$-spin model. We assume that $\xi$ has radius of convergence strictly bigger than one. We will take the convention that ``pure $p$-spin'' refers to the structure function
\[
    \xi(q) = \frac{q^p}{p!}.
\]

For $\gamma\ge0$, let $\cN_N^{\mathrm{mix}}(\xi, V,\gamma)$ be the random number of $\gamma$-gapped local minima, that is, the number of configurations satisfying
\begin{align}\label{eq:mixed-gap}
H_N(\sigma^{(i)})-H_N(\sigma)\ge\gamma
\qquad\text{for every }i\le N,
\end{align}
where $\sigma^{(i)}$ denotes the one-spin flip perturbation 
\[
    \gs_j^{(i)} = \begin{cases} \gs_j,  \quad \text{for $j\neq i$},\\
-\gs_i \quad \text{for $j=i$.}\end{cases}
\]
Our main result is an explicit variational formula for the exponential asymptotics of $\E[\cN_N^{\mathrm mix}(\xi, V, \gamma)]$ for the mixed $p$-spin model. To prepare the statement of main results, we introduce some necessary notations. Let
\[
    a = \xi'(1) > 0, \qquad b = \xi''(1) > 0,
\]
as well as the binary entropy in magnetization coordinates
\begin{align}\label{eq:mixed-entropy}
\fh(m)=-\frac{1+m}2\log\frac{1+m}2-\frac{1-m}2\log\frac{1-m}2.
\end{align}
so that $\fh(m)=\cH\bigl(\tfrac{1+m}2,\tfrac{1-m}2\bigr)$ in the notation of Section~\ref{subsec:framework}. 

Throughout the paper, we will use $\phi$ and $\Phi$ for the density and cumulative distribution function, respectively, of a standard Gaussian variable.

\begin{theorem}
\label{thm:complexity-p-spin}
For every potential function $V \in C^1([-1,1])$ and every $\gamma\ge 0$, the annealed complexity for the mixed $p$-spin models is given by
\begin{equation}
\label{eqn:complexity-p-spin}
\begin{split}
    \Sigma^{\textup{mix}}(\xi,V,\gamma) &\defeq \lim_{N \to \infty} \frac{1}{N} \log \E \mc{N}_N^{\textup{mix}}(\xi,V,\gamma) \\
    &= \sup_{m \in [-1,1]} \sup_{z \geq 0} \Big\{ 
    \frac{1+m}{2} \log\left(\Phi\left(\sqrt{\frac{b}{a}}z - \frac{1}{\sqrt{a}}\left(\frac{\gamma}{2} + V'(m)\right)\right)\right) \\
    &\qquad \qquad \qquad \qquad + \frac{1-m}{2} \log\left(\Phi\left(\sqrt{\frac{b}{a}}z - \frac{1}{\sqrt{a}}\left(\frac{\gamma}{2} - V'(m)\right)\right)\right) \\
    &\qquad \qquad \qquad \qquad + \fh(m) - \frac{z^2}{2} \Big\}.
\end{split}
\end{equation}
In the case $V(m) = hm$ for some $h\in \bR$, the above formula simplifies as 
\begin{equation}
\label{eqn:p-spin-simplified-formula}
\begin{split}
    &\Sigma^{\textup{mix}}(\xi,h,\gamma) \\
    &= \sup_{z \geq 0} \left\{ \log\left(\Phi\left(\sqrt{\frac{b}{a}}z - \frac{1}{\sqrt{a}}\left(\frac{\gamma}{2} + h\right)\right) + \Phi\left(\sqrt{\frac{b}{a}}z - \frac{1}{\sqrt{a}}\left(\frac{\gamma}{2} - h\right)\right)\right) - \frac{z^2}{2} \right\}.
\end{split}
\end{equation}
\end{theorem}

For the special external field potential case, the above variational formula is much simpler, and one can obtain some further implications based on this. (When discussing pure $p$-spin, we will abusively write $p$ instead of $\xi$.)

\begin{cor}
\label{cor:p-spin-analysis}
If $V(m) = hm$ for some $h \in \R$, then we have the following results.
\begin{enumerate}
\item The function $(h,\gamma)\mapsto $ $\Sigma^{\textup{mix}}(\xi,h,\gamma)$ is continuous on $\bR \times [0,\infty)$. \label{item:p-spin-analysis-1}
\item For a fixed mixture $\xi$ and fixed $h$, the function $\Sigma^{\textup{mix}}(\xi,h,\gamma)$ is strictly decreasing for $\gamma \geq 0$, positive when $\gamma = 0$, and negative for sufficiently large $\gamma$. In particular, there exists a unique $\gamma_c(\xi,h) \in (0,\infty)$ satisfying
\begin{equation}
\label{eqn:gamma-c}
    \Sigma^{\textup{mix}}(\xi,h,\gamma_c(\xi,h)) = 0.
\end{equation}
If $\gamma > \gamma_c(\xi,h)$, then 
\[
    \limsup_{N \to \infty} \frac{1}{N} \log \P(\text{there exists a $\gamma$-gapped state}) < 0.
\]
\label{item:p-spin-analysis-2}
\item For a fixed mixture $\xi$, the function $\Sigma^{\textup{mix}}(\xi,h,0)$ is strictly decreasing for $h \geq 0$. Its asymptotic behavior for large $h$ is given by 
\begin{equation}
\label{eqn:large-h-asymptotics}
    \lim_{h \to +\infty} \left(\Sigma^{\textup{mix}}(\xi,h,0) e^{\frac{h^2}{a}} \right) = \frac{b}{\pi a},
\end{equation}
and in particular it is positive for all $h$.
\label{item:p-spin-analysis-3}
\item In the case of pure $p$-spin glasses, the annealed complexity $\Sigma^{\textup{mix}}(p,0,0)$ is non-decreasing in $p$, and 
\begin{equation}
\label{eqn:rem-limit}
    \lim_{p \to \infty} \Sigma^{\textup{mix}}(p,0,0) = \log 2.
\end{equation}
\label{item:p-spin-analysis-4}
\end{enumerate}
\end{cor}

Finally, we make various remarks about interpretations.

\begin{rem}
[On gapped states]
\label{rem:gapped}
As previously mentioned, Sellke \cite{Sel2025} recently gave a short proof that low-temperature Glauber dynamics for the SK model mixes exponentially slowly from worst-case initial conditions, based on the existence of gapped states (proved in this context by Dandi, Gamarnik, and Zdeborov\'a \cite{DandiGamarnikZdeborova2023}). Roughly speaking, if $\sigma^\ast$ is a gapped state, then one can find $\rho > 0$ such that the Hamming sphere of radius $\rho N$ around $\sigma^\ast$ has much smaller Gibbs measure than the corresponding Hamming ball, so Glauber dynamics started at $\sigma^\ast$ is bottlenecked.

For the SK model with an external field, the existence of gapped states would imply, by Sellke's argument essentially verbatim, that low-temperature Glauber dynamics also mixes exponentially slowly. Since our results concern the annealed complexity, they only show that gapped states \emph{do not} exist for $\gamma > \gamma_c(2,h)$. Showing that they do exist for $\gamma < \gamma_c(2,h)$ would require a second-moment computation, which is beyond the scope of this paper. Nevertheless, the fact that $\gamma_c(2,h) > 0$ suggests that this conclusion could be correct.
\end{rem}

\begin{table}[ht]
\centering
\caption{Numerical approximations for pure $p$-spin models in zero external field: the critical gap $\gamma_c(p,0)$ from \eqref{eqn:gamma-c}, and the annealed complexity of local minima $\Sigma^{\textup{mix}}(p,0,0)$.}
\label{tab:pure-p-spin}
\medskip
\begin{tabular}{l|ccccc}
\hline
$p$ & $2$ & $3$ & $4$ & $5$ & $6$ \\
\hline
$\gamma_c(p,0)$ & $0.7026$ & $0.9080$ & $0.7325$ & $0.4604$ & $0.2446$ \\
$\Sigma^{\textup{mix}}(p,0,0)$ & $0.1992$ & $0.2956$ & $0.3552$ & $0.3965$ & $0.4273$ \\
\hline
\end{tabular}
\end{table}

\begin{rem}[On the threshold $\gamma_c$]\label{rem:gamma-c}
Numerical values of $\gamma_c(p,0)$ for $p = 2, \ldots, 6$ are listed in Table~\ref{tab:pure-p-spin}. The value $\gamma_c(2,0) \approx 0.7026$ is consistent with recent work of Dandi, Gamarnik, and Zdeborov\'a \cite{DandiGamarnikZdeborova2023}, who found a numerically matching critical threshold for a slightly different problem. Informally, they study ``approximately $\gamma$-gapped bisections'': configurations with approximately equal numbers of plus and minus spins (that is, bisections of $\{1, \ldots, N\}$ into a ``plus'' half and a ``minus'' half) which satisfy \eqref{eq:mixed-gap} for all but a small number of $i$. Their condition is thus stricter than ours in one direction (only certain $\sigma$'s are eligible) and weaker in another (the gap only needs to hold in most directions). Their Theorem 3 shows that an implicitly defined threshold, numerically $0.7026$ after adjusting the normalization\footnote{Their normalization differs from ours by a factor of two (compare their eq.~(1.3) with our \eqref{eq:mixed-gap}), so they write their critical threshold as approximately $0.3513 = 0.7026/2$.}, is critical for the existence of such configurations.

The bisection restriction is easy to explain. In \eqref{eqn:complexity-p-spin}, the parameter $m$ plays the role of the magnetization, and in deriving \eqref{eqn:p-spin-simplified-formula} we show that the optimal $m$ is zero when the external field vanishes. Thus the largest contribution to our count of gapped states comes from configurations with magnetization $o(1)$, that is, from approximate bisections, although \cite{DandiGamarnikZdeborova2023} requires a much more precise balance between plus and minus spins. When the external field is nonzero, the optimal $m$ is nonzero, and typical gapped states are no longer bisections.

The results of \cite{DandiGamarnikZdeborova2023} also hold beyond the SK model, for example for non-Gaussian couplings, although only for $p = 2$. Moreover, their second-moment method gives both sides of the threshold, whereas our first-moment computation gives only the non-existence side (see Remark~\ref{rem:gapped}).
\end{rem}

\begin{rem}[Recovering known results for pure models]\label{rem:recover-SK}
For the SK model, $\xi(q) = q^2/2$, we have $a = b = 1$; writing $\Sigma^{\textup{SK}}(h,\gamma)$ for its complexity, \eqref{eqn:p-spin-simplified-formula} becomes
\[
    \Sigma^{\textup{SK}}(h,\gamma) = \sup_{z \geq 0} \left\{ \log\left( \Phi\left( z - \left(\frac{\gamma}{2}+h\right) \right) + \Phi\left( z - \left(\frac{\gamma}{2}-h\right)\right)\right) - \frac{z^2}{2}\right\}.
\]
The special case $h = \gamma = 0$ gives
\[
    \Sigma^{\textup{SK}}(0,0) = \sup_{z \geq 0} \left\{ \log(2\Phi(z)) - \frac{z^2}{2}\right\} \approx 0.1992,
\]
a classical result (this exact formula appears as eqs.~(4.2)--(4.4) in \cite{TanEdw1980}), whose numerical value also matches the recent rigorous result of Addario-Berry, Devroye, Lugosi, and Oliveira \cite{AddarioBerryEtAl2019}. More generally, the formula for pure $p$-spin models with $h = \gamma = 0$, given in \eqref{eqn:pure-p-spin-formula} below (see also Table~\ref{tab:pure-p-spin}), has already appeared in the physics literature, in work of Gross and M\'ezard (eqs.~(53)--(55) in \cite{GroMez1984}) and of Rieger (eqs.~(28)--(29) in \cite{Rie1992}).
\end{rem}

\begin{rem}[Lack of topological trivialization]
\label{rem:topological-trivialization}
Many continuous Gaussian fields, in particular spherical spin glasses, are known to undergo a ``topological trivialization'': if one adds $h$ times a ``simple'' random function, with few critical points, to a ``complicated'' one, with many critical points, then there is a finite $h_c$ such that the sum is complicated for $h < h_c$ and simple for $h > h_c$ (the ``topologically trivial'' phase). Usually, ``complicated'' and ``simple'' mean positive and zero complexity. For example, if $\Sigma(h)$ denotes the complexity of critical points of a mixed spherical spin glass with an external field of size $h$, then Belius, \v{C}ern{\'y}, Nakajima, and Schmidt \cite{BelCerNakSch2022} found a mixture-dependent $h_c$ such that $\Sigma(h)$ is positive for $h < h_c$ and zero for $h \geq h_c$. (In fact, they showed that for $h > h_c$ there are, with high probability, only two critical points on the sphere, the minimum permitted by Morse theory.) Huang and Sellke later extended this to multi-species spherical models, and also showed a ``strong topological trivialization,'' in the sense of absence of even approximate critical points \cite{HuaSel2026}. For some of these models, the same $h_c$ is also critical for replica symmetry breaking \cite{CheSen2017} and for algorithmic hardness \cite{HuaSel2024,HuaSel2023}.

Surprisingly, item \eqref{item:p-spin-analysis-3} shows that $\Sigma^{\textup{mix}}(\xi,h,0)$ is positive for every $h$ and every $\xi$: no Ising mixture exhibits topological trivialization in the external field, at least at the annealed level. For the SK model, this was already predicted by Dean \cite{Dea1994}. Since the quenched complexity $\lim_{N \to \infty} \frac{1}{N} \E[\log \mc{N}_N^{\textup{mix}}(\xi,h,0)]$ is bounded above by the annealed one, it could still exhibit a topological trivialization in $h$. However, quenched complexities are poorly understood even in the continuous case: they can be computed rigorously when they coincide with the annealed ones \cite{Sub2017,Kiv2023}, but to the best of our knowledge there is no rigorous computation of the quenched complexity of any spherical spin glass for which it differs from the annealed one.
\end{rem}

\begin{rem}[Connections to the random energy model]
Derrida's random energy model (REM) \cite{Der1980} has long been understood as the $p \to \infty$ limit of Ising $p$-spin models. Equation \eqref{eqn:rem-limit} shows that this also holds at the level of configurational entropy: its left-hand side is the configurational entropy of the Ising $p$-spin model, and its right-hand side is that of the REM. Indeed, the REM Hamiltonian assigns an i.i.d.\ standard Gaussian to each $\sigma \in \Omega_N = \{-1,+1\}^N$, so a given configuration is a local minimum exactly when its value is the smallest among $N+1$ i.i.d.\ standard Gaussians, which has probability $1/(N+1)$ by symmetry. Thus the expected number of local minima of the REM is $2^N/(N+1)$, whose normalized logarithm tends to $\log 2$. Observations of this type were previously made in the physics literature, for example by Gross and M\'ezard \cite{GroMez1984} and Rieger \cite{Rie1992}.
\end{rem}


\subsection{Establishing the variational formula}
\label{subsec:p-spin-establishing}
We now apply the general results in Section~\ref{sec:general} to prove Theorem \ref{thm:complexity-p-spin}.
\begin{proof}[Proof of Theorem \ref{thm:complexity-p-spin}]
We first define the local fields
\begin{align}\label{eq:mixed-G}
G_i(\sigma)=-\sum_{p=2}^{\infty}\frac{\beta_p}{N^{(p-1)/2}}
\sum_{\substack{I\in\binom{[N]}p\\ i\in I}}J_I^{(p)}\sigma_I.
\end{align}
In the mixed $p$-spin Hamiltonian, every monomial containing $i$ changes sign when the $i$-th spin is flipped, thus $H_N^{\mathrm G}(\sigma^{(i)})-H_N^{\mathrm G}(\sigma)=2G_i(\sigma)$, and with $s=\sigma_i$, $m=m_N(\sigma)$,
\begin{align}\label{eq:mixed-increment}
H_N(\sigma^{(i)})-H_N(\sigma)
=2G_i(\sigma)+N\Bigl[V\Bigl(m-\frac{2s}N\Bigr)-V(m)\Bigr].
\end{align}
Therefore, $\sigma$ is $\gamma$-gapped if and only if $G_i(\sigma)\ge\theta_{\sigma_i,N}(m)$ for every $i\in [N]$, where
\begin{align}\label{eq:mixed-theta}
\theta_{s,N}(m)
\defeq\frac\gamma2-\frac N2\Bigl[V\Bigl(m-\frac{2s}N\Bigr)-V(m)\Bigr].
\end{align}
Since we assumed that $V$ is $C^1$ on $[-1,1]$, this converges pointwise to 
\begin{equation}
\label{eq:theta-limit}
    \theta_s(m) = \frac{\gamma}{2} + sV'(m),
\end{equation}
and actually
\begin{equation}
\label{eqn:theta-taylor}
    \sup_{m \in [-1,1]} \sup_{s \in \mc{A}} \abs{\theta_{s,N}(m) - \theta_s(m)} = o(1),
\end{equation}
by mean value theorem and continuity of $V'$. (Notice that $m_N(\sigma)-\frac{2\sigma_i}{N} = m_N(\sigma^{(i)})$ lies in $[-1,1]$ for each $\sigma$, so that we only need $V$ to be $C^1$ on $[-1,1]$ and not on some $[-1-\epsilon,1+\epsilon]$.) For fixed $\sigma$, counting shared hyperedges gives
\begin{align}\label{eq:mixed-cov}
\Cov(G(\sigma))=a_N I_N+b_N\one\one^{\mathsf T},
\qquad
a_N=\sum_{p=2}^{\infty}\frac{\beta_p^2}{N^{p-1}}\binom{N-2}{p-1},
\quad
b_N=\sum_{p=2}^{\infty}\frac{\beta_p^2}{N^{p-1}}\binom{N-2}{p-2},
\end{align}
independently of $\sigma$, and therefore
\begin{align}\label{eq:mixed-rep}
G_i(\sigma)\equald\sqrt{a_N}\,X_i+\sqrt{b_N}\,U,
\qquad
a_N\To a=\xi'(1),
\quad
Nb_N\To b=\xi''(1),
\end{align}
with $X_1,\dots,X_N,U$ independent standard Gaussians. In the notation of Section~\ref{subsec:framework}, we have
\[
\cA=\{-1,+1\},\quad d=r=1,\quad Y_i=G_i,\quad
\mu_{s,N}=0,\quad A_{s,N}=\sqrt{a_N},\quad B_{s,N}=\sqrt{Nb_N},
\]
\[
\cK_{s,N}(\pi)=[\theta_{s,N}(m(\pi)),\infty),
\qquad m(\pi)=\pi_+-\pi_-,
\]
and thus the Gaussian representation \eqref{eq:CF} holds. In this context, we have
\begin{align*}
    p_{s,N}(\pi,z) = \P(\sqrt{a_N} X + \sqrt{Nb_N} z \geq \theta_{s,N}(m(\pi))) &= \Phi\left( \sqrt{\frac{Nb_N}{a_N}} z - \frac{1}{\sqrt{a_N}} \theta_{s,N}(m(\pi)) \right),
\end{align*}
which converges pointwise for fixed $\pi$ and fixed $z$ to
\[
    p_s(\pi,z) = \Phi\left( \sqrt{\frac{b}{a}} z - \frac{1}{\sqrt{a}} \theta_s(m(\pi)) \right).
\]
To apply Corollary \ref{thm:main}, as suggested in \eqref{eqn:lambda-as-ps}, we define $\Lambda(\pi,z) = \sum_s \pi_s \log p_s(\pi,z)$. To check the main assumption \eqref{eqn:main-assumption}, it suffices to show that the convergence of $p_{s,N}$ to $p_s$ is uniform in the sense of \eqref{eqn:lambda-vs-p}.

Fix $M>0$, and let $N$ be large enough that $Nb_N \le 2b, a_N \ge \frac{a}{2}$, and the l.h.s of \eqref{eqn:theta-taylor} is at most one. Then uniformly over $\abs{z}\le M,\pi \in \gD_k$ and $s \in \cA$, we have 
\begin{equation}
\label{eqn:mixed-p-spin-l-estimate}
\abs{\sqrt{\frac{Nb_N}{a_N}} z - \frac{1}{\sqrt{a_N}} \theta_{s,N}(m(\pi))} \leq L, 
   \qquad  \abs{\sqrt{\frac{b}{a}} z - \frac{1}{\sqrt{a}} \theta_{s}(m(\pi))} \leq L,
\end{equation}
where $L:=2M \sqrt{\frac{b}{a}} + \frac{2}{\sqrt{a}} (\frac{\gamma}{2} + \|V'\|_\infty + 1)$. 
Besides that, 
\begin{align*}
    &\abs{ \left( \sqrt{\frac{Nb_N}{a_N}} z - \frac{1}{\sqrt{a_N}} \theta_{s,N}(m(\pi)) \right) - \left(\sqrt{\frac{b}{a}} z - \frac{1}{\sqrt{a}} \theta_{s}(m(\pi)) \right) } \\
    &\leq M \abs{ \sqrt{\frac{Nb_N}{a_N}} - \sqrt{\frac{b}{a}}} + \abs{\frac{1}{\sqrt{a_N}} - \frac{1}{\sqrt{a}}} \left( \frac{\gamma}{2} + \|V'\|_\infty + 1 \right) + o(1).
\end{align*}
Denoting the r.h.s bound as $\epsilon_N$, one can check that it tends to zero as $N\to \infty$, and we note again that it is uniform in $\abs{z} \leq M$ and $\pi \in \Delta_k$ and $s \in \mc{A}$. Thus, if we let $g(s) = \log(\Phi(s))$, we have
\begin{align*}
    &\sup_{\abs{z} \leq M} \sup_{\pi \in \Delta_k} \sup_{s \in \mc{A}} \abs{\log p_{s,N}(\pi,z) - \log p_s(\pi,z)} \\
    &\leq \sup\left\{ \abs{g(x) - g(y)} : \abs{x}, \abs{y} \leq L \text{ and } \abs{x - y} \leq \epsilon_N \right\} \\
    &\leq \epsilon_N \sup\{\abs{g'(s)} : \abs{s} \leq L\},
\end{align*}
which converges to zero, thus verifying \eqref{eqn:main-assumption}. The functions $p_s(\cdot,\cdot)$ are continuous in $\pi$ and $z$ and positive everywhere, so Corollary \ref{thm:main} applies, and gives
\[
    \Sigma^{\textup{mix}}(\xi,V,\gamma) = \sup_{\pi \in \Delta_2} \sup_{z \in \R} \left\{ \pi_{+1} \log p_{+1}(\pi,z) + \pi_{-1} \log p_{-1}(\pi,z) - \frac{z^2}{2} + \mc{H}(\pi)\right\}.
\]
To get from here to \eqref{eqn:complexity-p-spin}, we just change variables by parametrizing $\pi \in \Delta_2$ by $\pi_+ = \frac{1+m}{2}$ and $\pi_- = \frac{1-m}{2}$ for $m \in [-1,1]$. We can replace $\sup_{z \in \R}$ with $\sup_{z \geq 0}$ because the function being optimized is, for each $z > 0$, larger at $z$ than at $-z$.

When $V(m) = hm$, we have $V'(m) = h$, the variational problem in $m$ for fixed $z$ has the form
\[
    \sup_{m \in [-1,1]} \left\{\frac{1+m}{2} \alpha + \frac{1-m}{2} \beta + \fh(m)\right\},
\]
for some $\alpha, \beta < 0$, and one can easily check that the maximum is obtained at the interior point $m_\ast = \tanh(\frac{\alpha-\beta}{2})$, at which point the value of the problem becomes $\log(e^{\alpha}+e^{\beta})$. This gives \eqref{eqn:p-spin-simplified-formula}.
\end{proof}


\subsection{Analyzing the variational formula}
\label{subsec:p-spin-analyzing}

\begin{proof}[Proof of Corollary \ref{cor:p-spin-analysis}]
 We define
\begin{equation}
\label{eqn:f-definition}
\begin{split}
    f_{\xi,h,\gamma}(z) &= \Phi\left(\sqrt{\frac{b}{a}} z - \frac{1}{\sqrt{a}} \left(\frac{\gamma}{2} + h\right) \right) + \Phi\left(\sqrt{\frac{b}{a}} z - \frac{1}{\sqrt{a}} \left(\frac{\gamma}{2} - h\right) \right), \\
    g_{\xi,h,\gamma}(z) &= \log(f_{\xi,h,\gamma}(z)), \\
    k_{\xi,h,\gamma}(z) &= g_{\xi,h,\gamma}(z) - \frac{z^2}{2}.
\end{split}
\end{equation}
The function $k$ is clearly continuous in $\gamma$, $h$, and $z$, and we have $\lim_{z \to \infty} k_{\xi,h,\gamma}(z) = -\infty$, locally uniformly in $h$ and $\gamma$ for fixed $\xi$. Thus $\Sigma$ inherits continuity in $h$ and $\gamma$ from $k_{\xi,h,\gamma}$, proving part (1).

Next, we fix $\xi$ and $h$ and note that $k_{\xi,h,\gamma}(z)$ is strictly decreasing in $\gamma$ for each $z$. A priori this only tells us that $\Sigma^{\textup{mix}}(\xi,h,\cdot)$ is nonincreasing in $\gamma$, but again since  $\lim_{z \to +\infty} k_{\xi,h,\gamma}(z) = -\infty$, the defining supremum of $\Sigma^{\textup{mix}}(\xi,h,\gamma)$ must be achieved at a finite $z$ for each $\gamma$, and thus $\Sigma^{\textup{mix}}(\xi,h,\gamma)$ is actually strictly decreasing in $\gamma$. Since
\begin{align*}
    k_{\xi,h,0}(0) &= \log\left( \Phi\left( - \frac{h}{\sqrt{a}} \right) + \Phi\left( \frac{h}{\sqrt{a}} \right)\right) = 0, \\
    k'_{\xi,h,0}(0) &= \sqrt{\frac{b}{a}} \left( \phi\left( -\frac{h}{\sqrt{a}} \right) + \phi\left( \frac{h}{\sqrt{a}} \right) \right) > 0,
\end{align*}
we conclude that $k_{\xi,h,0}$ is strictly positive in a small neighborhood to the right of zero, and in particular $\Sigma^{\textup{mix}}(\xi,h,0) > 0$. Since $k_{\xi,h,\gamma}(z) \leq \log(2) - z^2/2$, we can show that $\Sigma^{\textup{mix}}(\xi,h,\gamma)$ is negative just by looking at $k_{\xi,h,\gamma}(z)$ just for $z$ in some interval containing $[0,\sqrt{2\log 2}]$ plus a neighborhood to the right, say $[0,2]$ for concreteness. Since $g_{\xi,h,\gamma}(z)$ is increasing in $z$, on this restricted domain we have 
\[
    k_{\xi,h,\gamma}(z) \leq g_{\xi,h,\gamma}(2) \quad \text{for $z \in [0,2]$}
\]
which tends to negative infinity as $\gamma \to +\infty$. This shows that $\Sigma^{\textup{mix}}(\xi,h,\gamma)$ is negative for sufficiently large $\gamma$. The last claim in part (2) just follows from Markov's inequality, since the probability that there exists a $\gamma$-gapped state is bounded above by the expected number of them, which is exponentially small in this regime. Now we show the monotone decreasing property in $h$. For $z>0$ and $h>0$,
\[
\partial_h f_{\xi,h,0}(z)
=
\frac1{\sqrt a}
\left[
\phi\left(\sqrt{\frac ba}\,z+\frac h{\sqrt a}\right)
-
\phi\left(\sqrt{\frac ba}\,z-\frac h{\sqrt a}\right)
\right]<0.
\]
Thus $k_{\xi,h,0}(z)$ is strictly decreasing in $h\ge0$
for every $z>0$.
The supremum defining $\Sigma^{\mathrm{mix}}(\xi,h,0)$ is attained. 
Every maximizer $z_0$ is positive, since
$k_{\xi,h,0}(z_0)=\Sigma^{\mathrm{mix}}(\xi,h,0)>0=k_{\xi,h,0}(0)$.
For $0\le h_1<h_2$, choose a maximizer $z_2>0$ at $h_2$.
Then
\[
\Sigma^{\mathrm{mix}}(\xi,h_1,0)
\ge k_{\xi,h_1,0}(z_2)
>k_{\xi,h_2,0}(z_2)
=\Sigma^{\mathrm{mix}}(\xi,h_2,0).
\]

Before proving \eqref{eqn:large-h-asymptotics}, we mention briefly an informal heuristic for this result. Let $z(h) = z(\xi,h)$ be the optimizer in the definition of $\Sigma^{\textup{mix}}(\xi,h,0)$, and take as an Ansatz that $\lim_{h \to \infty} z(h) = 0$ (we will show this during the proof). Since $g_{\xi,h,0}(0) = 0$, one can linearize it as $g_{\xi,h,0}(z) \approx z g'_{\xi,h,0}(0)$; this suggests that
\[
    \Sigma^{\textup{mix}}(\xi,h,0) \approx \sup_{z \geq 0} \left\{ z g'_{\xi,h,0}(0) - \frac{z^2}{2} \right\} = \frac{g'_{\xi,h,0}(0)^2}{2}.
\]
Computing $g'_{\xi,h,0}(0)^2 = \frac{2b}{\pi a} e^{-h^2/a}$ suggests \eqref{eqn:large-h-asymptotics}.

Now we present the formal proof. For the lower bound in \eqref{eqn:large-h-asymptotics}, we can pick a test $z$; as suggested by the heuristic above, we choose
\[
    z_\ast = z_\ast(h) = g'_{\xi,h,0}(0) = \sqrt{\frac{2b}{\pi a}} e^{-\frac{h^2}{2a}}.
\]
Thus
\begin{align*}
    \Sigma^{\textup{mix}}(\xi,h,0) &\geq g_{\xi,h,0}(z_\ast) - \frac{z_\ast^2}{2} = g_{\xi,h,0}(z_\ast) - z_\ast g'_{\xi,h,0}(0) + z_\ast g'_{\xi,h,0}(0) - \frac{z_\ast^2}{2} \\
    &= g_{\xi,h,0}(z_\ast) - z_\ast g'_{\xi,h,0}(0) + \frac{z_\ast^2}{2}.
\end{align*}
Since $\lim_{h \to \infty} e^{\frac{h^2}{a}} \frac{z_\ast^2}{2} = \frac{b}{\pi a}$, it suffices to consider the first two terms together; since $e^{\frac{h^2}{a}} = z_\ast^{-2} \cdot \frac{2b}{\pi a}$, it suffices to show that
\begin{equation}
\label{eqn:g-lower-bound}
    \liminf_{h \to \infty} \left( \frac{g_{\xi,h,0}(z_\ast) - z_\ast g'_{\xi,h,0}(z_\ast)}{z_\ast^2} \right) \geq 0.
\end{equation}
From Taylor's theorem, we know that the numerator equals $\frac{z_\ast^2}{2} g''_{\xi,h,0}(s)$ for some $s \in [0,z_\ast]$. Since $z_\ast(h) \to 0$, we can fix $\epsilon > 0$ and assume that $h$ is large enough that such an $s$ lies in $[0,\epsilon]$; then \eqref{eqn:g-lower-bound} follows from Lemma \ref{lem:g-second-order-bounds} below.

For the upper bound, we first show that the true optimizer $z_0(h)$ tends to zero as $h \to \infty$. Recall our argument above that
\[
    \Sigma^{\textup{mix}}(\xi,h,0) = \sup_{0 \leq z \leq \sqrt{2\log 2}} k_{\xi,h,0}(z).
\]
For each $\epsilon > 0$, we have
\[
    \sup_{\epsilon \leq z \leq \sqrt{2\log 2}} k_{\xi,h,0}(z) \leq g_{\xi,h,0}(\sqrt{2\log 2}) - \epsilon^2/2.
\]
As $h \to +\infty$, the right-hand side tends to $-\epsilon^2/2$, but we know that $\Sigma^{\textup{mix}}(\xi,h,0)$ is positive; thus, for each $\epsilon > 0$, the optimizer $z_0(h)$ must be in $[0,\epsilon]$ for all sufficiently large $h$. Now fix such an $\epsilon$. For $z \in [0,\epsilon]$, Taylor's theorem again tells us that there exists $s \in [0,z] \subseteq [0,\epsilon]$ such that
\[
    g_{\xi,h,0}(z) - zg'_{\xi,h,0}(0) = \frac{z^2}{2} g''_{\xi,h,0}(s) \leq \frac{z^2}{2} \sup_{s \in [0,\epsilon]} g''_{\xi,h,0}(s).
\]
From \eqref{eqn:g-second-order-bounds} below, we know that the supremum on the right-hand side tends to zero as $h \to \infty$; thus, for any fixed $\delta > 0$, it is below $\delta$ for large enough $h$. If we take such an $h$ (which is also large enough relative to $\epsilon$), we find that
\begin{align*}
    \Sigma^{\textup{mix}}(\xi,h,0) &= \sup_{0 \leq z \leq \epsilon} \left\{ g_{\xi,h,0}(z) - \frac{z^2}{2} \right\} = \sup_{0 \leq z \leq \epsilon} \left\{ zg'_{\xi,h,0}(0) + \frac{z^2}{2}(-1+\delta) \right\} \\
    &\leq \sup_{z \in \R} \left\{ z g'_{\xi,h,0}(0) + \frac{z^2}{2}(-1+\delta) \right\} = \frac{g'_{\xi,h,0}(0)^2}{2(1-\delta)}.
\end{align*}
Since $g'_{h,0,p}(0)^2 = \frac{2b}{\pi a} e^{-h^2/a}$, this proves that
\[
    \limsup_{h \to +\infty} (\Sigma^{\textup{mix}}(\xi,h,0) \cdot e^{h^2/a}) \leq \frac{1}{1-\delta} \cdot \frac{b}{\pi a},
\]
and we take $\delta \to 0$ to conclude.

Finally, we have
\begin{equation}
\label{eqn:pure-p-spin-formula}
    \Sigma^{\textup{mix}}(p,0,0) = \log 2 + \sup_{z \geq 0} \left\{ \log\left( \Phi\left( z \sqrt{p-1} \right) \right) - \frac{z^2}{2}\right\}.
\end{equation}
The argument of the supremum is nondecreasing in $p$ for each fixed $z$, so $\Sigma^{\textup{mix}}(p,0,0)$ is also nondecreasing in $p$. The argument of the supremum is also nonpositive for each $z$, so the upper bound in \eqref{eqn:rem-limit} is clear; for the lower bound, by taking the test point $z = s/\sqrt{p-1}$ for fixed $s > 0$ we find
\[
    \Sigma^{\textup{mix}}(p,0,0) \geq \log(2) + \log(\Phi(s)) - \frac{s^2}{2(p-1)}.
\]
Taking $p \to +\infty$ and then $s \to +\infty$ completes the proof of \eqref{eqn:rem-limit}.
\end{proof}

\begin{lemma}
\label{lem:g-second-order-bounds}
For every $\epsilon > 0$, we have
\begin{equation}
\label{eqn:g-second-order-bounds}
    \liminf_{h \to \infty} \left( \inf_{s \in [0,\epsilon]} g''_{\xi,h,0}(s) \right) = \limsup_{h \to \infty} \left( \sup_{s \in [0,\epsilon]} g''_{\xi,h,0}(s) \right) = 0.
\end{equation}
\end{lemma}
\begin{proof}
Write $\alpha = \sqrt{\frac{b}{a}}$ and $\beta = \frac{1}{\sqrt{a}}$ for simplicity. Since $\phi'(t) = -t\phi(t)$, we have
\begin{equation}
\label{eqn:g''expansion}
\begin{split}
    \frac{1}{\alpha} g''_{\xi,h,0}(s) &= -\left( \frac{(\alpha s + \beta h) \phi(\alpha s + \beta h) + (\alpha s - \beta h) \phi(\alpha s - \beta h)}{\Phi(\alpha s + \beta h) + \Phi(\alpha s - \beta h)} \right) \\
    &- \left( \frac{\phi(\alpha s + \beta h) + \phi(\alpha s - \beta h)}{\Phi(\alpha s + \beta h) + \Phi(\alpha s - \beta h)} \right)^2.
\end{split}
\end{equation}
The upper bound (i.e., the $\limsup$ bound in \eqref{eqn:g-second-order-bounds}) is easier: For $s \in [0,\epsilon]$, we have $\Phi(\alpha s + \beta h) + \Phi(\alpha s - \beta h) \geq \Phi(\beta h) + \Phi(-\beta h) = 1$, and therefore for all $s \in [0,\epsilon]$ and all $h$ we have
\[
    \frac{1}{\alpha} g''_{h,0,p}(s) \leq \frac{(\beta h - \alpha s) \phi(\alpha s - \beta h)}{\Phi(\alpha s + \beta h) + \Phi(\alpha s - \beta h)} \leq \beta h \phi(\alpha s - \beta h) \leq \beta h \phi(\alpha \epsilon - \beta h),
\]
which tends to zero as $h \to +\infty$. This completes the proof of the upper bound. For the lower bound, one considers the terms in \eqref{eqn:g''expansion} separately: The (negative) first term in \eqref{eqn:g''expansion} can be bounded (above) for sufficiently large $h$ and $s \in [0,\epsilon]$ as 
\begin{align*}
    &\left( \frac{(\alpha s + \beta h) \phi(\alpha s + \beta h) + (\alpha s - \beta h) \phi(\alpha s - \beta h)}{\Phi(\alpha s + \beta h) + \Phi(\alpha s - \beta h)} \right) \\
    &\leq \frac{(\alpha \epsilon + \beta h) \phi(\beta h)}{\Phi(\alpha s + \beta h) + \Phi(\alpha s - \beta h)} \leq \frac{(\alpha \epsilon + \beta h) \phi(\beta h)}{\Phi(\beta h) + \Phi(-\beta h)} = (\alpha \epsilon + \beta h) \phi(\beta h),
\end{align*}
which tends to zero as $h \to +\infty$, and the (negative) second term in \eqref{eqn:g''expansion} can be bounded (above) for sufficiently large $h$ and $s \in [0,\epsilon]$ as
\begin{align*}
    \frac{\phi(\alpha s + \beta h) + \phi(\alpha s - \beta h)}{\Phi(\alpha s + \beta h) + \Phi(\alpha s - \beta h)} \leq \frac{\phi(\beta h) + \phi(\alpha \epsilon - \beta h)}{\Phi(\beta h) + \Phi(-\beta h)} = \phi(\beta h) + \phi(\alpha \epsilon - \beta h),
\end{align*}
which also tends to zero as $h \to +\infty$. This finishes the $\liminf$ bound in \eqref{eqn:g-second-order-bounds}.
\end{proof}


\section{Application II: Disorder persistence of gapped states}
\label{sec:disorder-persistence}


\subsection{Setting and main results.}

Let the mixture function $\xi$ (and its relevant values $a$ and $b$), the mixed spin glass with zero potential $H^G_N$, and the mixed spin glass with potential $H_N$ be as in \eqref{eq:mixed-xi} and \eqref{eq:mixed-H}. Let $\widetilde{H^G_N}$ be an independent copy of $H^G_N$, and for $t \geq 0$ consider the interpolated Hamiltonian
\[
    H_{t,N}(\sigma) = e^{-t} H^G_N(\sigma) + \sqrt{1-e^{-2t}} \cdot \widetilde{H^G_N}(\sigma) + NV(m_N(\sigma)).
\]
Notice that $H_{t,N}$ has the same law as $H_N$, that it equals $H_N$ when $t = 0$, and that it is an independent copy of $H_N$ when (formally) $t = +\infty$. 

In this section, we are interested in high-probability upper bounds on the following random variable.

\begin{defn}
For each $\xi, V, \gamma, \gamma'$, we define the \emph{persistence time}
\begin{align*}
    T_N &= T_N(\xi,V,\gamma,\gamma') \\
    &= \sup\{t \geq 0 : \text{ there exists $\sigma \in \{-1,+1\}^N$ which is } \\
    &\qquad \qquad \qquad \quad \text{$\gamma$-gapped for $H_N$ and $\gamma'$-gapped for $H_{t,N}$}\}.
\end{align*}
\end{defn}

Since $T_N$ is only defined as a supremum, \emph{a priori} it could be poorly named, e.g., it could be that there are small times at which such $\sigma$'s do exist, intermediate times at which they do not, and large times at which they do exist again. The following lemma shows that this is almost-surely not the case when the potential is an external field $V(m) = hm$ and $\gamma \geq \gamma' \geq 2\abs{h}$, and that this is even true $\sigma$-wise. In it, we need the sets
\begin{align*}
    \mc{G}_N(\sigma) &= \{t \geq 0 : \sigma \text{ is $\gamma$-gapped for $H_N$ and $\gamma'$-gapped for $H_{t,N}$}\} \\
    \mc{G}_N &= \{t \geq 0 : \text{ there exists $\sigma \in \{-1,+1\}^N$ which is} \\
    &\qquad \qquad \qquad \text{$\gamma$-gapped for $H_N$ and $\gamma'$-gapped for $H_{t,N}$}\} \\
    &= \bigcup_{\sigma \in \{-1,+1\}^N} \mc{G}_N(\sigma)
\end{align*}
which are random subsets of $[0,\infty)$ such that $T_N = \sup \mc{G}_N$.

\begin{lem}
\label{lem:persistence-times}
Suppose that $V$ is the external field $V(m) = hm$ and that 
\begin{equation}
\label{eqn:gamma-gamma'-2h}
    \gamma \geq \gamma' \geq 2\abs{h}.
\end{equation}
Then the following holds with probability one:
For each $\sigma \in \{-1,+1\}^N$ and each $N$, the set $\mc{G}_N(\sigma)$ is either empty, or equal to $[0,\infty)$, or equal to a closed interval of the form $[0,T]$, and for each $N$, the set $\mc{G}_N$ has the same trichotomy. 
\end{lem}

Thus the term ``persistence time'' is reasonable under these restrictions: Almost surely, such configurations exist for all times up to some point, then cease to exist forever, plus degenerate cases. In other words, there are no so-called ``exceptional times.'' (The restriction that $\gamma \geq \gamma' \geq 2\abs{h}$ might be more than a technical assumption; see Remark \ref{rem:almost-sure-persistence} below.)

To study the random variable $T_N = \sup \mc{G}_N$, we will compute the annealed complexity of one-time marginals, meaning we will study the random variable 
\[
    \mc{N}_N^{\textup{dc}}(t,\xi,V,\gamma,\gamma') = \{\sigma : \sigma \text{ is $\gamma$-gapped for $H_N$ and $\gamma'$-gapped for $H_{t,N}$}\}.
\]

To state our result about the complexity of $\mc{N}_N^{\textup{dc}}$, we will need to extend the notation \eqref{eq:theta-limit} to include $\gamma$ in the argument, namely
\[
    \theta_{s,\gamma}(m) = \frac{\gamma}{2} + sV'(m).
\]
We will also need the following shorthand: For $s \in \{-1,+1\}$, write
\begin{align*}
    &f^{(s)}_{t,a,b,\gamma,\gamma'}(z_1,z_2,m) \\
    &= \Phi_2\left(\sqrt{\frac{b}{a}} z_1 - \frac{1}{\sqrt{a}} \theta_{s,\gamma}(m), \sqrt{\frac{b}{a}} (e^{-t}z_1 + \sqrt{1-e^{-2t}} z_2) - \frac{1}{\sqrt{a}} \theta_{s,\gamma'}(m); e^{-t}\right),
\end{align*}
where $\Phi_2(x,y;\rho)$ is the CDF of the standard bivariate Gaussian with correlation $\rho$ (i.e., of $\mc{N}(0,(\begin{smallmatrix} 1 & \rho \\ \rho & 1 \end{smallmatrix}))$). 

\begin{thm}[Disorder persistence for mixed spin glasses]
\label{thm:disorder-persistence}
 If $V$ is $C^1$ on $[-1,1]$, then
\begin{equation}
\label{eqn:complexity-disorder-persistence}
\begin{split}
    &\Sigma^{\textup{dc}}(t,\xi,V,\gamma,\gamma') \defeq \lim_{N \to \infty} \frac{1}{N} \log \E \mc{N}_N^{\textup{dc}}(t,\xi,V,\gamma,\gamma') \\
    &= \sup_{m \in [-1,1]} \sup_{z_1, z_2 \geq 0} \Big\{ 
    \frac{1+m}{2} \log\left(f^{(+1)}_{t,a,b,\gamma,\gamma'}(z_1,z_2,m)\right) + \frac{1-m}{2} \log\left(f^{(-1)}_{t,a,b,\gamma,\gamma'}(z_1,z_2,m) \right) \\
    &\qquad \qquad \qquad \qquad + \fh(m) - \frac{z_1^2}{2} - \frac{z_2^2}{2} \Big\}.
\end{split}
\end{equation}
Additionally, the quantity $\E \mc{N}^{\textup{dc}}_N(t,\xi,V,\gamma,\gamma')$ is nonincreasing in $t$ for all other parameters held fixed, which implies that $\Sigma^{\textup{dc}}(t,\xi,V,\gamma,\gamma')$ is nonincreasing in $t$ for all other parameters held fixed. 
\end{thm}

Suppose that $V(m) = hm$ for some $h \in \R$. Then $\theta_{s,\gamma}(m) = \gamma/2 + sh$ is actually a constant in $m$, which we will write $\theta_{s,\gamma}$, and we can abusively drop $m$ from the notation $f^{(s)}_{t,a,b,\gamma,\gamma_t}(z_1,z_2)$. In this case we find the following.

\begin{cor}[Analysis for special case of external fields]
\label{cor:disorder-persistence-analysis}
 Suppose that $V(m) = hm$ for some $h \in \R$. Then the formula \eqref{eqn:complexity-disorder-persistence} simplifies as
\begin{equation}
\label{eqn:disorder-persistence-simplified-formula}
\begin{split}
    &\Sigma^{\textup{dc}}(t,\xi,h,\gamma,\gamma') = \sup_{z_1,z_2 \geq 0} \left\{ \log\left( f^{(+1)}_{t,a,b,\gamma,\gamma'}(z_1,z_2) + f^{(-1)}_{t,a,b,\gamma,\gamma'}(z_1,z_2) \right) - \frac{z_1^2}{2} - \frac{z_2^2}{2} \right\}.
\end{split}
\end{equation}
Additionally:
\begin{enumerate}
\item The function $\Sigma^{\textup{dc}}$ is continuous in all of its parameters. \label{item:disorder-persistence-analysis-1}
\item When $t = 0$, the function reduces to the one from Corollary \ref{cor:p-spin-analysis} as 
\begin{equation}
\label{eqn:dc-limit-t-zero}
    \Sigma^{\textup{dc}}(0,\xi,h,\gamma,\gamma') = \Sigma^{\textup{mix}}(\xi,h,\max(\gamma,\gamma')).
\end{equation}
In the special case $h = 0$, the $t \to +\infty$ limit also factors:
\begin{equation}
\label{eqn:dc-limit-t-infinity}
    \lim_{t \to +\infty} \Sigma^{\textup{dc}}(t,\xi,0,\gamma,\gamma') = \Sigma^{\textup{mix}}(\xi,0,\gamma) + \Sigma^{\textup{mix}}(\xi,0,\gamma') - \log(2).
\end{equation}
\label{item:disorder-persistence-analysis-2}
\item Let
\[
    t_c(\xi,h,\gamma,\gamma') = \inf\{t : \Sigma^{\textup{dc}}(t,\xi,h,\gamma,\gamma') < 0\} \in [0,\infty],
\]
and assume that $\gamma \geq \gamma' \geq 2\abs{h}$. Then the (deterministic quantity) $t_c$ is an effective high-probability upper bound for (the random quantity) $T_N$ in the sense that, if $t_c < +\infty$, then for each $\epsilon > 0$, we have
\[
    \limsup_{N \to \infty} \frac{1}{N} \log \P(T_N(\xi,h,\gamma,\gamma') > t_c(\xi,h,\gamma,\gamma') + \epsilon) < 0.
\]
\label{item:disorder-persistence-analysis-3}
\item For local minima of pure $p$-spin in zero external field, we have
\begin{equation}
\label{eqn:tc-p-3-4}
    t_c(p,0,0,0) \begin{cases} \in (0,\infty) & \text{ for } p = 2, 3, \\ = +\infty & \text{ for } p \geq 4. \end{cases}
\end{equation}
\label{item:disorder-persistence-analysis-4}
\end{enumerate}
\end{cor}

\begin{rem}
A priori, the quantity $\Sigma^{\textup{dc}}$ is only counting configurations which are $\gamma$-gapped at time zero and $\gamma'$-gapped at time $t$. However, in the special case $\gamma \geq \gamma' \geq 2\abs{h}$, Lemma \ref{lem:persistence-times} implies that all such points must actually have been $\gamma'$-gapped at all times up to $t$. 
\end{rem}

\begin{rem}
\label{rem:persistence}
Numerically, for the case of pure $p$-spin, one can estimate
\begin{align*}
    t_c(2,0,0,0) \approx 0.3782, \\
    t_c(3,0,0,0) \approx 1.3542.
\end{align*}
It would be interesting to understand whether the dichotomy suggested in \eqref{eqn:tc-p-3-4} is real or just an artifact of the annealed computation. In other words, a weak consequence of our results is that $\P(T_N(p,0,0,0) = \infty) \to 0$ for $p = 2, 3$, i.e., that all local minima die out asymptotically almost surely; but we have not shown that local minima can live forever when $p \geq 4$ (in other words, we have not shown that, e.g., $\liminf_{N \to \infty} \P(T_N(p,0,0,0) = \infty) > 0$ for $p \geq 4$).
\end{rem}

\begin{rem}
\label{rem:persistence-ldp}
We emphasize that the ``persistence'' in our terminology is only about the existence of \emph{at least one} gapped state which survives for order-one times. (Here we are assuming $\gamma = \gamma'$ for simplicity.) If $V \equiv 0$, then for a \emph{fixed} configuration $\sigma$, Theorems \ref{thm:complexity-p-spin} and \ref{thm:disorder-persistence} imply that
\[
    \lim_{N \to \infty} \frac{1}{N} \log \P(\sigma \text{ is $\gamma$-gapped for $H_t$} \, | \, \sigma \text{ is $\gamma$-gapped for $H$}) = \Sigma^{\textup{dc}}(t,\xi,0,\gamma,\gamma) - \Sigma^{\textup{mix}}(\xi,0,\gamma).
\]
For $t > 0$, Corollary \ref{cor:disorder-persistence-analysis} shows that the right-hand side is nonpositive (which was anyway obvious from the left-hand side), but numerically it seems to be the case that it is strictly negative. If so, this would imply that persistence of a given $\gamma$-gapped state is a large deviations event. In other words, there would be disorder chaos at the level of \emph{typical} gapped states, but annealed disorder persistence for \emph{atypical} gapped states.
\end{rem}

\begin{rem}
We prove \eqref{eqn:dc-limit-t-zero} and \eqref{eqn:dc-limit-t-infinity} by analyzing the variational formulas, as a coherence check, but one can also see them immediately at the level of expectations, before taking large-$N$ limits. Indeed, at time zero we have $H_{0,N} = H_N$, so we are just asking for points to be both $\gamma$-gapped and $\gamma'$-gapped. At $t = +\infty$, we are counting $\sigma$'s which are $\gamma$-gapped for $H_N$ and $\gamma'$-gapped for the independent copy $H_{\infty,N}$. The probability $p_N(\sigma,\gamma,\gamma')$ of this is just the product $p_N(\sigma,\gamma) p_N(\sigma,\gamma') = \P(\sigma \text{ is $\gamma$-gapped for } H_N) \P(\sigma \text{ is $\gamma'$-gapped for } H_N)$, and in fact this probability does not depend on $\sigma$, so that
\begin{align*}
    &\E[\# \{\sigma : \sigma \text{ is $\gamma$-gapped for $H_N$ and $\gamma'$-gapped for $H_{\infty,N}$}\}] \\
    &= 2^N \P(\sigma \text{ is $\gamma$-gapped for } H_N) \P(\sigma \text{ is $\gamma'$-gapped for } H_N),
\end{align*}
which makes \eqref{eqn:dc-limit-t-infinity} clear. (When $h \neq 0$, we have $\sum_\sigma p_N(\sigma,\gamma,\gamma') = \sum_\sigma p_N(\sigma,\gamma)p_N(\sigma,\gamma')$, but now the terms depend nontrivially on $\sigma$, so even in the simplest case $\gamma = \gamma'$ it is not clear how to write $\sum_\sigma p_N(\sigma,\gamma)^2$ in terms of $(\sum_\sigma p_N(\sigma,\gamma))^2$.)
\end{rem}


\subsection{Establishing the variational formula.}

We start with a lemma based on the Gaussian interpolation inequality, which requires the following definition: Say that $A \subseteq \R^N$ is ``coordinatewise nondecreasing'' if, whenever $x = (x_1,\ldots,x_N) \in A$ and $y_i \geq x_i$ for all $i$, also $y = (y_1,\ldots,y_N) \in A$.

\begin{lem}
\label{lem:gaussian-interpolation}
Fix $\rho \in [0,1]$. Let $(X,Y_\rho)$ be a centered Gaussian vector in $\R^{N + N}$ with covariance matrix given in block form by
\[
    \begin{pmatrix} C & \rho C \\ \rho C & C \end{pmatrix},
\]
where $C \in \R^{N \times N}$ is a nonsingular, entrywise nonnegative matrix. Then, whenever $A, B \subseteq \R^N$ are coordinatewise nondecreasing, the quantity $\P(X \in A, Y_\rho \in B)$ is nondecreasing in $\rho$. 
\end{lem}
\begin{proof}
Say that a function $f : \R^N \to \R$ is coordinatewise nondecreasing if, whenever $x \leq y$ entrywise, also $f(x) \leq f(y)$. If $f$ is smooth, this is equivalent to $f$ having entrywise nonnegative gradient everywhere. We first claim that, whenever $f$ and $g$ are smooth, coordinatewise-nondecreasing functions, which are bounded and have bounded partial derivatives, we have that $\E[f(X)g(Y_\rho)]$ is nondecreasing in $\rho$. Since $A$ and $B$ are coordinatewise nondecreasing as sets, their indicators are coordinatewise nondecreasing as functions, and can be bounded from below by smooth, coordinatewise-nondecreasing functions which are bounded and have bounded partial derivatives; this will give the claimed result by approximation.

Let $G$ and $H$ be centered Gaussian vectors in $\R^{N + N}$ with respective covariance matrices $\Sigma^G = (\begin{smallmatrix} C & 0 \\ 0 & C \end{smallmatrix})$ and $\Sigma^H =(\begin{smallmatrix} C & C \\ C & C \end{smallmatrix})$, which are defined on the same space and independent. Let $Z(u) = \sqrt{u} G + \sqrt{1-u} H$ be an interpolated Gaussian vector for $u \in [0,1]$. By Gaussian interpolation (see, e.g., \cite[Lemma 7.2.5]{Ver2026} for the exact formula we are using), for any twice-differentiable function $k : \R^{2N} \to \R$, we have
\[
    \frac{\diff}{\diff u} \E[ k(Z(u))] = \frac{1}{2} \sum_{i,j=1}^{2N} (\Sigma^G_{ij} - \Sigma^H_{ij}) \E\left[ \frac{\partial^2 k}{\partial x_i \partial x_j} Z(u) \right]
\]
whenever all expectations exist. Now split $\R^{2N} = \R^{N+N}$ and choose $k$ to have the form $k(x,y) = f(x)g(y)$ for some smooth, coordinatewise-nondecreasing functions $f, g : \R^N \to \R$ which are bounded and have bounded partial derivatives (the latter assumptions ensure that all expectations exist). If we write $G = (G_1,G_2)$ and similarly for $H$ and $Z$, then we have
\[
    \frac{\diff}{\diff u} \E[f(Z_1(u)) g(Z_2(u))] = - \sum_{i,j=1}^N C_{ij} \E\left[ \frac{\partial}{\partial x_i} f(Z_1(u)) \cdot \frac{\partial}{\partial x_j} g(Z_2(u)) \right]. 
\]
Since $C$ is entrywise nonnegative, and since $f$ and $g$ are coordinatewise nondecreasing, the right-hand side is nonpositive as a sum of nonpositive terms. This means that $\E[f(Z_1(u)) g(Z_2(u))]$ is nonincreasing as a function of $u$. Since $Z(u)$ has covariance matrix $(\begin{smallmatrix} C & (1-u) C \\ (1-u) C & C \end{smallmatrix})$, we complete the proof by taking $1-u = \rho$.
\end{proof}

\begin{proof}[Proof of Theorem \ref{thm:disorder-persistence}]
Let $G_i(\sigma)$ be as in \eqref{eq:mixed-G}, and let $\widetilde{G_i}(\sigma)$ be the corresponding local field for $\widetilde{H^G_N}$; then 
\[
    G_{t,i}(\sigma) = e^{-t} G_i(\sigma) + \sqrt{1-e^{-2t}} \widetilde{G_i}(\sigma)
\]
is the local field for the interpolation $H_{t,N}$. Extend the notation \eqref{eq:mixed-theta} by including $\gamma$ in the argument, namely
\[
    \theta_{s,N,\gamma}(m) \defeq \frac{\gamma}{2} - \frac{N}{2}\left[ V \left( m - \frac{2s}{N} \right) - V(m) \right].
\]
Then, as in \eqref{eq:mixed-increment}, $\sigma$ is $\gamma$-gapped for $H_N$ and $\gamma_t$-gapped for $H_{t,N}$ if and only if both $G_i(\sigma) \geq \theta_{\sigma_i,N,\gamma}(m)$ for every $i$, and $G_{t,i}(\sigma) \geq \theta_{\sigma_i,N,\gamma_t}(m)$ for every $i$. 

Write $C_N = a_N I_N + b_N \one \one^{\mathsf T}$ for the $N \times N$ covariance matrix from \eqref{eq:mixed-cov}. Then, for each $\sigma \in \Omega_N$, the Gaussian vector $(G_1(\sigma),\ldots,G_N(\sigma),G_{t,1}(\sigma),\ldots,G_{t,N}(\sigma))$ is centered with covariance matrix in block form
\[
    M = M_N = \begin{pmatrix} C_N & e^{-t} C_N \\ e^{-t} C_N & C_N \end{pmatrix}.
\]
Thus, if $X_1,\ldots,X_N$, $X'_1, \ldots, X'_N$, $U$, and $U'$ are independent standard Gaussians, we have
\begin{align*}
    G_i(\sigma) &\overset{d}{=} \sqrt{a_N} X_i + \sqrt{b_N} U, \\
    G_{t,i}(\sigma) &\overset{d}{=} e^{-t}(\sqrt{a_N} X_i + \sqrt{b_N} U) + \sqrt{1-e^{-2t}} (\sqrt{a_N} X'_i + \sqrt{b_N} U'),
\end{align*}
jointly in $i$. In the language of Section \ref{subsec:framework}, we have
\[
    \cA=\{-1,+1\}, \quad d=r=2, \quad Y_i=(G_i,G_{t,i}), \quad
\mu_{s,N}=0,
\]
\[
    A_{s,N}= \sqrt{a_N} \begin{pmatrix} 1 & 0 \\ e^{-t} & \sqrt{1-e^{-2t}} \end{pmatrix},\quad B_{s,N} = \sqrt{Nb_N} \begin{pmatrix} 1 & 0 \\ e^{-t} & \sqrt{1-e^{-2t}} \end{pmatrix} ,
\]
\[
    \cK_{s,N}(\pi)=[\theta_{s,N,\gamma}(m(\pi)),\infty) \times [\theta_{s,N,\gamma'}(m(\pi))),
\qquad m(\pi)=\pi_+-\pi_-.
\]
If $T \in \R^2$ is a standard Gaussian vector, then $W = (\begin{smallmatrix} 1 & 0 \\ e^{-t} & \sqrt{1-e^{-2t}} \end{smallmatrix}) T$ is standard bivariate normal with correlation $\rho = e^{-t}$, i.e., it is centered Gaussian with covariance matrix $(\begin{smallmatrix} 1 & \rho \\ \rho & 1 \end{smallmatrix})$. Thus for $z = (z_1,z_2) \in \R^2$ we have
\begin{align*}
    p_{s,N}(\pi,z) &= \P(\sqrt{a_N} W_1 + \sqrt{Nb_N} z_1 \geq \theta_{s,N,\gamma}(m(\pi)), \\
    &\qquad \qquad \sqrt{a_N} W_2 + \sqrt{Nb_N} (\rho z_1 + \sqrt{1-\rho^2} z_2) \geq \theta_{s,N,\gamma'}(m(\pi))) \\
    &= \Phi_2\left(\frac{\sqrt{Nb_N} z_1 - \theta_{s,N,\gamma}(m(\pi))}{\sqrt{a_N}}, \frac{\sqrt{Nb_N}(\rho z_1 + \sqrt{1-\rho^2} z_2) - \theta_{s,N,\gamma'}(m(\pi))}{\sqrt{a_N}}; \rho\right).
\end{align*}
For fixed $\pi$ and $z$, this converges pointwise to
\[
    p_s(\pi,z) = \Phi_2\left( \sqrt{\frac{b}{a}} z_1 - \frac{1}{\sqrt{a}} \theta_{s,\gamma}(m(\pi)), \sqrt{\frac{b}{a}} (\rho z_1 + \sqrt{1-\rho^2} z_2) - \frac{1}{\sqrt{a}} \theta_{s,\gamma'}(m(\pi)); \rho\right).
\]
We will show that this convergence is uniform in the sense of \eqref{eqn:lambda-vs-p} in essentially the same way that we did in the proof of Theorem \ref{thm:complexity-p-spin}, then define $\Lambda(\pi,z)$ by \eqref{eqn:lambda-as-ps} and apply Corollary \ref{thm:main}. Indeed, define $x_1 = x_{1,N,s,\pi,z}$ and $y_1 = y_{1,N,s,\pi,z}$ implicitly through $p_{s,N}(\pi,z) = \Phi_2(x_1,y_1;\rho)$, and similarly $p_s(\pi,z) = \Phi_2(x_2,y_2;\rho)$. Then arguments as in \eqref{eqn:mixed-p-spin-l-estimate} show that, for fixed $M > 0$ and $N$ large enough, uniformly over $\|z\|_2 \leq M$ and $\pi \in \Delta_k$ and $s \in \mc{A}$, we have 
\[
    \abs{x_1}, \abs{y_1}, \abs{x_2}, \abs{y_2} \leq L
\]
(here we are using that $|\rho z_1 + \sqrt{1-\rho^2} z_2| \leq \|z\|_2 \leq M$). Still using the notation from Theorem \ref{thm:complexity-p-spin}, we have
\[
    \sup_{\abs{z} \leq M} \sup_{\pi \in \Delta_k} \sup_{s \in \mc{A}} \abs{x_1-x_2} \leq \epsilon_N,
\]
and similarly for $\abs{y_1-y_2}$, so that the same argument as before goes through (using that $g(x,y) = \log(\Phi(x,y;\rho))$ is locally Lipschitz). This verifies \eqref{eqn:main-assumption}, and the functions $p_s(\cdot,\cdot)$ are continuous in $\pi$ and $z$ and positive everywhere, so Corollary \ref{thm:main} applies, and gives
\[
    \Sigma^{\textup{dc}}(t,\xi,V,\gamma,\gamma') = \sup_{\pi \in \Delta_2} \sup_{z \in \R^2} \left\{ \pi_{+1} \log p_{+1}(\pi,z) + \pi_{-1} \log p_{-1}(\pi,z) - \frac{\|z\|_2^2}{2} + \mc{H}(\pi) \right\}.
\]
We conclude the proof of \eqref{eqn:complexity-disorder-persistence} by changing variables as before. 

What remains are the claims about monotonicity in $t$. By earlier arguments in the proof, we note that
\begin{align*}
    \P(\sigma \text{ is $\gamma$-gapped for $H_N$ and $\gamma'$-gapped for $H_{t,N}$}) = \P(G(\sigma) \in A_\sigma, \widetilde{G}(\sigma) \in B_\sigma),
\end{align*}
where $(G(\sigma),\widetilde{G}(\sigma))$ is a centered Gaussian vector with covariance $(\begin{smallmatrix} C_N & e^{-t} C_N \\ e^{-t} C_N & C_N \end{smallmatrix})$, and where 
\begin{align*}
    A_\sigma = \prod_{i=1}^N [\theta_{\sigma_i,N,\gamma}(m_N(\sigma)),+\infty), \qquad  
    B_\sigma = \prod_{i=1}^N [\theta_{\sigma_i,N,\gamma'}(m_N(\sigma)),+\infty).
\end{align*}
For each $\sigma$, the sets $A_\sigma$ and $B_\sigma$ are coordinatewise nondecreasing, and $C_N$ is nonsingular and entrywise nonnegative; thus Lemma \ref{lem:gaussian-interpolation} shows that $\P(\sigma \text{ is $\gamma$-gapped for $H_N$ and $\gamma'$-gapped for $H_{t,N}$})$ is nonincreasing in $t$. Since $\E[\mc{N}_N^{\textup{dc}}(t,\xi,V,\gamma,\gamma')]$ is a sum of these probabilities over $\sigma$, we conclude that it is also nonincreasing in $t$, and then the claim about $\Sigma^{\textup{dc}}(t,\xi,V,\gamma,\gamma')$ is immediate. 
\end{proof}


\subsection{Basics about persistence times}

\begin{proof}[Proof of Lemma \ref{lem:persistence-times}]
Say that a subset of $[0,\infty)$ is ``good'' if it is either empty, or equal to $[0,\infty)$, or equal to $[0,T]$ for some $T$. Since there are countably many $N$'s and $\sigma$'s, it suffices to prove that a single $\mc{G}_N(\sigma)$ is almost-surely good; since goodness is closed under finite unions, this immediately implies the almost-sure result for $\mc{G}_N$. 

From the proof of Theorem \ref{thm:disorder-persistence}, we see that, if we write $a_i = a_i(\sigma) = \theta_{\sigma_i,N,\gamma}(m(\sigma))$ and $b_i = b_i(\sigma) = \theta_{\sigma_i,N,\gamma'}(m(\sigma))$, then 
\begin{align*}
    \mc{G}_N(\sigma) &= \bigcap_{i=1}^N \{t \geq 0 : G_i(\sigma) \geq a_i \text{ and } e^{-t} G_i(\sigma) + \sqrt{1-e^{-2t}} \widetilde{G_i}(\sigma) \geq b_i\} =: \bigcap_{i=1}^N \mc{G}_{N,i}(\sigma).
\end{align*}
Since goodness is closed under finite intersections, it suffices to show that a single $\mc{G}_{N,i}(\sigma)$ is almost-surely good. Let us drop $i$ and $\sigma$ from the notation, so that we are considering the set $\{t \geq 0 : G \geq a \text{ and } e^{-t} G + \sqrt{1-e^{-2t}} \widetilde{G} \geq b\}$. Notice that, from our assumptions, we have
\[
    a = \frac{\gamma}{2} + \sigma_i h \geq \frac{\gamma'}{2} + \sigma_i h = b \geq 0.
\]

If $G < a$, this set is empty, and we are done.

Otherwise, we have $G \geq a \geq b \geq 0$. Since we are only considering almost-sure events, we can assume $G > a$ strictly (in particular, $G$ is strictly positive), and we can assume $\widetilde{G} \not\in \{0,b\}$. Now consider the function $f(x) = x G + \sqrt{1-x^2} \widetilde{G}$ for $x \in [0,1]$. This is not a convex combination, and it can take values not between $G$ and $\widetilde{G}$. However, we compute $f'(x) = G - \frac{x}{\sqrt{1-x^2}} \widetilde{G}$. The function $x \mapsto -\frac{x}{\sqrt{1-x^2}}$ is a monotone bijection between $[0,1]$ and $(-\infty,0]$; thus, if $\widetilde{G} < 0$, then $f'$ is positive for all arguments, so $f$ is strictly increasing for all arguments, meaning the map $g(t) = f(e^{-t})$ is strictly decreasing; since it equals $G > a \geq b$ when $t = 0$ and tends to $f(0) = \widetilde{G} < 0$ when $t \to +\infty$, and since it is continuous, the set $\mc{G}_{N,i}(\sigma) = \{t : g(t) \geq b\}$ must be good. If $\widetilde{G} > 0$, then $f'$ is positive until some $x_c$ and negative thereafter; this means that $g(t) = f(e^{-t})$ is increasing up to some $t_c$ and decreasing thereafter. But its value at $t = 0$ is $G > b$, and $\lim_{t \to +\infty} g(t) = \widetilde{G}$. If $\widetilde{G} > b$, then $\{t \geq 0 : g(t) \geq b\} = [0,\infty)$ and we are done; if $\widetilde{G} < b$, then $\{t \geq 0 : g(t) \geq b\}$ must be an interval of the form $[0,T]$, and we are also done. Since we assmumed $\widetilde{G} \neq b$, this completes the proof. 
\end{proof}

\begin{rem}
\label{rem:almost-sure-persistence}
It would also be interesting to understand the structure of the set $\mc{G}_N$ of persistence times outside the assumption $\gamma \geq \gamma' \geq 2\abs{h}$. The proof above strongly relies on $a_i$ being positive; otherwise it may fail substantially. For example, if $a = G = \widetilde{G} = -1$ and $b = -1-\epsilon$ for small $\epsilon$, then the set $\{t \geq 0 : -e^{-t} - \sqrt{1-e^{-2t}} \geq -1-\epsilon\}$ is the union of two disjoint intervals. 
\end{rem}


\subsection{Analyzing the variational formula}

\begin{proof}[Proof of Corollary \ref{cor:disorder-persistence-analysis}]
Equation \eqref{eqn:disorder-persistence-simplified-formula} follows as in the previous section, since it is still of the form $\sup_{m \in [-1,1]} \{\frac{1+m}{2} \alpha + \frac{1-m}{2}\beta + \fh(m)\}$. 

The function being optimized in \eqref{eqn:disorder-persistence-simplified-formula} is clearly continuous in all of its parameters, and it tends to $-\infty$ as $\min(z_1,z_2) \to +\infty$, locally uniformly in all of its parameters. Thus the supremum over $z_1, z_2$ can be locally restricted to a supremum over a compact set, allowing $\Sigma^{\textup{dc}}$ to inherit its continuity. This completes the proof of item \ref{item:disorder-persistence-analysis-1}. A priori this only gives continuity in $t$ for only fixed $t$, whereas momentarily we will need it as $t \to +\infty$; for this extension, we simply write $\exp(-t) = \rho \in [0,1]$, and run the argument again with $\rho$ in place of $t$, to find continuity at $\rho = 0$, thus at $t = +\infty$.

When $t = 0$, we have $\Phi_2(x,y;1) = \Phi(\min(x,y))$; thus
\[
    f^{(s)}_{t,a,b,\gamma,\gamma'}(z_1,z_2) = \Phi\left( \sqrt{\frac{b}{a}} z_1 - \frac{1}{\sqrt{a}}\theta_{s,\max(\gamma,\gamma')} \right).
\]
The only remaining dependence in $z_2$ is $-z_2^2/2$, which is clearly optimized at zero, and what remains takes the form \eqref{eqn:p-spin-simplified-formula}, proving \eqref{eqn:dc-limit-t-zero}. When $t \to +\infty$,
\begin{align*}
    \lim_{t \to +\infty} f^{(s)}_{t,a,b,\gamma,\gamma'} &= \Phi_2\left(\sqrt{\frac{b}{a}} z_1 - \frac{1}{\sqrt{a}} \theta_{s,\gamma}, \sqrt{\frac{b}{a}} z_2 - \frac{1}{\sqrt{a}} \theta_{s,\gamma'}; 0 \right) \\
    &= \Phi\left(\sqrt{\frac{b}{a}} z_1 - \frac{1}{\sqrt{a}} \theta_{s,\gamma} \right) \Phi\left( \sqrt{\frac{b}{a}} z_2 - \frac{1}{\sqrt{a}} \theta_{s,\gamma'} \right).
\end{align*}
When $h = 0$ this no longer depends on $s$. From the previous argument, we can exchange the limit $t \to +\infty$ with the supremum; thus the formula \eqref{eqn:disorder-persistence-simplified-formula} becomes
\begin{align*}
    &\lim_{t \to +\infty} \Sigma^{\textup{dc}}(t,\xi,0,\gamma,\gamma') + \log(2) \\
    &= \sup_{z_1, z_2 \geq 0} \left\{ \log\left(2 \Phi\left(\sqrt{\frac{b}{a}} z_1 - \frac{\gamma}{2\sqrt{a}}\right) \Phi\left(\sqrt{\frac{b}{a}} z_2 - \frac{\gamma'}{2\sqrt{a}}\right) \right) - \frac{z_1^2}{2} - \frac{z_2^2}{2} \right\} + \log(2) \\
    &= \sup_{z_1 \geq 0} \left\{ \log\left(2 \Phi\left(\sqrt{\frac{b}{a}} z_1 - \frac{\gamma}{2\sqrt{a}}\right) \right) - \frac{z_1^2}{2} \right\} + \sup_{z_2 \geq 0} \left\{ \log\left(2 \Phi\left(\sqrt{\frac{b}{a}} z_2 - \frac{\gamma'}{2\sqrt{a}}\right) \right) - \frac{z_2^2}{2} \right\} \\
    &= \Sigma^{\textup{mix}}(\xi,0,\gamma) + \Sigma^{\textup{mix}}(\xi,0,\gamma'),
\end{align*}
which proves item \eqref{item:disorder-persistence-analysis-2}.

For item \eqref{item:disorder-persistence-analysis-3}, we note that, from the assumption $\gamma \geq \gamma' \geq 2\abs{h}$, we may apply Lemma \ref{lem:persistence-times}, which tells us that $\{T_n > t_c + \epsilon\} \subseteq \{\mc{N}^{\textup{dc}}_N(t,\xi,h,\gamma,\gamma') \geq 1\}$. (Otherwise, $\{T_n > t_c + \epsilon\}$ is a union over uncountably many times.) Since $\E[\mc{N}^{\textup{dc}}_N(t,\xi,h,\gamma,\gamma')]$ is exponentially small for $t > t_c$, the proof of \eqref{item:disorder-persistence-analysis-3} is just Markov's inequality, as in item \eqref{item:p-spin-analysis-2} for mixed spin glasses.

Finally, we address item \eqref{item:disorder-persistence-analysis-4}. From \eqref{eqn:dc-limit-t-zero}, we have $\Sigma^{\textup{dc}}(0,p,0,0,0) = \Sigma^{\textup{mix}}(p,0,0)$, which is strictly positive for all $p$ by item \eqref{item:p-spin-analysis-2} of Corollary \ref{cor:p-spin-analysis}. Since $\Sigma^{\textup{dc}}(t,p,0,0,0)$ is continuous in $t$, it is also positive for $t$ in some small neighborhood of zero; thus $t_c(p,0,0,0)$ is strictly positive for each $p$.

The same continuity tells us that $t_c$ is finite if $\lim_{t \to +\infty} \Sigma^{\textup{dc}}(t,p,0,0,0)$ is negative, and infinite if $\lim_{t \to +\infty} \Sigma^{\textup{dc}}(t,p,0,0,0)$ is positive. Thus, comparing to \eqref{eqn:dc-limit-t-infinity}, we see that we need only show
\[
    \Sigma^{\textup{mix}}(p,0,0) \begin{cases} < \frac{\log(2)}{2} & \text{ for } p = 2, 3, \\ > \frac{\log(2)}{2} & \text{ for } p \geq 4. \end{cases}
\]
Since the left-hand side is nondecreasing in $p$ (from item \eqref{item:p-spin-analysis-4} of Corollary \ref{cor:p-spin-analysis}), it suffices just to check the cases $p = 3$ and $p = 4$. Using the formula \eqref{eqn:pure-p-spin-formula} for $\Sigma^{\textup{mix}}(p,0,0)$, we see that it suffices to show 
\[
    \sup_{z \geq 0} \left\{ \log(\Phi(z\sqrt{p-1})) - \frac{z^2}{2} \right\} \begin{cases} < -\frac{\log(2)}{2} & \text{ for } p = 3, \\ > -\frac{\log(2)}{2} & \text{ for } p = 4. \end{cases}
\]
This is clear from plotting the functions, but we can also see it analytically. The $p = 4$ case is easier, since we can just exhibit a $z$; the choice $z = 1/2$ works, since then the left-hand side is $\log(\Phi(\frac{\sqrt{3}}{2})) - 1/8 \approx -0.3397$, and the right-hand side is $\approx -0.3466$. When $p = 3$, we change variables to write the variational problem as $\sup_{u \geq 0} \{\log(\Phi(u)) - \frac{u^2}{4}\}$. The argument is bounded above by $-u^2/4$, which is smaller than $-\log(2)/2$ when $u > \sqrt{2\log 2}$; on the interval $u \in [0,\sqrt{2\log 2} + \epsilon]$, we need to show that $\Phi(u) e^{-u^2/4} < \frac{1}{\sqrt{2}}$. We have $\Phi(u) = 1/2 + \int_0^u \diff s \exp(-s^2/2)/\sqrt{2\pi} \leq 1/2 + \frac{u}{\sqrt{2\pi}}$, so that it suffices to show that 
\[
    r(u) = \left(\frac{1}{2} + \frac{u}{\sqrt{2\pi}} \right) e^{-\frac{u^2}{4}} < \frac{1}{\sqrt{2}}
\]
for $u \in [0,\sqrt{2\log 2} + \epsilon]$. We have $r'(u) = (\frac{1}{\sqrt{2\pi}} -\frac{u}{4} - \frac{u^2}{2\sqrt{2\pi}}) \exp(-u^2/4)$. The bracket is a quadratic in $u$, which is positive for $0 \leq u \leq u_\ast = \frac{\sqrt{2(\pi+16)}-\sqrt{2\pi}}{4} \approx 0.92$ and negative otherwise; thus $r(u) \leq r(u_\ast) \approx 0.702$, whereas $1/\sqrt{2} \approx 0.707$, completing the proof.
\end{proof}


\section{Application III: Ghatak-Sherrington model}
\label{sec:GS}

We now apply the general method to the Ghatak-Sherrington (GS)
model. The main result is a variational formula for the annealed complexity
of local minima and gapped states, allowing arbitrary crystal and
external fields. In zero external field case, we further show the existence of an annealed topological trivialization transition as the crystal field changes. We begin with the definition of the GS model.


\subsection{Setting and main results}

For $N>1$, the configuration space of GS model is $\cA^N = \{-1,0,+1\}^N$. Its associated Hamiltonian is 
\[
H_N^{\rm{GS}}(\gs) = \frac{1}{\sqrt{N}}\sum_{i<j} J_{ij} \gs_i \gs_j + D \sum_{i=1}^N \gs_i^2 + h \sum_{i=1}^N \gs_i
\]
where $D,h\in \bR$ are the crystal and external field parameter respectively, and $(J_{ij})_{i,j \in [N]}$ is symmetric matrix with zero diagonals and i.i.d. standard Gaussians in the upper triangular entries.  
For a configuration $\gs \in \cA^N$, we define the following one-spin flip perturbation: for $t \in \cA$ and $t \neq \gs_i$,
\begin{align}
    (\gs^{(i,t)})_j := \begin{cases}
        t,  & \text{ if $j=i$}\\
        \gs_j, & \text{ if $j\neq i$}.
    \end{cases}
\end{align}
In other words, this replaces the $i$-th coordinate $\sigma_i$ with $t$. This is a natural generalization of the standard one-spin perturbation of the binary value case. Similarly, for $\gamma \ge 0$, we define $\gs$ as $\gamma$-gapped local minimum if 
\begin{align}\label{eq:local-min-gs}
H_N^{\rm{GS}}(\gs^{(i,t)}) - H_N^{\rm{GS}}(\gs)\ge \gamma \ \ \text{for all $i\in [N], \quad t\in \cA\setminus \{\gs_i\}$}. 
\end{align}
We use the following notation to denote the number of local minima in the GS case: 
\begin{align}
    \cN_N^{\rm{GS}}(D,h,\gamma) = |\{\gs \in \cA^N: \gs \ \text{satisfies}\ \eqref{eq:local-min-gs} \}|.
\end{align}

To state our main results, we specify the notation of Section~\ref{sec:general} for the GS model. For $\gs \in \cA^N$, we write
\(
q_N(\gs)=\frac1N\sum_{i=1}^N\gs_i^2
\)
for the fraction of its entries in $\{\pm1\}$. For a probability measure $\pi$ on $\cA$, identified with $(\pi_-,\pi_0,\pi_+)\in\gD_3$ where $\pi_\pm=\pi(\{\pm1\})$ and $\pi_0=\pi(\{0\})$, we write the analogous quantity as $q(\pi)=\pi_++\pi_-$, so that $q_N(\gs)=q(\pi(\gs))$ for the empirical profile \eqref{eq:profile}. Finally, we write the binary entropy as
\begin{align}
\label{eq:GS-q}
 h_2(q)=-q\log q-(1-q)\log(1-q).
\end{align}
As before, we use the convention $0\log0=0$ throughout the paper. Set
\begin{equation}
\label{eq:GS-theta}
 \theta=\theta(D,\gamma):=\max\{D+\gamma,\gamma/2\}.
\end{equation}
For $q\in(0,1]$ and $z\in\R$, define
\begin{align}
 p_s(q,z)
 &=\Phi\left(\frac{z-\theta-sh}{\sqrt q}\right),
 \qquad s\in\{-1,+1\},
 \label{eq:GS-ps}\\
 p_0(q)
 &=\begin{cases}
 \displaystyle
 \Phi\left(\frac{h+D-\gamma}{\sqrt q}\right)
 -\Phi\left(\frac{h-D+\gamma}{\sqrt q}\right),&D>\gamma,\\[1ex]
 0,&D\le\gamma.
 \end{cases}
 \label{eq:GS-pzero}
\end{align}
In the sequel, we write $p_+$ and $p_-$ for $p_{+1}$ and $p_{-1}$.
Finally, define
\begin{equation}
\label{eq:GS-zero-rate}
 \Sigma_0(D,h,\gamma)=
 \begin{cases}
 0,&D\ge\gamma+|h|,\\
 -\infty,&D<\gamma+|h|.
 \end{cases}
\end{equation}
This is the contribution of the all-zero configuration, whose energy
increments are $D+h$ and $D-h$.

\begin{thm}[Gapped complexity for the GS model]
\label{thm:complexity-GS}

For every $D,h\in\R$ and $\gamma\ge0$, the limit
\[
 \Sigma^{\mathrm{GS}}(D,h,\gamma)
 =\lim_{N\to\infty}\frac1N\log\E\cN_N^{\mathrm{GS}}(D,h,\gamma)
\]
exists and is given by 
\[
\max\Big\{\gS_0(D,h,\gc), \sup_{\pi \in \gD_3: q(\pi)>0} \sup_{z\ge 0}\Big[  \cH(\pi) +\pi_0 \log p_0(q) + \sum_{s \in \{-1,+1\}} \pi_s \log p_s(q,z) - \frac{z^2}{2}  \Big] \Big\},
\]
where we recall the notation $\cH(\pi) = -\sum_{s\in \{0,-1,+1\}} \pi_s \log \pi_s$ from \eqref{eqn:entropy}.
\end{thm}
Note that comparing to the results of binary-valued model in Section~\ref{sec:p-spin}, there is new term $\gS_0(D,h,\gc)$ contributing to the annealed complexity. By analyzing the above variational formula, one can deduce that, when the external field vanishes, the GS model exhibits topological trivialization for strong crystal field. We defer the proof of Theorem~\ref{thm:complexity-GS} to Section~\ref{ssec:GS-complexity-proof}. Now we state the topological trivialization results as follows.

\begin{prop}[Topological trivialization of GS model]\label{prop:GS-topo-trivial} 
Consider the number of ordinary local minima for the GS model in the zero external field, namely $\gamma = 0, h=0$ case. Then there exists a constant $D_c \in (0,\infty)$ such that the following holds:
\begin{enumerate}
    \item [(i)] For $D< D_c$, we have $\gS^{GS}(D,0,0)>0$.
    \item [(ii)] For $D\ge D_c$, we have $\gS^{GS}(D,0,0)=0$.
\end{enumerate}
\end{prop}

\begin{rem}[Comparison with~\cite{daCostaDeAraujo2000}]
\label{rem:da-costa-de-araujo}
Da Costa and de Ara\'ujo \cite{daCostaDeAraujo2000}\footnote{The arXiv version contains some typos, after corrections, the complexity formulas match exactly.} computed, non-rigorously, the annealed complexity of solutions of the zero-temperature TAP equations for the GS model in zero external field. These solutions are precisely the local minima \eqref{eq:local-min-gs} with $h=\gc=0$, so they computed $\gS^{\mathrm{GS}}(D,0,0)$. Indeed, after eliminating a Lagrange multiplier, their saddle-point formula \cite[eqs. (17)--(18)]{daCostaDeAraujo2000} coincides with ours in Lemma~\ref{lem:GS-reduction-K}. Both give a complexity equal to the SK value $\gS^{\mathrm{mix}}(2,0,0)\approx0.1992$ for $D\le0$, which increases to a maximum and then decreases to zero at a finite crystal field. Proposition~\ref{prop:GS-topo-trivial} proves rigorously that such a threshold $D_c$ exists and an implicit definition is given in~\eqref{eq:GS-annealed-critical-definition}. The numerical values differ slightly. In our normalization, \cite{daCostaDeAraujo2000} reports a maximum near $D\approx0.550$ and vanishing near $D\approx1.225$, whereas evaluating \eqref{eq:GS-complexity-simple} gives a maximum near $D\approx0.596$ and $D_c\approx1.288$, for instance, $\gS^{\mathrm{GS}}(1.225,0,0)\approx0.045>0$. Since the formulas agree, the difference might be due to the numerical computations.
\end{rem}


\subsection{Establishing the variational formula}\label{ssec:GS-complexity-proof}

\begin{proof}[Proof of Theorem~\ref{thm:complexity-GS}]

Fix $D,h\in\R$ and $\gamma\ge0$. For $\gs_i=s$ and
$t\in\cA\setminus\{s\}$,
\[
 H_N^{\mathrm{GS}}(\gs^{(i,t)})-H_N^{\mathrm{GS}}(\gs)
 =(t-s)\left(\frac1{\sqrt N}\sum_{j\ne i}J_{ij}\gs_j+h\right)
   +D(t^2-s^2).
\]
Define the full local field
\begin{equation}
\label{eq:GS-G}
 Y_i(\gs):=-\frac1{\sqrt N}\sum_{j\ne i}J_{ij}\gs_j-h-2D\gs_i.
\end{equation}
Then
\begin{equation}
\label{eq:GS-increment}
 H_N^{\mathrm{GS}}(\gs^{(i,t)})-H_N^{\mathrm{GS}}(\gs)
 =(s-t)Y_i(\gs)+D(t-s)^2.
\end{equation}
In the notation of Section~\ref{sec:general}, we have
$\cA=\{-1,0,+1\}$, $\abs{\cA}=3$, $d=r=1$, and
\[
 \cK_{s,N}(\pi)=\cK_s
 :=\bigcap_{t\in\cA\setminus\{s\}}
 \{y\in\R:(s-t)y+D(t-s)^2\ge\gamma\}.
\]
For $s=\pm1$, the two constraints are
\[
 sy+D\ge\gamma,\qquad 2sy+4D\ge\gamma
 \quad\Longleftrightarrow\quad
 sy\ge\max\{\gamma-D,\gamma/2-2D\}=\theta-2D,
\]
where $\theta$ is defined in~\eqref{eq:GS-theta}. 

For $s=0$, they are $D-y\ge\gamma$ and $D+y\ge\gamma$. Hence
\begin{equation}
\label{eq:GS-stability-sets}
 \cK_+=[\theta-2D,\infty),\qquad
 \cK_-=(-\infty,2D-\theta],\qquad
 \cK_0=[\gamma-D,D-\gamma],
\end{equation}
where $\cK_0=\varnothing$ when $D<\gamma$. In particular,
\[
 \gs\text{ satisfies }\eqref{eq:local-min-gs}
 \quad\Longleftrightarrow\quad
 Y_i(\gs)\in\cK_{\gs_i}\text{ for every }i.
\]
Thus the count in \eqref{eq:count} equals
$\cN_N^{\mathrm{GS}}(D,h,\gamma)$. To apply
Theorem~\ref{thm:general}, we need to verify the conditions \eqref{eq:CF},
\eqref{cond:gene-upper}, and \eqref{cond:gene-lowers}.

--\noindent\textbf{Verification of \eqref{eq:CF}.}
For a configuration with profile $\pi$, denote $q=q(\pi)$. One can compute the mean and variance of the local fields:
\begin{equation}
\label{eq:GS-covariance}
\begin{aligned}
 \E Y_i(\sigma) &=-h-2D\gs_i,\\
 \Var(Y_i(\sigma))&=\frac1N\sum_{j\ne i}\gs_j^2
                  =q-\frac{\gs_i^2}{N},\\
 \Cov(Y_i(\sigma),Y_j(\sigma))&=\frac1N\sum_{u\ne i}\sum_{v\ne j}
        \gs_u\gs_v\E[J_{iu}J_{jv}]
                  =\frac{\gs_i\gs_j}{N},\qquad i\ne j.
\end{aligned}
\end{equation}
The fields are jointly Gaussian, so \eqref{eq:GS-covariance} verifies
\eqref{eq:CF} with
\begin{equation}
\label{eq:GS-dictionary}
\begin{gathered}
 \mu_{s,N}(\pi)=-h-2Ds,\\
 (A_{s,N}(\pi),B_{s,N}(\pi))=
 \begin{cases}
 \left(\sqrt{q-2s^2/N},s\right),&Nq\ge2,\\[1mm]
 \left(\sqrt q\,\mathbf1_{\{s=0\}},0\right),&Nq\le1.
 \end{cases}
\end{gathered}
\end{equation}
In particular, for $Nq\ge2$,
\begin{equation}
\label{eq:GS-common-factor}
 (Y_i(\gs))_{i=1}^N\equald
 \left(-h-2D\gs_i+\sqrt{q-\frac{2\gs_i^2}{N}}\,X_i
                   +\frac{\gs_i}{\sqrt N}Z\right)_{i=1}^N.
\end{equation}
For $s=\pm1$ and $Nq>2$, substitution in \eqref{eq:p-def} gives
\begin{align*}
 p_{s,N}(\pi,z)
 =\P\big(-sh-2D+\sqrt{q-2/N}\,(sX)+z\ge\theta-2D\big)
 =\Phi\left(\frac{z-\theta-sh}{\sqrt{q-2/N}}\right).
\end{align*}
For the zero type and every $q>0$,
\begin{align*}
 p_{0,N}(\pi,z)
 &=\P\big(-h+\sqrt q\,X\in[\gamma-D,D-\gamma]\big)\\
 &=\P\big(\sqrt q\,X\in[h-D+\gamma,h+D-\gamma]\big)=p_0(q).
\end{align*}
Thus for $s = \pm 1$ and $Nq>2$, we have
\begin{equation}
\label{eq:GS-finite-probabilities}
 p_{s,N}(\pi,z)=\Phi\left(\frac{z-\theta-sh}{\sqrt{q-2/N}}\right)
 ,\qquad
 p_{0,N}(\pi,z)=p_0(q).
\end{equation}
At the remaining occupations the same coefficients give, for $s=\pm1$,
\[
 p_{s,N}(\pi,z)=
 \begin{cases}
 \mathbf1_{\{z\ge\theta+sh\}},&Nq=2,\\
 \mathbf1_{\{0\ge\theta+sh\}},&Nq\le1.
 \end{cases}
\]
At $\pi^{(0)}=(0,1,0)$, only the zero type is present, and
\[
 p_{0,N}(\pi^{(0)},z)
 =\mathbf1_{\{D\ge\gamma+|h|\}},\qquad
 \Lambda_N(\pi^{(0)},z)=\Sigma_0(D,h,\gamma).
\]
Define the candidate function by
\begin{equation}
\label{eq:GS-Lambda}
 \Lambda(\pi,z)=
 \begin{cases}
 \displaystyle\pi_0\log p_0(q)+\sum_{s=\pm1}\pi_s\log p_s(q,z),
                         &q=q(\pi)>0,\\[1mm]
 \Sigma_0(D,h,\gamma),&\pi=\pi^{(0)}.
 \end{cases}
\end{equation}
By the conventions for zero probabilities,
$\Lambda_N,\Lambda\in[-\infty,0]$.

--\noindent\textbf{Verification of \eqref{cond:gene-upper}.}
Let $\pi_N\in\cP_N$ and $z_N\in\R$ satisfy
$(\pi_N,z_N)\to(\pi,z)$, and put $q_N=q(\pi_N)$ and $q=q(\pi)$.
If $q>0$, then $Nq_N>2$ for all sufficiently large $N$, and
\eqref{eq:GS-finite-probabilities} gives
\begin{equation}
\label{eq:GS-positive-occupation-convergence}
 p_{s,N}(\pi_N,z_N)
 =\Phi\left(\frac{z_N-\theta-sh}{\sqrt{q_N-2/N}}\right)
 \longrightarrow
 \Phi\left(\frac{z-\theta-sh}{\sqrt q}\right)
 =p_s(q,z)>0,\qquad s=\pm1.
\end{equation}
If $D>\gamma$, then also $p_0(q_N)\to p_0(q)>0$, so
\begin{align*}
 \Lambda_N(\pi_N,z_N)
 &=\pi_{N,0}\log p_0(q_N)
     +\sum_{s=\pm1}\pi_{N,s}\log p_{s,N}(\pi_N,z_N)\\
 &\to
   \pi_0\log p_0(q)+\sum_{s=\pm1}\pi_s\log p_s(q,z)
   =\Lambda(\pi,z).
\end{align*}
If $D\le\gamma$ and $0<q<1$, then eventually $0<q_N<1$ and
$p_0(q_N)=0$. Hence
\[
 \limsup_{N\to\infty}\Lambda_N(\pi_N,z_N)
 =-\infty=\Lambda(\pi,z).
\]
If $D\le\gamma$ and $q=1$, the zero-type term is nonpositive, and
\eqref{eq:GS-positive-occupation-convergence} yields
\begin{align*}
 \limsup_{N\to\infty}\Lambda_N(\pi_N,z_N)
 &\le\lim_{N\to\infty}
       \sum_{s=\pm1}\pi_{N,s}\log p_{s,N}(\pi_N,z_N)\\
 &=\sum_{s=\pm1}\pi_s\log p_s(1,z)=\Lambda(\pi,z).
\end{align*}
Thus \eqref{cond:gene-upper} holds at positive occupation. It remains to consider $q=0$, so $\pi=\pi^{(0)}$. If
$D\ge\gamma+|h|$, then
\[
 \limsup_{N\to\infty}\Lambda_N(\pi_N,z_N)
 \le0=\Lambda(\pi^{(0)},z).
\]
If $D<\gamma+|h|$ and $D\le\gamma$, then $q_N<1$ eventually and
\[
 \Lambda_N(\pi_N,z_N)=
 \begin{cases}
 \Sigma_0(D,h,\gamma)=-\infty,&q_N=0,\\
 -\infty,&0<q_N<1.
 \end{cases}
\]
Finally, if $\gamma<D<\gamma+|h|$, put $a=|h|-(D-\gamma)>0$.
For $q'>0$, the interval $[h-D+\gamma,h+D-\gamma]$ lies on one side
of zero at distance $a$. Hence
\[
 p_0(q')\le\Phi\left(-\frac{a}{\sqrt{q'}}\right)
 \le\inf_{t\ge0}e^{-ta/\sqrt{q'}+t^2/2}
 =e^{-a^2/(2q')}.
\]
For all sufficiently large $N$,
\[
 \Lambda_N(\pi_N,z_N)\le
 \begin{cases}
 \displaystyle(1-q_N)\log p_0(q_N)
       \le-\frac{(1-q_N)a^2}{2q_N},&q_N>0,\\[1mm]
 -\infty,&q_N=0.
 \end{cases}
\]
Since $q_N\to0$, this gives
$\limsup_N\Lambda_N(\pi_N,z_N)=-\infty=\Lambda(\pi^{(0)},z)$.
Thus \eqref{cond:gene-upper} holds on the whole simplex.

--\noindent\textbf{Verification of \eqref{cond:gene-lowers}.}
Fix $(\pi,z)$ with $\Lambda(\pi,z)>-\infty$ and $\varepsilon>0$.
Suppose first that $q=q(\pi)>0$. Then either $D>\gamma$ or $q=1$.
Choose empirical profiles $\pi_N\to\pi$, keeping $q(\pi_N)=1$ when
$D\le\gamma$. The probability limits above, now with $z_N=w$ fixed,
give
\[
 \Lambda_N(\pi_N,w)\longrightarrow\Lambda(\pi,w)\in\R
 \qquad\text{for every fixed }w\in\R.
\]
The function $w\mapsto\Lambda(\pi,w)$ is continuous, so choose
$\delta>0$ such that
\[
 \Lambda(\pi,z-\delta)\ge\Lambda(\pi,z)-\varepsilon/2.
\]
For all sufficiently large $N$, $w\mapsto\Lambda_N(\pi_N,w)$ is
non-decreasing: the occupied-type probabilities in
\eqref{eq:GS-finite-probabilities} increase with $w$, while the
zero-type probability is independent of $w$. Using convergence at
the single point $w=z-\delta$, for all sufficiently large $N$ we obtain
\begin{align*}
 \inf_{|w-z|\le\delta}\Lambda_N(\pi_N,w)
 =\Lambda_N(\pi_N,z-\delta) \ge\Lambda(\pi,z-\delta)-\varepsilon/2 
 \ge\Lambda(\pi,z)-\varepsilon.
\end{align*}
If $q(\pi)=0$, finiteness means $D\ge\gamma+|h|$.
Taking $\pi_N=\pi^{(0)}$ for every $N$ and any $\delta>0$ gives
\[
 \inf_{|w-z|\le\delta}\Lambda_N(\pi^{(0)},w)
 =0=\Lambda(\pi^{(0)},z)
 \ge\Lambda(\pi^{(0)},z)-\varepsilon.
\]
This proves \eqref{cond:gene-lowers} at every finite value of
$\Lambda$. Now applying Theorem~\ref{thm:general}, it gives
\[
 \Sigma^{\mathrm{GS}}(D,h,\gamma)
 =\sup_{\substack{\pi\in\Delta_3\\z\in\R}}
 \left\{\cH(\pi)+\Lambda(\pi,z)-\frac{z^2}{2}\right\}.
\]
At $\pi^{(0)}$, $\cH(\pi^{(0)})=0$ and
\[
 \sup_{z\in\R}\left\{\Lambda(\pi^{(0)},z)-\frac{z^2}{2}\right\}
 =\Sigma_0(D,h,\gamma).
\]
For $q>0$, $p_s(q,|z|)\ge p_s(q,z)$ for $s=\pm1$. Substitution of \eqref{eq:GS-Lambda} therefore proves Theorem~\ref{thm:complexity-GS}.    
\end{proof}


\subsection{Analyzing the variational formula}

To prepare the proof of Proposition~\ref{prop:GS-topo-trivial}, we first
simplify the complexity formula in Theorem~\ref{thm:complexity-GS}. As in Proposition~\ref{prop:GS-topo-trivial}, throughout this subsection
we take $\gamma=h=0$. Recall that $q=\pi_-+\pi_+$ is the fraction of
nonzero spins and $h_2(q)=-q\log q-(1-q)\log(1-q)$.

\begin{lem}
\label{lem:GS-reduction-K}
For $\gamma=h=0$ and $D\in\R$, we have
\begin{equation}
\label{eq:GS-complexity-simple}
 \Sigma^{\mathrm{GS}}(D,0,0)=
 \begin{cases}
  \Sigma^{\mathrm{mix}}(2,0,0),&D\le0,\\[1mm]
  \displaystyle\max\left\{0,\sup_{q\in(0,1]}
        K(D/\sqrt q,q)\right\},&D>0,
 \end{cases}
\end{equation}
where $\Sigma^{\mathrm{mix}}(2,0,0)$ is the SK complexity in
Corollary~\ref{cor:p-spin-analysis}, and, for $t>0$ and $q\in(0,1]$,
\begin{equation}
\label{eq:GS-rescaled-objective}
 K(t,q):=h_2(q)+(1-q)\log(2\Phi(t)-1)
       +q\max_{u\ge0}\left\{\log\bigl(2\Phi(u-t)\bigr)
                              -\frac{u^2}{2}\right\}.
\end{equation}
\end{lem}

\begin{proof}
    The only configuration with $q=0$ is $\mathbf0$. By the definition of local minima, since $h= 0, \gamma = 0$, replacing one
of its spins by $\pm 1$ changes the energy by exactly $D$, so $\mathbf 0$ is a local minimum if and only if $D\ge0$. This is the term $\Sigma_0(D,0,0)$ in Theorem~\ref{thm:complexity-GS}, which by \eqref{eq:GS-zero-rate} equals $0$ if $D\ge0$ and $-\infty$ otherwise. For $\pi\in\gD_3$ with $q=q(\pi)>0$, since $h=\gamma=0$, we have $p_-=p_+$. At fixed $q,z$,
\[
 \pi_0\log p_0(q)+\sum_{s=\pm1}\pi_s\log p_s(q,z)
 =(1-q)\log p_0(q)+q\log p_+(q,z),
\]
so the objective in Theorem~\ref{thm:complexity-GS} depends on $\pi$ beyond $q$ only through
\[
 \cH(\pi)=h_2(q)+q\,h_2(\pi_+/q)\le h_2(q)+q\log2.
\]
Equality holds at $\pi_+=\pi_-=q/2$. Hence the supremum of the objective over $\pi\in\gD_3$ with $q(\pi)=q$ is
\begin{equation}
\label{eq:GS-fixed-q-full-maximum}
 h_2(q)+(1-q)\log p_0(q)
       +q\log\bigl(2p_+(q,z)\bigr)-\frac{z^2}{2}.
\end{equation}
If $p_0(q)=0$ and $q<1$, the objective equals $-\infty$ for every such $\pi$.
For $q=1$, For $q=1$, the term $\pi_0\log p_0(q)$ vanishes by the convention $0\log0=0$. If $D\le0$, then $p_0(q)=0$ for $q>0$ and
$p_+(q,z)=\Phi(z/\sqrt q)$. Only $q=1$ contributes to the supremum over $\pi$ with $q(\pi)>0$, so Theorem~\ref{thm:complexity-GS} gives
\[
 \Sigma^{\mathrm{GS}}(D,0,0)
 =\max\left\{\Sigma_0(D,0,0),
   \sup_{z\ge0}\left[\log(2\Phi(z))-\frac{z^2}{2}\right]\right\}
 =\Sigma^{\mathrm{mix}}(2,0,0)>0.
\]
Here the second entry is the positive SK complexity by
Corollary~\ref{cor:p-spin-analysis}, while the first is at most zero. If $D>0$, then $\Sigma_0(D,0,0)=0$ and
\[
 p_0(q)=2\Phi(D/\sqrt q)-1,\qquad
 p_+(q,z)=\Phi((z-D)/\sqrt q).
\]
For each fixed $q>0$, substitute $t=D/\sqrt q$ and $z=\sqrt q\,u$
in \eqref{eq:GS-fixed-q-full-maximum}. It becomes
\[
 h_2(q)+(1-q)\log(2\Phi(t)-1)
 +q\left[\log(2\Phi(u-t))-\frac{u^2}{2}\right].
\]
The substitution maps $z\ge0$ bijectively to $u\ge0$.
The inner maximum is attained because its continuous objective
tends to $-\infty$ as $u\to\infty$.
Maximizing over $u$, then over $q$, gives
\eqref{eq:GS-complexity-simple}.
\end{proof}

\begin{lem}
\label{lem:GS-positivity-below-pair}
Let $D_0>0$. If there exists $q_0\in(0,1]$ such that
\(
 K(D_0/\sqrt{q_0},q_0)\ge0,
\)
then
\[
 \Sigma^{\mathrm{GS}}(D,0,0)>0\qquad\text{for every }0<D<D_0.
\]
\end{lem}

\begin{proof}
We first identify an explicit $q$ where $K$ is strictly positive.
For any $t>0$, take $q=2\Phi(-t)\in(0,1)$. Then
$2\Phi(t)-1=1-q$, and the expression defining $K$ before maximization,
evaluated at $u=0$, equals
\[
 h_2(q)+(1-q)\log(2\Phi(t)-1)+q\log(2\Phi(-t))
 =h_2(q)+(1-q)\log(1-q)+q\log q=0.
\]
Its derivative with respect to $u$ at zero is
\[
 q\left.\frac{\mathrm d}{\mathrm du}
 \left[\log(2\Phi(u-t))-\frac{u^2}{2}\right]\right|_{u=0}
 =q\frac{\Phi'(-t)}{\Phi(-t)}>0.
\]
A sufficiently small $u>0$ therefore gives a positive value, so $ K(t,2\Phi(-t))>0$.

Fix $q_0\in(0,1]$ satisfying the hypothesis. For $D\in(0,D_0)$, let
\[
 t_0=\frac{D_0}{\sqrt{q_0}},\qquad \bar q=2\Phi(-t_0).
\]
Then $D_0=t_0\sqrt{q_0}$, and $K(t_0,q_0)\ge0$, and
$K(t_0,\bar q)>0$. By \eqref{eq:GS-complexity-simple}, it suffices
to construct $t>0$, $q\in(0,1]$ with $D=t\sqrt q$ and $K(t,q)>0$.
If $D\ge t_0\sqrt{\bar q}$, keep $t=t_0$.
Set $q=(D/t_0)^2$. Then
\[
 \bar q\le q<q_0\le1,\qquad D=t_0\sqrt q.
\]
Thus $q_0-\bar q>0$. The formula for $K$ is binary entropy
plus affine terms in $q$, so $K(t_0,\cdot)$ is concave on $(0,1]$:
$\partial_q^2K(t_0,q)=-1/(q(1-q))<0$ for $0<q<1$.
Therefore
\begin{align*}
 K(t_0,q)
 \ge\frac{q_0-q}{q_0-\bar q}K(t_0,\bar q)
     +\frac{q-\bar q}{q_0-\bar q}K(t_0,q_0)
 \ge\frac{q_0-q}{q_0-\bar q}K(t_0,\bar q)>0.
\end{align*}
The last inequality uses $q<q_0$ and $K(t_0,\bar q)>0$. If $D<t_0\sqrt{\bar q}$, use $q=2\Phi(-t)$ and vary $t$.
The function $t\mapsto t\sqrt{2\Phi(-t)}$ is continuous on $[0,t_0]$.
It equals zero at $t=0$ and $t_0\sqrt{\bar q}>D$ at $t=t_0$.
The intermediate value theorem gives $t_1\in(0,t_0)$ such that
\(
 D=t_1\sqrt{2\Phi(-t_1)}.
\)
Taking $t=t_1$ and $q=2\Phi(-t_1)$ gives $D=t\sqrt q$ and
$K(t,q)>0$. In either case, \eqref{eq:GS-complexity-simple} yields
\(
 \Sigma^{\mathrm{GS}}(D,0,0)\ge K(D/\sqrt q,q)=K(t,q)>0.
\)
\end{proof}

\begin{lem}\label{lem:GS-finite-threshold}
We define the following threshold:
\begin{equation}
\label{eq:GS-annealed-critical-definition}
 D_c=\sup\left\{
 D_0>0:\ \exists\,q_0\in(0,1]\text{ with }
 K(D_0/\sqrt{q_0},q_0)\ge0\right\}.
\end{equation}
Then $0<D_c<\infty$.
\end{lem}

\begin{proof}
For the lower bound, take $q_0=2\Phi(-1)$ and $D_0=\sqrt{q_0}>0$.
As shown in the previous lemma, we have
\[
 K(D_0/\sqrt{q_0},q_0)=K(1,2\Phi(-1))>0.
\]
Thus this $D_0$ belongs to the set in
\eqref{eq:GS-annealed-critical-definition}, and
$D_c\ge\sqrt{2\Phi(-1)}>0$.

For the upper bound, $\Phi(-x)\le e^{-x^2/2}$ for $x\ge0$ gives,
for every $t>0$ and $u\ge0$,
\[
 \log(2\Phi(u-t))-\frac{u^2}{2}
 \le\log2-\frac{u^2+(\max\{t-u,0\})^2}{2}
 \le\log2-\frac{t^2}{4}.
\]
Here $u^2+(\max\{t-u,0\})^2\ge t^2/2$, with equality at $u=t/2$.
Since $\log(2\Phi(t)-1)\le0$ and $h_2(q)+q\log2\le\log3$,
\begin{equation}
\label{eq:GS-zero-field-upper-bound}
 K(t,q)\le h_2(q)+q\log2-\frac{qt^2}{4}
          \le\log3-\frac{qt^2}{4}.
\end{equation}
Here $h_2(q)+q\log2$ is the entropy of the three-point profile
$(q/2,1-q,q/2)$, whose entropy is at most $\log3$.
For every $D_0,q_0$ satisfying the condition in
\eqref{eq:GS-annealed-critical-definition}, substituting
$t=D_0/\sqrt{q_0}$ gives
\[
 0\le K(D_0/\sqrt{q_0},q_0)\le\log3-\frac{D_0^2}{4},
 \qquad D_0\le2\sqrt{\log3}.
\]
Taking the supremum over exactly those $D_0$ proves
$D_c\le2\sqrt{\log3}<\infty$.
\end{proof}

\begin{lem}
\label{lem:GS-vanishing-above}
For every $D\ge D_c$ and $q\in(0,1]$,
\[
 K(D/\sqrt q,q)\le0.
\]
Consequently, $\Sigma^{\mathrm{GS}}(D,0,0)=0$ for every
$D\ge D_c$.
\end{lem}

\begin{proof}
Fix $D\ge D_c$ and suppose $K(D/\sqrt q,q)>0$
for some $q\in(0,1]$. By Corollary~\ref{cor:p-spin-analysis}, the
inner maximum in \eqref{eq:GS-rescaled-objective} equals
$\Sigma^{\mathrm{mix}}(2,0,2t)$ and is continuous in $t>0$.
Thus, for this fixed $q$, the defining formula for $K$ makes
$D'\mapsto K(D'/\sqrt q,q)$ continuous on $(0,\infty)$.
There is therefore $\varepsilon>0$ such that
\[
 K((D+\varepsilon)/\sqrt q,q)>0.
\]
Thus $D+\varepsilon$ belongs to the set in
\eqref{eq:GS-annealed-critical-definition}, and
\(
 D+\varepsilon\le D_c\le D,
\)
a contradiction. Therefore every $K(D/\sqrt q,q)$ is nonpositive.
Since $D\ge D_c>0$, \eqref{eq:GS-complexity-simple} applies and gives
\[
 \Sigma^{\mathrm{GS}}(D,0,0)
 =\max\left\{0,\sup_{0<q\le1}K(D/\sqrt q,q)\right\}=0.
\]
This includes $D=D_c$.
\end{proof}

\begin{proof}[Proof of Proposition~\ref{prop:GS-topo-trivial}]
Lemma~\ref{lem:GS-finite-threshold} gives $0<D_c<\infty$, and Lemma~\ref{lem:GS-vanishing-above}
proves part~(ii). For part~(i), Lemma~\ref{lem:GS-reduction-K} covers $D\le0$.
Fix $0<D<D_c$. By the definition of $D_c$ in
\eqref{eq:GS-annealed-critical-definition}, there exist $D_0>D$ and
$q_0\in(0,1]$ such that
\[
 K(D_0/\sqrt{q_0},q_0)\ge0.
\]
Otherwise every $D_0$ in the defining set would be at most $D$, forcing
$D_c\le D$. Lemma~\ref{lem:GS-positivity-below-pair}
now gives $\Sigma^{\mathrm{GS}}(D,0,0)>0$, completing part~(i).
\end{proof}


\section{Application IV: Complexity at a given energy level}
\label{sec:energy-constrained}


\subsection{Setting and main results}

Although our general result, Theorem \ref{thm:general}, is written to count all critical points, in this section we explain how to extend it to count critical points lying in certain energy bands, which are frequently studied in the physics literature.

Consider the case of pure $p$-spin in zero potential, in which case the function from \eqref{eqn:p-spin-simplified-formula} reduces to
\begin{equation}
\label{eqn:energy-constrained-complexity-beginning}
\begin{split}
    \Sigma^{\textup{mix}}(p,0,\gamma) &= \sup_{z \geq 0} \left\{ \log\left(2\Phi\left( \sqrt{p-1} z - \frac{\gamma}{2}\sqrt{(p-1)!} \right) \right) - \frac{z^2}{2} \right\} \\
    &= \sup_{z \geq 0} \left\{ \log\left(2\Phi\left( \sqrt{\frac{b}{a}} z - \frac{\gamma}{2\sqrt{a}}\right)\right) - \frac{z^2}{2}\right\}.
\end{split}
\end{equation}
The function being optimized in \eqref{eqn:energy-constrained-complexity-beginning} is strictly concave on $[0,\infty)$ and increasing near zero, and thus the supremum is achieved at a unique $z_\ast = z_\ast(p,\gamma) \in (0,\infty)$. An important threshold for us is the rescaling
\begin{equation}
\label{eqn:e-infty-definition}
    E_\infty = E_\infty(p,\gamma) = \frac{1}{\sqrt{(p-1) \cdot (p-1)!}} z_\ast(p,\gamma) \in (0,\infty).
\end{equation}

We will need the function
\begin{equation}
\label{eqn:energy-constrained-sigma-p}
    \Sigma_p(E,\gamma) = \sup_{u \geq 0} \inf_{t \geq 0} \left\{ \log\left(2\Phi\left( t + \sqrt{\frac{b}{a}} u - \frac{\gamma}{2\sqrt{a}}\right)\right) - \frac{u^2}{2} + \frac{t^2}{2} - t\left( \frac{pE - \sqrt{b}u}{\sqrt{a}}\right) \right\}.
\end{equation}

Notice that, by taking $t = 0$, we find
\[
    \Sigma_p(E,\gamma) \leq \Sigma^{\textup{mix}}(p,0,\gamma).
\]

\begin{thm}
\label{thm:energy-constrained}
For pure $p$-spin in zero potential and $E \geq 0$, we have
\begin{equation}
\label{eqn:thm-energy-constrained}
    \lim_{N \to \infty} \frac{1}{N} \log \E\left[\#\left\{ \sigma : \sigma \text{ is $\gamma$-gapped for $H_N$ and } \frac{H_N(\sigma)}{N} \leq -E\right\} \right] = \Sigma_p(E,\gamma).
\end{equation}
\end{thm}

\begin{lem}
\label{lem:e-infinity}
We have 
\begin{equation}
\label{eqn:energy-constrained-complexity}
    \Sigma_p(E,\gamma) \begin{cases} = \Sigma^{\textup{mix}}(p,0,\gamma) & \text{if } 0 \leq E \leq E_\infty(p,\gamma), \\ < \Sigma^{\textup{mix}}(p,0,\gamma) & \text{if } E > E_\infty(p,\gamma). \end{cases}
\end{equation}
The function $\Sigma_p(\cdot,\gamma)$ is continuous, and strictly decreasing for $E > E_\infty(p,\gamma)$. Additionally, when $E > E_\infty(p,\gamma)$, it can alternatively be written as
\begin{equation}
\label{eqn:energy-constrained-explicit}
    \Sigma_p(E,\gamma) = \log(2\Phi(x)) + \frac{1}{2}\left(\frac{\phi(x)}{\Phi(x)}\right)^2 - \frac{p!E^2}{2},
\end{equation}
where $x$ is the unique solution to 
\begin{equation}
\label{eqn:energy-constrained-implicit}
    x + \frac{\phi(x)}{\Phi(x)} = \frac{pE-\gamma/2}{\sqrt{a}}.
\end{equation}
\end{lem}

\begin{rem}
The phenomenon that we find here -- namely that the annealed complexity function for sublevel sets is flat up to a point and decreasing thereafter -- is well-documented in spherical spin glasses (see, e.g., \cite{AufBenCer2013,McK2024}), and indeed the notation $E_\infty$ comes from this literature.

In spherical spin glasses, the threshold $-E_\infty$ is known to correspond to algorithmic phase transitions (\cite{HuaSel2023,HuaSel2024,Sel2024}, among others), in the sense that local algorithms seeking to descend in energy typically get trapped around energy $-E_\infty$. In our models, the connection to algorithms is less clear.

We note also that, for local minima of the SK model, one can numerically estimate $E_\infty(2,0) = z_\ast(2,0) \approx 0.506$, which has previously appeared many times in the literature (see, for example, \cite{AddarioBerryEtAl2019,TanEdw1980}). 
\end{rem}

\begin{rem}
\label{rem:BraMoo1981shortranged}
In the case of the SK model, when $E > E_\infty(2,0)$, the formula \eqref{eqn:energy-constrained-explicit} simplifies to 
\[
    \Sigma_2(E,0) = \log(2\Phi(x)) + E^2 - 2Ex + \frac{x^2}{2},
\]
where $x$ solves $x+\frac{\phi(x)}{\Phi(x)} = 2E$. This formulation previously appeared in work of Bray and Moore in the physics literature \cite{BraMoo1981shortranged} (precisely, by taking the first line of their (3.11) with $u = -(\epsilon+\tau)$, per their (3.8), one finds $g(\epsilon) = \log(2\Phi(\tau))+\epsilon^2+2\epsilon \tau + \frac{\tau^2}{2}$, where $\tau$ solves $-2\epsilon = \tau + \frac{\phi(\tau)}{\Phi(\tau)}$; setting $-2\epsilon = 2E$ and $\tau = x$ matches the formulas). 
\end{rem}

\begin{rem}
\label{rem:deOFon1997}
The special case $\gamma = 0$ of our result, in the regime $E > E_\infty$, recovers results of de Oliveira and Fontanari in the physics literature \cite{deOFon1997} (although they also consider the case of nonzero external field, which we do not). Indeed, our formula can be rewritten in this case as
\[
    \Sigma^{\textup{MW}}(E,0) = \log(2\Phi(x)) + \frac{1}{2}\left(x-\frac{p}{\sqrt{a}}E\right)^2 - \frac{p!E^2}{2},
\]
where $x$ solves $x + \frac{\phi(x)}{\Phi(x)} = \frac{pE}{\sqrt{a}}$. If we write our $p$-spin Hamiltonian as $H^{MW}_N$, they consider the Hamiltonian $H^{\textup{dOF}}_N = -\sqrt{\frac{p!}{2}} H^{\textup{MW}}_N$; thus we should change variables as $E = -\sqrt{\frac{2}{p!}} \epsilon$, in other words $\frac{p}{\sqrt{a}}E = -\sqrt{2p}\epsilon$, in which our formula becomes
\[
    \Sigma^{\textup{MW}}(\epsilon,0) = \log(2\Phi(x)) + \frac{1}{2}(x + \sqrt{2p}\epsilon)^2 - \epsilon^2 = \log(2\Phi(x)) + \frac{x^2}{2} + \sqrt{2p}\epsilon x + \epsilon^2(p-1).
\]
De Oliveira and Fontanari find that the complexity of critical points of their Hamiltonian near energy levels $N\epsilon$ is
\begin{equation}
\label{eqn:sigma-dof}
    \Sigma^{\textup{dOF}}(\epsilon) = \frac{\epsilon \nu}{\sqrt{p}} - \frac{1}{p-1} \left( \mu^2 - \mu \nu + \frac{\nu^2}{4p}\right) + \log(\erfc(-\mu)),
\end{equation}
(see their eqs. (16) and (17)) where $\mu$ and $\nu$ are saddle-point parameters found by taking the $\nu$- and $\mu$-partials of the right-hand side and setting them simultaneously to zero. We have $\log(\erfc(-\mu)) = \log(2\Phi(\eta))$, where we introduced the notation $\eta = \mu\sqrt{2}$ (notice that we can thus take partials in $\nu$ and $\eta$). Using this, and the fact that \eqref{eqn:sigma-dof} is quadratic in $\nu$ so that $\nu$ can be found easily, we see that
\[
    \Sigma^{\textup{dOF}}(\epsilon) = \log(2\Phi(\eta)) + \frac{\eta^2}{2} + \sqrt{2p}\epsilon\eta + \epsilon^2(p-1),
\]
where $\eta + \frac{\phi(\eta)}{\Phi(\eta)} + \sqrt{2p}\epsilon = 0$. Thus $\eta = x$ and the formulas match.
\end{rem}

\begin{rem}
\label{rem:AddarioBerryEtAl2019}
In the mathematics literature, Addario-Berry \emph{et al.} \cite{AddarioBerryEtAl2019} gave a formula for the energy-constrained complexity of the SK model with $\gamma = 0$, and with minimal work one can extend their methods to consider the pure $p$-spin model with $\gamma = 0$. This special case of our result matches theirs; see Section \ref{subsubsec:addario-berry-energy-constrained} for details.
\end{rem}


\subsection{Proof sketch}

We first sketch the proof informally, focusing especially on where $E_\infty(p,\gamma)$ comes from. The formal proof is essentially routine given this sketch; thus we postpone it to Appendix \ref{appx:energy-constrained}. To keep the paper at a manageable length, in the formal proof we do not prove all the formulas used in this proof sketch, but only weaker versions of them which suffice to prove the result, although we see no obstacle to proving them fully.

The fundamental observation is simply the computation
\begin{equation}
\label{eqn:sum-of-the-local-fields}
    -\sum_{i=1}^N G_i(\sigma) = pH_N(\sigma).
\end{equation}
This allows us to write our probability of interest just in terms of the local fields, namely as
\begin{align*}
    &\P\left(\sigma \text{ is $\gamma$-gapped for $H_N$ and } \frac{H_N(\sigma)}{N} \leq -E \right) \\
    &= \P\left(G_i(\sigma) \geq \frac{\gamma}{2} \text{ for all $i$ and } \frac{1}{N} \sum_{i=1}^N G_i(\sigma) \geq pE \right) \\
    &= \E_U\left[ \E_{X_1,\ldots,X_N}\left[ \prod_{i=1}^N \one\{\sqrt{a_N}X_i + \sqrt{Nb_N}\frac{U}{\sqrt{N}} \geq \frac{\gamma}{2}\} \right. \right. \\
    &\qquad \qquad \qquad \qquad \qquad \times \left. \left. \one\{ \sqrt{a_N} \left(\frac{1}{N} \sum_{i=1}^N X_i\right) + \sqrt{Nb_N} \frac{U}{\sqrt{N}} \geq pE \} \right] \right] \\
    &= \E_U\left[ p_N\left(\frac{U}{\sqrt{N}}, E \right) \Phi\left( \sqrt{\frac{Nb_N}{a_N}} \frac{U}{\sqrt{N}} - \frac{\gamma}{2\sqrt{a_N}} \right)^N \right]
\end{align*}
where
\[
    p_N(u, E) = \P\left( \left. \sqrt{a_N} \left(\frac{1}{N} \sum_{i=1}^N X_i\right) + \sqrt{Nb_N} u \geq pE \right| \sqrt{a_N} X_i + \sqrt{Nb_N} u \geq \frac{\gamma}{2} \text{ for all } i \right).
\]
(We use $a$ and $b$ throughout just for notational brevity, but we emphasize that \eqref{eqn:sum-of-the-local-fields} does not hold in the mixed case, and thus it is not clear how to compute energy-constrained complexity for mixed spin glasses with our methods.) Thus we need to understand
\[
    \lim_{N \to \infty} \frac{1}{N} \log \E_U\left[ \exp\left\{ N \Lambda_N\left(\frac{U}{\sqrt{N}} \right) \right\} \right]
\]
where
\[
    \Lambda_N(u) = \frac{1}{N} \log p_N(u,E) + \log \Phi\left(\sqrt{\frac{Nb_N}{a_N}} u - \frac{\gamma}{2\sqrt{a_N}} \right).
\]
We know how to handle the second term, from Section \ref{sec:p-spin}. The difficulty is the first term. Suppose that $p_N(u)$ is replaced with
\[
    \widetilde{p_N}(u,E) = \P\left( \left. \sqrt{a} \left(\frac{1}{N} \sum_{i=1}^N X_i\right) + \sqrt{b} u \geq pE \right| \sqrt{a} X_i + \sqrt{b} u \geq \frac{\gamma}{2} \text{ for all } i \right).
\]
This probability considers the empirical mean of a collection of i.i.d. random variables, whose distribution is truncated Gaussian; thus it is either close to one, or exponentially small as described by Cram\'er's theorem, according to whether our lower bound on the empirical mean is above or below its expectation. Explicitly, if $Y$ has the distribution of a standard Gaussian variable conditioned to lie above some threshold $-s$, then
\begin{equation}
\label{eqn:truncated-gaussian-basics}
\begin{split}
    \E[Y] &= \frac{\phi(s)}{\Phi(s)}, \\
    \E[e^{tY}] &= e^{\frac{t^2}{2}} \frac{\Phi(t+s)}{\Phi(s)}
\end{split}
\end{equation}
It will be useful to have the notations
\[
    s = s(u) = \frac{\sqrt{b}u - \gamma/2}{\sqrt{a}}, \qquad s_N = s_N(u) = \frac{\sqrt{Nb_N}u - \gamma/2}{\sqrt{a_N}}.
\]
If we let $\widetilde{Y_i}$ have the distribution of standard Gaussian conditioned to lie above $-s$, then \eqref{eqn:truncated-gaussian-basics} tells us that
\begin{align*}
    \E\left[ \sqrt{a} \left( \frac{1}{N} \sum_{i=1}^N \widetilde{Y_i} \right) + \sqrt{b}u\right] &= \sqrt{a} \frac{\phi(s)}{\Phi(s)} + \sqrt{b}u = \sqrt{a}g(s) + \frac{\gamma}{2},
\end{align*}
where 
\[
    g(x) = x + \frac{\phi(x)}{\Phi(x)}.
\]
This function is strictly increasing on $\R$ with $\lim_{x \to -\infty} g(x) = 0$ and $\lim_{x \to +\infty} g(x) = +\infty$. If this expectation is strictly bigger then $pE$, then we expect $\widetilde{p_N}(u) \approx 1$ from the weak law of large numbers. Otherwise, equivalently if 
\begin{equation}
\label{eqn:lln-boundary}
    u < \frac{1}{\sqrt{b}} \left[ \sqrt{a} g^{-1}\left( \frac{1}{\sqrt{a}} \left( pE - \frac{\gamma}{2}\right)\right) + \frac{\gamma}{2} \right], 
\end{equation}
Cram\'er's theorem and \eqref{eqn:truncated-gaussian-basics} tell us that $\frac{1}{N} \sum_{i=1}^N \widetilde{Y_i}$ satisfies an LDP with the good rate function
\[
    I(x) = \sup_{t \in \R} \left\{ tx - \frac{t^2}{2} - \log\Phi(t+s) + \log\Phi(s)\right\}
\]
(and for deviations above the mean, as here, we can restrict to $t \geq 0$), and thus in this case we expect
\[
    \frac{1}{N} \log p_N(u,E) \approx \frac{1}{N} \log \widetilde{p_N}(u,E) \approx -I\left( \frac{pE - \sqrt{b}u}{\sqrt{a}}\right).
\]
Putting this together, we imagine that
\[
    \Lambda_N(u) \approx -J(u,E,\gamma) + \log\Phi\left(s(u)\right),
\]
where
\begin{equation}
\label{eqn:energy-constrained-j-formula}
    J(u,E,\gamma) = \begin{cases} 0 & \text{if } u > \frac{1}{\sqrt{b}} \left[ \sqrt{a} g^{-1}\left( \frac{1}{\sqrt{a}} \left( pE - \frac{\gamma}{2}\right)\right) + \frac{\gamma}{2} \right], \\ I\left( \frac{pE-\sqrt{b}u}{\sqrt{a}}\right) & \text{otherwise.} \end{cases}
\end{equation}
Altogether, this suggests that
\[
    \Sigma_p(E,\gamma) = \sup_{u \in \R} \left\{ \log\Phi\left(s(u)\right) - J(u,E,\gamma) - \frac{u^2}{2} \right\}.
\]
In the regime where the $J$ term is nonzero, this has the form
\[
    \log\Phi\left(s(u)\right) - \frac{u^2}{2} + \inf_{t \geq 0} \left\{ \frac{t^2}{2} - t\left( \frac{pE - \sqrt{b}u}{\sqrt{a}} \right) + \log\Phi(t+s(u)) - \log(\Phi(s(u)))\right\}.
\]
Notice that the $\log(\Phi(s(u)))$ terms cancel, which suggests the formula \eqref{eqn:energy-constrained-sigma-p}. The question is whether the $J$ term is ``active.'' If the optimizer $z_\ast$ of this formula \emph{without} $J(u,E,\gamma)$ satisfies $J(z_\ast,E,\gamma) = 0$, then clearly $z_\ast$ will still be the optimizer of the formula with $J$, since $-J$ is nonpositive; otherwise, one naturally guesses (and it turns out to be true \emph{a posteriori}) that the formula with $J$ has a new optimizer. Since $z_\ast$ is the optimizer of \eqref{eqn:energy-constrained-complexity-beginning}, it satisfies the first-order condition
\begin{equation}
\label{eqn:z-ast-first-order}
    z_\ast = \sqrt{\frac{b}{a}} \frac{\phi\left(\sqrt{\frac{b}{a}} z_\ast - \frac{\gamma}{2\sqrt{a}}\right)}{\Phi\left(\sqrt{\frac{b}{a}} z_\ast - \frac{\gamma}{2\sqrt{a}}\right)}.
\end{equation}
Comparing to \eqref{eqn:energy-constrained-j-formula}, we see that $J(z_\ast,E,\gamma)$ vanishes when
\begin{equation}
\label{eqn:e-infinity-manipulation-1}
    \frac{1}{\sqrt{a}} (pE - \frac{\gamma}{2}) < g\left( \sqrt{\frac{b}{a}}z_\ast - \frac{\gamma}{2\sqrt{a}}\right) = \sqrt{\frac{a}{b}} z_\ast + \sqrt{\frac{b}{a}} z_\ast - \frac{\gamma}{2\sqrt{a}} 
\end{equation}
or in other words when
\begin{equation}
\label{eqn:e-infinity-manipulation-2}
    E < \frac{\sqrt{a}}{p} \left( \sqrt{\frac{a}{b}} + \sqrt{\frac{b}{a}} \right)z_\ast = E_\infty,
\end{equation}
where we finally substituted the values of $a$ and $b$. 


\appendix 


\section{Comparison to existing methods}
\label{sec:comparison}

In \cite{AddarioBerryEtAl2019}, Addario-Berry, Devroye, Lugosi, and Oliveira introduced a very nice method that reduces counting local minima to large deviations for the $L^1$ norm of a standard Gaussian vector. Their paper was written only for the SK model with zero potential, but their methods extend immediately to mixed $p$-spin models with zero potential. With non-trivial additional work, one can also add an external field and study gapped states (and understand disorder persistence, although in the interest of space we do not pursue this in details). However, it seems to us that the resulting complexity formulas (which are different variational problems than the ones we find, although \emph{a posteriori} equal) are typically higher-dimensional, so that it is not clear to us how to use them to establish qualitative results like those in Corollary \ref{cor:p-spin-analysis}. 


\subsection{The method and immediate generalizations}
Here is a rough sketch of their method, translated into our notation and scaling conventions, and generalized to mixed $p$-spin glasses with zero potential (i.e., $V = 0$). Throughout, we use $\approx$ informally to denote things like dropping subexponential factors, replacing $N-1$ with $N$, and replacing $a_N$ and $b_N$ from \eqref{eq:mixed-cov} with their limits.

The vector $G(\sigma)$ defined in \eqref{eq:mixed-G} is centered Gaussian with covariance matrix $C = a_N I_N + b_N \one\one^{\mathsf T}$, which has inverse 
\begin{equation}
\label{eqn:addario-berry-c-inverse}
    C^{-1} = \frac{1}{a_N} \left( I_N - \left( \frac{b_N}{a_N + Nb_N} \right) \one \one^{\mathsf T} \right) \approx \frac{1}{a} \left( I_N - \frac{1}{N} \frac{b}{a+b} \one \one^{\mathsf T} \right).
\end{equation}
Thus, as in \eqref{eq:mixed-rep} in our argument,
\begin{align*}
    \P(\sigma \text{ is a local minimum}) &= \P(G_i \geq 0 \text{ for all } i) \\
    &= \frac{1}{(2\pi)^{N/2} \det(C)^{1/2}} \int_{[0,\infty)^N} \exp\left(-\frac{1}{2} \ip{x,C^{-1}x} \right) \diff x \\
    &\approx \frac{1}{(2\pi a)^{N/2}} \int_{[0,\infty)^N} \exp\left( -\frac{\|x\|_2^2}{2a} + \frac{b}{2a(a+b)} \frac{\|x\|_1^2}{N} \right) \diff x \\
    &= \frac{1}{(2\pi)^{N/2}} \int_{[0,\infty)^N} \exp\left( -\frac{\|x\|_2^2}{2} + \frac{b}{2(a+b)} \frac{\|x\|_1^2}{N} \right) \diff x \\
    &= 2^{-N} \frac{1}{(2\pi)^{N/2}} \int_{\R^N} \exp\left( -\frac{\|x\|_2^2}{2} + \frac{b}{2(a+b)} \frac{\|x\|_1^2}{N} \right) \diff x \\
    &= 2^{-N} \E\left[ \exp\left(N \cdot \frac{b}{2(a+b)} Y_N^2\right)\right],
\end{align*}
where $Y \in \R^N$ is a standard Gaussian vector and $Y_N \defeq \|Y\|_1/N$. The variable $Y_N$ is the average of $N$ i.i.d. random variables, so Cram\'er's theorem implies that it satisfies an LDP at speed $N$ with the good rate function
\begin{equation}
\label{eqn:mu-ast-definition}
    \mu^\ast(x) = \begin{cases} +\infty & \text{if } x < \sqrt{2/\pi}, \\ \sup_{\lambda \geq 0} \{\lambda x - \lambda^2/2 - \log(2\Phi(\lambda))\} & \text{if } x \geq \sqrt{2/\pi}. \end{cases}
\end{equation}
Thus Varadhan's lemma suggests that we should have
\begin{equation}
\label{eqn:addario-berry-varadhan}
    \lim_{N \to \infty} \frac{1}{N} \log \E \mc{N}_N^{\textup{mix}}(\xi,0,0) = \sup_{y \geq \sqrt{2/\pi}} \left\{ \frac{b}{2(a+b)} y^2 - \mu^\ast(y) \right\}.
\end{equation}
(The function $y \mapsto \frac{b}{a+b} y^2$ is unbounded, so classical Varadhan's lemma does not apply, but \cite{AddarioBerryEtAl2019} does hands-on estimates to show that \eqref{eqn:addario-berry-varadhan} nevertheless holds for the SK model.) A posteriori, this must match our formla \eqref{eqn:p-spin-simplified-formula}; for the special case of $p$-spin, the argument in Section \ref{subsubsec:addario-berry-energy-constrained} below shows that this can also be seen analytically. 


\subsection{External field}
So far, our argument essentially follows that of \cite{AddarioBerryEtAl2019}; generalizing from the SK model to mixed spin glasses involves only cosmetic parameter changes. However, one can still adapt their methods to count local minima for a mixed spin glass with external field $V(m) = hm$, although it involves more work. Since our purpose here is just to point out that the resulting variational problems seem to be complicated, and in the interest of space, we give only heuristics rather than full proofs (although we see no obstruction to rigorously proving formula \eqref{eqn:addario-berry-sk-external-field} below). For the same reason, we restrict to the case of the SK model ($\xi(q) = q^2/2$), although one could extend to mixtures just by changing the parameters (as in the case $h = 0$), and one could count gapped states by putting $\gamma$ parameters everywhere.

Indeed, comparing to \eqref{eq:mixed-increment}, we see that we need to compute
\begin{align*}
    \P(\sigma \text{ is a local minimum}) &= \P(G_i \geq h\sigma_i \text{ for each } i).
\end{align*}
We can think about this as entrywise nonnegativity for the non-centered Gaussian vector $\widetilde{G}(\sigma) = G(\sigma) - h\sigma$; in this way, we find
\begin{align*}
    &\P(\sigma \text{ is a local minimum}) \\
    &= \frac{1}{(2\pi)^{N/2} \det(C)^{1/2}} \int_{[0,\infty)^N} \exp\left( - \frac{1}{2} \ip{x + h\sigma, C^{-1} \left(x + h\sigma\right) } \right) \diff x \\
    &\approx \frac{1}{(2\pi)^{N/2}} \exp\left(-\frac{h^2}{2} \ip{\sigma, C^{-1}\sigma} \right) \int_{[0,\infty)^N} \exp\left(-\frac{1}{2}\ip{x,C^{-1}x} - h\ip{x,C^{-1}\sigma} \right) \diff x.
\end{align*}
From the display, one can see that $\P(\sigma \text{ is a local minimum})$ only depends on the number of plus and minus ones in $\sigma$, not their ordering. Thus, for $t \in [0,1]$, write $\sigma_t$ for the point in $\{-1,+1\}^N$ with plus ones in the first $tN$ positions and minus ones in the last $(1-t)N$ positions (assume these are integers). If we compute
\[
    S(t) = \lim_{N \to \infty} \frac{1}{N} \log \P(\sigma_t \text{ is a local minimum}),
\]
then since the number of $\sigma$'s with $tN$ plus ones is $\binom{N}{tN} \approx \exp(-I(t))$ with 
\begin{equation}
\label{eqn:addario-berry-I}
    I(t) = t\log t + (1-t)\log(1-t),
\end{equation}
we would expect to find
\begin{equation}
\label{eqn:addario-berry-sk-external-field}
    \lim_{N \to \infty} \frac{1}{N} \log \E \mc{N}_N = \sup_{t \in [0,1]} \{S(t) - I(t)\}.
\end{equation}
Thus it remains only to continue our computation of $\P(\sigma_t \text{ is a local minimum})$. Since $\ip{\sigma_t,\one} = (2t-1)N$, the approximation \eqref{eqn:addario-berry-c-inverse} (with $a = b = 1$ for the SK model) gives
\[
    \ip{\sigma,C^{-1}\sigma} \approx N \left( 1 - \frac{(2t-1)^2}{2} \right).
\]
To estimate the integrand, think of $x \in \R^N$ as the concatenation $x = (x_1,x_2)$, where $x_1 \in \R^{tN}$ stores the first $tN$ coordinates and $x_2$ stores the last $(1-t)N$ coordinates. Then, for $x \in [0,\infty)^N$, we have
\[
    \ip{x,C^{-1}\sigma} \approx \|x_1\|_1 - \|x_2\|_1 - \frac{2t-1}{2} (\|x_1\|_1 + \|x_2\|_1).
\]
Combining, we find that
\begin{align*}
    &\P(\sigma_t \text{ is a local minimum}) \\
    &\approx \frac{2^{-N}}{(2\pi)^{N/2}} \exp\left( - \frac{h^2N}{2} \left(1 - \frac{(2t-1)^2}{2}\right) \right) \\
    &\times \int_{\R^N} \exp\left(-\frac{\|x\|_2^2}{2} + \frac{(\|x_1\|_1 + \|x_2\|_1)^2}{4N} - h\left(\|x_1\|_1 - \|x_2\|_1 - \frac{2t-1}{2}(\|x_1\|_1 + \|x_2\|_1) \right) \right) \diff x \\
    &= 2^{-N} \exp\left( - \frac{h^2N}{2} \left(1 - \frac{(2t-1)^2}{2}\right) \right) \E\left[\exp\left( N \phi_{h,t}(Y^{(1)}_N,Y^{(2)}_N) \right) \right],
\end{align*}
where $Y$ is a standard Gaussian vector, $Y^{(1)}_N = \frac{1}{tN} \sum_{i=1}^{tN} \abs{Y_i}$, and $Y^{(2)}_N = \frac{1}{(1-t)N} \sum_{i=tN+1}^N \abs{Y_i}$, and
\[
    \phi_{h,t}(x,y) = \frac{(tx+(1-t)y)^2}{4} - h\left( tx-(1-t)y - \frac{2t-1}{2}(tx+(1-t)y)\right).
\]
This suggests that
\begin{equation}
\label{eqn:addario-berry-sk-S}
    S(t) = -\log(2) - \frac{h^2}{2} \left(1-\frac{(2t-1)^2}{2}\right) + \sup_{x,y \geq \sqrt{2/\pi}} \{\phi_{h,t}(x,y) - t\mu^\ast(x) - (1-t)\mu^\ast(y)\}.
\end{equation}
(The Varadhan's lemma corrections are $-t\mu^\ast(x)$ rather than $-\mu^\ast(x)$ because the ``natural speed'' for the LDP for $Y^{(1)}_N$ is $tN$, but the prefactor of $\phi_{h,t}$ is $N$, and similarly for $y$.) 

Thus the variational formula in \eqref{eqn:addario-berry-sk-external-field} is over either three or five parameters, depending on how one counts (it is over $t$, $x$, and $y$, plus the implicit variational problems defining $\mu^\ast(x)$ and $\mu^\ast(y)$), compared to the one-dimensional problem in \eqref{eqn:p-spin-simplified-formula}. It is not clear to us how to see, for example, that \eqref{eqn:addario-berry-sk-external-field} has the scaling given by \eqref{eqn:large-h-asymptotics} as $h \to +\infty$, nor even necessarily how to see that it is positive for all $h$.


\subsection{Energy-constrained complexity}
\label{subsubsec:addario-berry-energy-constrained}

One can also use the methods of \cite{AddarioBerryEtAl2019} to count local minima in certain energy bands. The authors of \cite{AddarioBerryEtAl2019} already noticed this for the SK model, but we can extend their arguments for pure $p$-spin glasses, $\xi(q) = q^p/p!$; by mimicking their arguments and putting $p$'s in the correct places, one would find
\begin{align*}
    &\lim_{N \to \infty} \frac{1}{N} \log \E[\#\{\sigma : \sigma \text{ is a local minimum and } H_N(\sigma) \leq -EN\}] \\
    &= \sup_{y \geq \max(\sqrt{2/\pi},p \sqrt{(p-1)!} \cdot E)} \left\{ \frac{p-1}{2p} y^2 - \mu^\ast(y)\right\} =: \Sigma^{ABDLO}(E).
\end{align*}
In this case, we are able to analytically match our result to theirs. Indeed, consider the definition \eqref{eqn:mu-ast-definition} of $\mu^\ast(y)$. Since 
\[
    \left. \frac{\partial}{\partial \lambda} \{\lambda y - \lambda^2/2 - \log(2\Phi(\lambda))\}\right|_{\lambda = 0} = y - \sqrt{2/\pi},
\]
for $y > \sqrt{2/\pi}$ we see that the optimizer $\lambda_\ast(x)$ in the definition of $\mu^\ast$ is in the interior, and thus that it satisfies the first-order conditions: We have
\[
    \mu^\ast(y) = \lambda_\ast(y) y - \lambda_\ast(y)^2/2 - \log(2\Phi(\lambda_\ast(y)))
\]
where $\lambda_\ast(y)$ satisfies
\begin{equation}
\label{eqn:lambda-ast-y}
    y = \lambda_\ast(y) + \frac{\phi(\lambda_\ast(y))}{\Phi(\lambda_\ast(y))}.
\end{equation}
Suppose we can show that the function
\[
    h_p(y) = \frac{p-1}{2p} y^2 - \mu^\ast(y)
\]
is strictly concave and achieves its maximum at $y = p\sqrt{(p-1)!} \, E_\infty$. Then, in the definition of $\Sigma^{ABDLO}(E)$, the supremum is achieved at $y = p\sqrt{(p-1)!} E \eqdef \widetilde{E}$ for $E > E_\infty$. If we let $x = \lambda_\ast(\widetilde{E})$, then we would have
\[
    \Sigma^{ABDLO}(E) = \frac{(p-1) \cdot p!}{2} E^2 - \lambda_\ast(\widetilde{E}) \widetilde{E} + \frac{\lambda_\ast( \widetilde{E})^2}{2} + \log\left(2\Phi\left(\lambda_\ast(\widetilde{E})\right)\right).
\]
On the other hand, writing our formula $\Sigma^{MW}$, we see that 
\[
    \Sigma^{MW}(E) = \log(2\Phi(x)) + \frac{1}{2}(\widetilde{E}-x)^2 - \frac{p!E^2}{2} = \frac{(p-1) \cdot p!}{2}E^2 - x\widetilde{E} + \frac{x^2}{2},
\]
where $x$ satisfies $x + \frac{\phi(x)}{\Phi(x)} = \widetilde{E}$. Matching $x = \lambda_\ast(\widetilde{E})$, we see that $\Sigma^{ABDLO}(E) = \Sigma^{MW}(E)$ for $E > E_\infty$. For $E \leq E_\infty$, the same claims about $h_p$ would show that $\Sigma^{ABDLO}(E)$ is constant for $0 \leq E \leq E_\infty$, and that it is right-continuous at $E_\infty$; since the same is true for $\Sigma^{MW}(E)$, and since $\Sigma^{MW}(E) = \Sigma^{ABDLO}(E)$ for $E > E_\infty$, we conclude that they agree everywhere.

Thus it remains only to check the claims about $h_p(y)$. Notice that
\[
    (\mu^\ast)'(y) = \lambda'_\ast(y)y + \lambda_\ast(y) - \lambda_\ast(y)\lambda_\ast'(y) - \frac{\phi(\lambda_\ast(y))}{\Phi(\lambda_\ast(y))} \lambda_\ast'(y) = \lambda_\ast(y),
\]
from the first order-conditions. Next, differentiating the implicit definition of $\lambda_\ast$, we find
\begin{equation}
\label{eqn:lambdaastprime}
    1 = \lambda_\ast'(y) \left[ 1 - \lambda_\ast(y) \frac{\phi(\lambda_\ast(y))}{\Phi(\lambda_\ast(y))} - \left(\frac{\phi(\lambda_\ast(y))}{\Phi(\lambda_\ast(y))}\right)^2 \right].
\end{equation}
Let $r(\lambda) = \phi(\lambda)/\Phi(\lambda)$. Now we claim that, for all $\lambda \geq 0$, we have $1 - \lambda r(\lambda) - r(\lambda)^2 > 0$. Indeed, if $Z$ has the distribution of a standard Gaussian variable conditioned to lie below $\lambda$, then we have
\[
    \E[e^{tZ}] = e^{\frac{t^2}{2}} \frac{\Phi(\lambda - t)}{\Phi(\lambda)},
\]
from which we have
\[
    \Var(Z) = 1 - \lambda r(\lambda) - r(\lambda)^2.
\]
This is strictly positive since $Z$ is not degenerate, and since $\lambda, r(\lambda) > 0$ we also have $\Var(Z) < 1$. Comparing to \eqref{eqn:lambdaastprime}, we see that $\lambda_\ast'(y) > 1$, and thus
\[
    \frac{\diff^2}{\diff y^2} \mu^\ast(y) = \lambda_\ast'(y) > 1
\]
for all $y$, so that the function $h_p(y) = \frac{p-1}{2p} y^2 - \mu^\ast(y)$ satisfies
\[
    h''_p(y) = \frac{p-1}{p} - \frac{\diff^2}{\diff y^2} \mu^\ast(y) < -\frac{1}{p}.
\]
Thus $h_p(y)$ is strictly convex, and it achieves its maximum when $y$ satisfies the first-order condition $\frac{p-1}{p} y = \lambda_\ast(y)$, or in other words -- comparing to \eqref{eqn:lambda-ast-y} -- when
\begin{equation}
\label{eqn:optimizing-y}
    y = \frac{p-1}{p} y + \frac{\phi(\frac{p-1}{p}y)}{\Phi(\frac{p-1}{p}y)}.
\end{equation}
When $y = y_\ast = p\sqrt{(p-1)!} \, E_\infty$, we have that $\frac{p-1}{p} y_\ast = (p-1) \sqrt{(p-1)!} \, E_\infty = \sqrt{p-1} \, z_\ast(p,0)$, in the notation of \eqref{eqn:e-infty-definition}; on the other hand, from the first-order condition \eqref{eqn:z-ast-first-order}, we see that $z_\ast(p,0)$ satisfies
\begin{equation}
\label{eqn:z-ast-first-order-differently}
    z_\ast(p,0) = \sqrt{p-1} \frac{\phi(\sqrt{p-1}z_\ast(p,0))}{\Phi(\sqrt{p-1}z_\ast(p,0))}.
\end{equation}
Comparing \eqref{eqn:e-infty-definition}, \eqref{eqn:z-ast-first-order-differently}, and \eqref{eqn:optimizing-y}, we see that $h_p(y)$ indeed achieves its maximum when $y = y_\ast = p\sqrt{(p-1)!} \, E_\infty$, which finishes the proof. 


\section{Proof details for energy-constrained complexity}
\label{appx:energy-constrained}

In this section, we prove the results of Section \ref{sec:energy-constrained}.

\begin{proof}
Now we prove these results formally.
\begin{itemize}
\item \textbf{Case 1 (Shallow energies, $E < E_\infty(p,\gamma)$):} Since $p_N$ is a probability, we have $\Lambda_N(u) \leq \log(\Phi(\sqrt{\frac{Nb_N}{a_N}} u - \frac{\gamma}{2\sqrt{a_N}}))$; then the upper bound in \eqref{eqn:energy-constrained-complexity} is exactly as in the proof of Theorem \ref{thm:complexity-p-spin} and Corollary \ref{cor:p-spin-analysis}. So we only need to understand the lower bound, which we can do by restricting $u$. Let $\delta > 0$ be as in Lemma \ref{lem:wlln-application} below; by lower-bounding $\Lambda_N$ on the event $\{U/\sqrt{N} \in [z_\ast-\delta,z_\ast+\delta]\}$, applying Lemma \ref{lem:wlln-application}, and then using Gaussian tail bounds, we find
\begin{align*}
    &\liminf_{N \to \infty} \frac{1}{N} \log \E_U\left[ \exp\left\{ N \Lambda_N\left(\frac{U}{\sqrt{N}}\right) \right\} \right] \\
    &\geq \liminf_{N \to \infty} \frac{1}{N} \log \E_U\left[ \exp\left\{ N \Lambda_N\left(\frac{U}{\sqrt{N}}\right) \right\} \mathds{1}_{\frac{U}{\sqrt{N}} \in [z_\ast-\delta,z_\ast+\delta]} \right] \\
    &\geq \liminf_{N \to \infty} \log \Phi\left( \sqrt{\frac{Nb_N}{a_N}} (z_\ast - \delta) - \frac{\gamma}{2\sqrt{a_N}}\right) + \liminf_{N \to \infty} \frac{1}{N} \log \P(U/\sqrt{N} \in [z_\ast-\delta,z_\ast+\delta]) \\
    &= \log\Phi\left( \sqrt{\frac{b}{a}} (z_\ast-\delta) - \frac{\gamma}{2\sqrt{a}} \right) - \frac{(z_\ast-\delta)^2}{2}.
\end{align*}
Letting $\delta \to 0$ shows that 
\begin{align*}
    &\liminf_{N \to \infty} \frac{1}{N} \log \E\left[\#\left\{ \sigma : \sigma \text{ is $\gamma$-gapped for $H_N$ and } \frac{H_N(\sigma)}{N} \leq -E\right\} \right] \\
    &\geq \log\left(2 \Phi\left(\sqrt{\frac{b}{a}} z_\ast - \frac{\gamma}{2\sqrt{a}}\right)\right) - \frac{z_\ast^2}{2} = \Sigma^{\textup{mix}}(p,0,\gamma).
\end{align*}
\item \textbf{Case 2 (Deep energies, $E > E_\infty(p,\gamma)$):} 
Standard Chernoff bounds yield
\[
    \Lambda_N(u) \leq \inf_{t \geq 0} \left\{ \log\Phi(t+s_N(u)) + \frac{t^2}{2} - t\left(\frac{pE-\sqrt{Nb_N}u}{\sqrt{a_N}} \right)\right\}.
\]
Let
\[
    \Lambda(u) = \inf_{t \geq 0} \left\{ \log\Phi(t+s(u)) + \frac{t^2}{2} - t\left(\frac{pE-\sqrt{b}u}{\sqrt{a}} \right)\right\}.
\]
Fix $u$, and let $t_\ast(u)$ be the optimizer in the definition of $\Lambda(u)$. By choosing $t = t_\ast(u)$, we find
\[
    \Lambda_N(u) \leq \log\Phi(t_\ast(u) + s_N(u)) + \frac{t_\ast(u)^2}{2} - t_\ast(u)\left(\frac{pE - \sqrt{Nb_N}u}{\sqrt{a_N}}\right) = \Lambda(u) + e_N(u),
\]
where
\begin{equation}
\label{eqn:energy-constrained-error}
    e_N(u) = \log\Phi(t_\ast(u) + s_N(u)) - \log\Phi(t_\ast(u) + s(u)) + t_\ast(u)\left( \frac{pE - \sqrt{Nb_N}u}{\sqrt{a_N}} - \frac{pE - \sqrt{b}u}{\sqrt{a}} \right).
\end{equation}
Notice also that, by taking $t = 0$, we have $\Lambda_N(u) \leq 0$; thus
\begin{align*}
    &\limsup_{N \to \infty} \frac{1}{N} \log \E_U\left[ \exp\left\{ N \Lambda_N\left(\frac{U}{\sqrt{N}} \right)\right\}\mathds{1}\left\{\abs{U/\sqrt{N}} > M\right\}\right] \\
    &\leq \limsup_{N \to \infty} \frac{1}{N} \log \P(\abs{U}/\sqrt{N} > M) = -\frac{M^2}{2}.
\end{align*}
Since $\Lambda(u) \leq 0$ as well, Varadhan's lemma gives
\begin{align*}
    \lim_{N \to \infty} \frac{1}{N} \log \E_U\left[ \exp\left\{ N\Lambda\left(\frac{U}{\sqrt{N}}\right) \right\} \right] &= \sup_{u \in \R} \left\{ \Lambda(u) - \frac{u^2}{2}\right\} = \sup_{u \geq 0} \left\{ \Lambda(u) - \frac{u^2}{2}\right\} \\
    &= \Sigma_p(E,\gamma) - \log(2),
\end{align*}
where we used that $\Lambda(u) \leq \Lambda(\abs{u})$. (The $-\log(2)$ factor is since we are here counting the probability of being an energy-constrained gapped state, not the number of such states.) Combining these two estimates with Lemma \ref{lem:energy-constrained-error} below, we find
\[
    \limsup_{N \to \infty} \frac{1}{N} \log \E_U\left[ \exp \left\{N \Lambda_N\left(\frac{U}{\sqrt{N}}\right)\right\} \right] \leq \max\left( - \frac{M^2}{2}, \Sigma_p(E,\gamma) - \log(2) \right)
\]
for any $M > 0$, and then take $M \to \infty$ to complete the upper bound.

Now we consider the lower bound. Since $E > E_\infty$ is fixed, the proof of Lemma \ref{lem:e-infinity} shows that the variational problem over $u$ is achieved at $\widetilde{u} = \widetilde{u}(E) = E\sqrt{b}/a$. With this choice, for each $\delta > 0$ we have
\[
    \E_U\left[ \exp\left\{ N \Lambda_N \left(\frac{U}{\sqrt{N}} \right) \right\} \right] \geq \E_U\left[ \exp\left\{ N \Lambda_N \left(\frac{U}{\sqrt{N}} \right) \right\} \mathds{1}_{\frac{U}{\sqrt{N}} \in [\widetilde{u}-\delta,\widetilde{u}+\delta]}\right]
\]
and thus, from Lemma \ref{lem:pointwise-cramer}, that
\begin{align*}
    &\liminf_{N \to \infty} \frac{1}{N} \log \E_U\left[ \exp\left\{ N \Lambda_N \left(\frac{U}{\sqrt{N}} \right) \right\} \right] \\
    &\geq \liminf_{N \to \infty} \frac{1}{N} \log \E_U\left[ \exp\left\{ N \left( -I\left(\frac{pE-\sqrt{b}\frac{U}{\sqrt{N}}}{\sqrt{a}}\right) + \log\Phi(s(U/\sqrt{N})) \right) \right\} \mathds{1}_{\frac{U}{\sqrt{N}} \in [\widetilde{u}-\delta,\widetilde{u}+\delta]} \right] \\
    &\qquad - \limsup_{N \to \infty} \sup_{u \in [\widetilde{u}-\delta,\widetilde{u}+\delta]} \abs{\log(\Phi(s_N(u))) - \log(\Phi(s(u)))}.
\end{align*}
By-now-familiar local-Lipschitz arguments show that the second term vanishes. For the first term, we note that
\[
    \liminf_{\delta \to 0} \inf_{u \in [\widetilde{u}-\delta,\widetilde{u}+\delta]} \left( - I \left( \frac{pE-\sqrt{b}u}{\sqrt{a}}\right) + \log\Phi(s(u)) \right) = - I \left( \frac{pE-\sqrt{b}\widetilde{u}}{\sqrt{a}}\right) + \log\Phi(s(\widetilde{u})).
\]
Gaussian tail bounds give $\lim_{\delta \to 0} \lim_{N \to \infty} \frac{1}{N} \log \P(U/\sqrt{N} \in [\widetilde{u}-\delta,\widetilde{u}+\delta]) = -\widetilde{u}^2/2$, so that (by taking $N \to \infty$ and then $\delta \to 0$) we find
\begin{align*}
    \liminf_{N \to \infty} \frac{1}{N} \log \E_U\left[ \exp\left\{ N \Lambda_N \left(\frac{U}{\sqrt{N}} \right) \right\} \right] &\geq -I\left(\frac{pE-\sqrt{b}\widetilde{u}}{\sqrt{a}}\right) + \log\Phi(s(\widetilde{u})) - \frac{\widetilde{u}^2}{2} \\
    &= \Sigma_p(E,\gamma)-\log(2),
\end{align*}
where we used Lemma \ref{lem:e-infinity} again.
\item \textbf{Case 3 (Critical energies, $E = E_\infty(p,\gamma)$):} This is cheap, given the previous two. The left-hand side of \eqref{eqn:thm-energy-constrained} is clearly non-increasing in $E$. From Lemma \ref{lem:e-infinity}, the right-hand side is also non-increasing in $E$, and is additionally continuous in $E$. Since the left- and right-hand sides agree for all $E \neq E_\infty$, these conditions force them to agree at $E = E_\infty$ as well.
\end{itemize}
\end{proof}

\begin{proof}[Proof of Lemma \ref{lem:e-infinity}]
Let
\begin{align*}
    k_E(u,t) &= \log(2\Phi\left( t + s(u)\right)) - \frac{u^2}{2} + \frac{t^2}{2} - t\left( \frac{pE - \sqrt{b}u}{\sqrt{a}}\right).
\end{align*}
We have
\begin{equation}
\label{eqn:first-order-in-t}
    \partial_t k_E(u,t) = g\left(t + s(u)\right) + \frac{\gamma}{2\sqrt{a}} - \frac{pE}{\sqrt{a}}.
\end{equation}
The right-hand side of \eqref{eqn:first-order-in-t} is strictly increasing in $t$ since $g$ is. Its value at $t = 0$ is $g(s(u)) + \frac{\gamma/2-pE}{\sqrt{a}}$. When $E \leq E_\infty$, by reversing the manipulations around \eqref{eqn:e-infinity-manipulation-1} and \eqref{eqn:e-infinity-manipulation-2} we see that $\partial_t k_E(z_\ast,0) \geq 0$; thus, when $E \leq E_\infty$, the function $\partial_t k_E(z_\ast,t)$ is nonnegative for all $t$, which implies
\[
    \Sigma_p(E,\gamma) \geq \inf_{t \geq 0} k_E(z_\ast,t) = k_E(z_\ast,0) = \Sigma^{\textup{mix}}(p,0,\gamma).
\]
This proves the first case of \eqref{eqn:energy-constrained-complexity} (recall that the corresponding upper bound is clear). In the remainder, we assume that $E > E_\infty$. 

We want to find a zero (in $u$) of $\partial_t k_E(u,t) = g(t+s(u)) + \frac{\gamma/2-pE}{\sqrt{a}}$. Since $g > 0$, this can only exist when $\gamma/2 - pE < 0$. So we want to show that $\frac{\gamma}{2p} < E_\infty$, or in other words (comparing to \eqref{eqn:e-infty-definition}) that $\frac{\sqrt{p-1}}{p} \cdot \frac{\gamma}{2\sqrt{a}} < z_\ast$. Let
\[
    h_E(u) = k_E(u,0) = \log(2\Phi(s(u))) - \frac{u^2}{2}.
\]
Since $\log \Phi$ is strictly concave and $s(\cdot)$ is linear, the function $h_E$ is strictly concave in $u$, and optimized at $z_\ast$, so one can show that $\frac{\sqrt{p-1}}{p} \cdot \frac{\gamma}{2\sqrt{a}} < z_\ast$ by showing that $h'_E(u) > 0$ when $u = \frac{\sqrt{p-1}}{p} \cdot \frac{\gamma}{2\sqrt{a}}$. Notice that
\[
    h'_E(u) = \sqrt{\frac{b}{a}} \left( g(s(u)) - \left(\sqrt{\frac{b}{a}} + \sqrt{\frac{a}{b}}\right) u + \frac{\gamma}{2\sqrt{a}}\right).
\]
One can compute $s(\frac{\sqrt{p-1}}{p} \cdot \frac{\gamma}{2\sqrt{a}}) = -\frac{\gamma}{2\sqrt{a}}$, and thus (after some algebra)
\begin{align*}
    h'_E\left(\frac{\sqrt{p-1}}{p} \cdot \frac{\gamma}{2\sqrt{a}}\right) = \sqrt{\frac{b}{a}} g \left( -\frac{\gamma}{2\sqrt{a}}\right) > 0.
\end{align*}
Thus indeed $E_\infty > \gamma/(2p)$, so that when $E \geq E_\infty$ we can safely pick some $u_\ast(E,\gamma)$ which satisfies
\begin{equation}
\label{eqn:u-ast-implicit}
    g(s(u_\ast(E,\gamma))) + \frac{\gamma/2-pE}{\sqrt{a}} = 0,
\end{equation}
which can also be written explicitly as
\begin{equation}
\label{eqn:u-ast}
    u_\ast(E,\gamma) = \sqrt{\frac{a}{b}} \left( g^{-1}\left( \frac{pE-\gamma/2}{\sqrt{a}} \right) + \frac{\gamma}{2\sqrt{a}} \right).
\end{equation}
To show that $u_\ast(E,\gamma) > 0$, we want to show that $g(s(0)) + \frac{\gamma/2-pE}{\sqrt{a}} < 0$; to show this for all $E \geq E_\infty$, it suffices to show it when $E = E_\infty$. If we let $r(x) = \frac{\phi(x)}{\Phi(x)}$, which is strictly positive and strictly decreasing, we can rearrange $g(s(0)) + \frac{\gamma/2-pE_\infty}{\sqrt{a}} < 0$ to $\frac{\sqrt{p-1}}{p} r(-\gamma/(2\sqrt{a})) < z_\ast$, which we can show again by showing that $h'_E(u) > 0$ when $u = \frac{\sqrt{p-1}}{p} r(-\gamma/(2\sqrt{a}))$. After some algebra, one finds that
\[
    h'_E\left( \frac{\sqrt{p-1}}{p} r\left(-\frac{\gamma}{2\sqrt{a}}\right)\right) = \sqrt{\frac{b}{a}} \left[ r \left( \frac{p-1}{p} r \left(-\frac{\gamma}{2\sqrt{a}} \right) - \frac{\gamma}{2\sqrt{a}} \right) - \frac{1}{p} r\left(-\frac{\gamma}{2\sqrt{a}}\right) \right].
\]
Let $c = \gamma/(2\sqrt{a})$ and $\alpha = r(-c)$ and $\beta = \frac{p-1}{p} \alpha - c$. Since $r'(\cdot) > -1$, we have $r(\beta) - r(-c) > -(\beta+c) = \frac{1-p}{p}\alpha$, and thus $r(\beta) > \alpha + \frac{1-p}{p}\alpha = \alpha/p$, which means that the above display is strictly positive. This completes the argument that $u_\ast(E,\gamma) > 0$.

Let
\[
    j_E(u) = \inf_{t \geq 0} k_E(u,t).
\]
Our goal is to compute $\sup_{u \geq 0} j_E(u)$. Since $u \mapsto g(s(u)) + \frac{\gamma/2-pE}{\sqrt{a}} = \partial_t k_E(u,0)$ is strictly increasing in $u$, and it vanishes by defition at $u = u_\ast$, we now have two cases:
\begin{itemize}
\item \textbf{Case 1 ($u \geq u_\ast(E,\gamma)$):} In this case, $\partial_t k_E(u,0) \geq 0$, so 
\[
    j_E(u) = \inf_{t \geq 0} k_E(u,t) = k_E(u,0) = h_E(u).
\]
Rearranging \eqref{eqn:e-infinity-manipulation-1} and \eqref{eqn:e-infinity-manipulation-2}, we find that, when $E > E_\infty$, we have that $g(s(z_\ast)) + \frac{\gamma/2-pE}{\sqrt{a}} < 0$; from the definition of $u_\ast$, this implies that $u_\ast > z_\ast$. Since $h_E(u)$ is strictly concave in $u$ and optimized at $u = z_\ast$, we obtain
\[
    \sup_{u \geq u_\ast(E,\gamma)} j_E(u) = \sup_{u \geq u_\ast(E,\gamma)} h_E(u) = h_E(u_\ast(E,\gamma)),
\]
as well as that $j_E$ has a strictly negative right-hand derivative at $u = u_\ast(E,\gamma)$, as long as $E > E_\infty$ strictly.
\item \textbf{Case 2 ($0 \leq u < u_\ast(E,\gamma)$):} In this case, $\partial_t k_E(u,0) < 0$, so there exists a unique $t_\ast = t_\ast(u) = t_\ast(u,E,\gamma) \in (0,\infty)$ solving the first-order condition
\begin{equation}
\label{eqn:first-order-in-t-specific}
    0 = \partial_t k_E(u,t_\ast(u)) = g(t_\ast(u) + s(u)) + \frac{\gamma/2-pE}{\sqrt{a}}. 
\end{equation}
Thus $j_E(u) = k_E(u,t_\ast(u))$, so we are interested in $\sup_{0 \leq u \leq u_\ast(E,\gamma)} j_E(u)$. We compute
\begin{equation}
\label{eqn:j-derivative}
\begin{split}
    \frac{\diff}{\diff u} j_E(u) &= \frac{\diff}{\diff u} k_E(u,t_\ast(u)) = \partial_u k_E(u,t_\ast(u)) + \partial_t k_E(u,t_\ast(u)) \cdot \partial_u t_\ast(u) \\
    &= \partial_u k_E(u,t_\ast(u)) = \frac{\phi(t_\ast(u)+s(u))}{\Phi(t_\ast(u)+s(u))} s'(u) - u + t_\ast(u)\sqrt{\frac{b}{a}} \\
    &= \sqrt{\frac{b}{a}} g(t_\ast(u)+s(u)) - \sqrt{\frac{b}{a}} s(u) - u \\
    &= \sqrt{\frac{b}{a}} \left( \frac{pE}{\sqrt{a}} - \frac{\gamma}{2\sqrt{a}} - s(u)\right) - u = \frac{pE\sqrt{b}}{a} - \left(\frac{b}{a}+1\right)u \\
    &= p \left( \frac{E\sqrt{b}}{a} - u \right),
\end{split}
\end{equation}
where we used twice that $t_\ast(u)$ satisfies \eqref{eqn:first-order-in-t-specific}, and used that $b/a+1=p$. This will tell us that $\sup_{0 \leq u < u_\ast(E,\gamma)} k_E(u,t_\ast(u))$ is achieved at $u = \widetilde{u} \defeq \widetilde{u}(E) \defeq E\sqrt{b}/a$, as long as we can show that this value of $u$ actually lies in the interval $0 \leq u \leq u_\ast(E,\gamma)$, for which the only nontrivial part is the upper bound. To show this, we want to show that $\widetilde{u}$ is at most the right-hand side of \eqref{eqn:u-ast}, which we can rearrange (using again that $b/a = p-1$) into the statement that
\[
    g\left( \frac{(p-1)E - \gamma/2}{\sqrt{a}} \right) \leq \frac{pE-\gamma/2}{\sqrt{a}}
\]
or, in other words, that
\begin{equation}
\label{eqn:finding-einfty-again!}
    \frac{\phi\left(\frac{(p-1)}{\sqrt{a}} E - \frac{\gamma}{2\sqrt{a}}\right)}{\Phi\left(\frac{(p-1)}{\sqrt{a}} E - \frac{\gamma}{2\sqrt{a}}\right)} \leq \frac{E}{\sqrt{a}}.
\end{equation}
However, $E_\infty$ satisfies this as an \emph{equality}, since it is defined (see \eqref{eqn:e-infty-definition}) in terms of $z_\ast$, which satisfies the first-order condition \eqref{eqn:z-ast-first-order}. Since the left-hand side of \eqref{eqn:finding-einfty-again!} is strictly decreasing in $E$, while the right-hand side is strictly increasing in $E$, we conclude that \eqref{eqn:finding-einfty-again!} holds for all $E \geq E_\infty$.
\end{itemize}
This tells us several things. The first is that $\widetilde{u}$ is, for $E > E_\infty$, in the interior of the interval $[0,u_\ast(E,\gamma)]$; comparing to \eqref{eqn:j-derivative}, we see that $j_E$ has a strictly negative left-hand derivative at $u_\ast(E,\gamma)$. We previously showed that $j_E$ has a strictly negative right-hand derivative at $u_\ast(E,\gamma)$; since $t_\ast(u_\ast(E,\gamma)) = 0$ (see \eqref{eqn:u-ast-implicit}), it is also continuous at $u_\ast(E,\gamma)$, and thus $\sup_{u \geq 0} j_E(u)$ must actually be achieved on the interval $[0,u_\ast(E,\gamma)]$.

The second is that 
\begin{equation}
\label{eqn:sigma-p-achieved-on-left}
    \Sigma_p(E,\gamma) = \sup_{0 \leq u \leq u_\ast(E,\gamma)} j_E(u) = j_E\left(\widetilde{u}\right) = k_E\left(\widetilde{u},t_\ast\left(\widetilde{u}\right)\right).
\end{equation}
From \eqref{eqn:first-order-in-t-specific}, one can write $t_\ast(u)$ as an explicit function using $g^{-1}$; thus it is continuous in $u$, and hence $t_\ast(\widetilde{u})$ is a continuous function of $E > E_\infty$, including a right-continuous extension $\lim_{E \downarrow E_\infty} t_\ast(\widetilde{u}(E)) = 0$. Since $\widetilde{u}(E)$ is also continuous in $E$, and $k_E$ is continuous in $E$, we conclude that $\Sigma_p(E,\gamma)$ is continuous at $E_\infty(p,\gamma)$.

Next, we try to simplify \eqref{eqn:sigma-p-achieved-on-left} as in the statement of Lemma \ref{lem:e-infinity}. The equation \eqref{eqn:energy-constrained-implicit} has a unique solution, since the left-hand side is $g(x)$ (which is strictly increasing in $x$) and the right-hand side is positive. Comparing to \eqref{eqn:first-order-in-t-specific}, we see that $x = t_\ast(u) + s(u)$ for all $u < u_\ast(E,\gamma)$, and in particular for $u = \widetilde{u}$. We also use \eqref{eqn:first-order-in-t-specific} to find that
\[
    0 = g(t_\ast(\widetilde{u})+s(\widetilde{u})) = r(t_\ast(\widetilde{u}) + s(\widetilde{u})) + t_\ast(\widetilde{u}) - \frac{E}{\sqrt{a}},
\]
so that 
\begin{equation}
\label{eqn:t-ast-vs-r-x}
    t_\ast(\widetilde{u}) - E/\sqrt{a} = -r(x).
\end{equation}
We also have $\frac{pE-\sqrt{b}\widetilde{u}}{\sqrt{a}} = E/\sqrt{a}$, and thus
\[
    -\frac{\widetilde{u}^2}{2} + \frac{t_\ast(\widetilde{u})^2}{2} - t_\ast(\widetilde{u}) \left( \frac{pE-\sqrt{b}\widetilde{u}}{\sqrt{a}}\right) = \frac{(t_\ast(u) - E/\sqrt{a})^2}{2} - \frac{pE^2}{2a} = \frac{(r(x))^2}{2} - \frac{p!E^2}{2},
\]
which, combined with \eqref{eqn:sigma-p-achieved-on-left}, yields \eqref{eqn:energy-constrained-explicit}.

Next we show that $\Sigma_p(E,\gamma)$ is \emph{strictly} below $\Sigma^{\textup{mix}}(p,0,\gamma)$ for $E > E_\infty$; for this, we use that $t_\ast(\widetilde{u}) > 0$, and thus
\begin{align*}
    \Sigma_p(E,\gamma) &= k_E\left(\widetilde{u},t_\ast\left(\widetilde{u}\right)\right) < k_E\left(\widetilde{u},0\right) = h_E\left(\widetilde{u}\right) \leq \sup_{u \geq 0} h_E(u) = \Sigma^{\textup{mix}}(p,0,\gamma).
\end{align*}
The only remaining claim is that $\Sigma_p$ is strictly decreasing for $E > E_\infty$. For this, we note first that $r'(x) = -r(x)g(x)$, and that $g'(x) = 1-r(x)g(x) \in (0,1)$. Thus, writing $x(E)$ for the solution to \eqref{eqn:energy-constrained-implicit}, we differentiate that equation to find that $\frac{p}{\sqrt{a}} = g'(x(E))x'(E)$. Differentiating \eqref{eqn:energy-constrained-explicit}, we see that
\begin{align*}
    \partial_E \Sigma_p(E,\gamma) &= r(x(E))x'(E) + r(x(E))r'(x(E))x'(E) - p!E \\
    &= r(x(E))x'(E)g'(x(E)) - p!E = \frac{p}{\sqrt{a}} r(x(E)) - p!E.
\end{align*}
This has the same sign as $r(x(E)) - E/\sqrt{a}$. But \eqref{eqn:t-ast-vs-r-x} says this equals $-t_\ast(\widetilde{u}(E))$, which is negative for all $E > E_\infty$, completing the proof. 
\end{proof}

\begin{lem}
\label{lem:wlln-application}
Fix $E < E_\infty(p,\gamma)$. There exists $\delta > 0$ such that
\[
    \liminf_{N \to \infty} \inf_{u \in [z_\ast-\delta,z_\ast+\delta]} \frac{1}{N} \log p_N(u,E) \geq 0.
\]
\end{lem}
\begin{proof}
Let $Y_i = Y_i^{(N,u)}$ have the distribution of $X_i$ conditioned on $X_i \geq -s_N(u)$, and let $Z_i = Y_i - \E[Y_i] = Y_i - g(s_N(u)) + s_N(u)$. Then
\begin{align*}
    p_N(u,E) = \P\left(\frac{1}{N} \sum_{i=1}^N Y_i \geq \frac{pE-\sqrt{Nb_N}u}{\sqrt{a_N}}\right) = \P\left( \frac{1}{N} \sum_{i=1}^N Z_i \geq \frac{pE-\gamma/2}{\sqrt{a_N}} - g(s_N(u)) \right).
\end{align*}
Recall from \eqref{eqn:e-infinity-manipulation-1} and \eqref{eqn:e-infinity-manipulation-2} that $\frac{pE_\infty-\gamma/2}{\sqrt{a_N}} - g(s(z_\ast)) = 0$. Since $E < E_\infty$ is strict, we can choose small enough $\delta > 0$ and $\epsilon' > 0$ such that 
\[
    \frac{pE-\gamma/2}{\sqrt{a_N}} - g(s_N(u)) \leq -\epsilon'
\]
uniformly over $u \in [z_\ast-\delta,z_\ast+\delta]$. Thus
\[
    \sup_{u \in [z_\ast-\delta,z_\ast+\delta]} (1-p_N(u,E)) \leq \P\left(\abs{\frac{1}{N}\sum_{i=1}^N Z_i} \geq \epsilon'\right) \leq \frac{\E[Z_i^2]}{N(\epsilon')^2},
\]
as in the proof of the weak law of large numbers. One can compute moments of $Z_i$ from \eqref{eqn:truncated-gaussian-basics}, and show softly that $\limsup_{N \to \infty} \sup_{u \in [z_\ast-\delta,z_\ast+\delta]} \E[Z_i^2] < \infty$; this completes the proof (actually in the strong sense that the probabilities $p_N$ tend to one). 
\end{proof}

\begin{lem}
\label{lem:energy-constrained-error}
If the error $e_N(u)$ is as defined in \eqref{eqn:energy-constrained-error}, then for each $M > 0$ we have
\[
    \limsup_{N \to \infty} \sup_{\abs{u} \leq M} \abs{e_N(u)} = 0.
\]
\end{lem}
\begin{proof}
First, we claim that there exists some $K_M$ such that $\abs{t_\ast(u)} \leq K_M$ for all $\abs{u} \leq M$. Indeed, if
\[
    h_E(u,t) = \log\Phi(t+s(u)) + \frac{t^2}{2} - t\left(\frac{pE-\sqrt{b}u}{\sqrt{a}}\right),
\]
then for such $u$ we have
\[
    h_E(u,t) \geq \log \Phi(s(-M)) + \frac{t^2}{2} - t\left(\frac{pE + \sqrt{b}M}{\sqrt{a}}\right)
\]
There exists $K_M$ such that the right-hand side is positive for all $t > K_M$; since $h_E(u,0) = \log(\Phi(s(u))) \leq 0$, this means that the infimizer in $t$ must indeed lie in $[-K_M,K_M]$.

Additionally, there exists another compact set containing all the values of $s(u)$ and $s_N(u)$ for all $N$ and all $\abs{u} \leq M$; thus, since $\log \Phi$ is locally Lipschitz, we have
\[
    \sup_{\abs{u} \leq M} \abs{\log\Phi(t_\ast(u) + s_N(u)) - \log\Phi(t_\ast(u) + s(u))} \leq C_M \sup_{\abs{u} \leq M} \abs{s_N(u)-s(u)}
\]
which tends to zero. The remaining terms in the error are even easier, since $\abs{\frac{1}{\sqrt{a_N}} - \frac{1}{\sqrt{a}}} \to 0$ and $\abs{\sqrt{\frac{Nb_N}{a_N}} - \sqrt{\frac{b}{a}}} \to 0$.
\end{proof}

\begin{lem}
\label{lem:pointwise-cramer}
Fix $E > E_\infty(p,\gamma)$. Then for any given $\widetilde{u}$ satisfying \eqref{eqn:lln-boundary} we have
\[
    \liminf_{\delta \to 0} \liminf_{N \to \infty} \inf_{u \in [\widetilde{u}-\delta,\widetilde{u}+\delta]} \left( \frac{1}{N} \log p_N(u,E) + I\left( \frac{pE - \sqrt{b}u}{\sqrt{a}} \right)\right) \geq 0.
\]
\end{lem}
\begin{proof}
Let $Y_i = Y_i^{(u)}$ have the distribution of $X_i$ conditioned on $X_i \geq -s_N(u)$, let
\begin{align*}
    \nu_N(u,\theta) &= \frac{\theta^2}{2} + \log\Phi(\theta + s_N(u)) - \log\Phi(s_N(u)),
\end{align*}
and consider the tilted probability measure
\[
    \frac{\diff \P_{N,u,\theta}}{\diff \P}(Y_1,\ldots,Y_N) = \exp\left(\theta\sum_{i=1}^N Y_i - N\nu_N(u,\theta)\right)
\]
The variables $Y_1, \ldots, Y_N$ are still i.i.d. under $\P_{u,\theta}$, but their mean becomes
\[
    \partial_\theta \nu_N(u,\theta) = \theta + \frac{\phi(\theta + s_N(u))}{\Phi(\theta + s_N(u))} = g(\theta + s_N(u)) - s_N(u).
\]
This is a strictly increasing function of $\theta$, ranging from $g(s_N(u)) - s_N(u)$ (the mean of $Y_i$ under $\P$) to $+\infty$ as $\theta$ ranges from $0$ to $+\infty$. Thus, given any threshold $t > g(s_N(u)) - s_N(u)$, there exists a unique $\theta_N$ such that the mean of $Y_i$ under $\P_{N,u,\theta}$ is $t$. Let us fix $\epsilon > 0$ and choose $\theta = \theta_{N,u}$ such that this mean is
\[
    \frac{pE-\sqrt{Nb_N}u}{\sqrt{a_N}} + \epsilon.
\]
(Since $\tilde{u}$ satisfies \eqref{eqn:lln-boundary}, one can choose $\delta$ small enough that this indeed lies eventually above $g(s_N(u)) - s_N(u)$, uniformly over $u \in [\widetilde{u}-\delta,\widetilde{u}+\delta]$.) For $\eta > 0$, consider the event
\[
    \mc{E}_N(\eta) = \left\{pE + \eta \geq \sqrt{a_N} \left( \frac{1}{N} \sum_{i=1}^N Y_i \right) + \sqrt{Nb_N} u \geq pE\right\}.
\]
We will abusively write $\mc{E}_N(+\infty)$ for the event $\{\sqrt{a_N}(N^{-1}\sum_{i=1}^N Y_i) + \sqrt{Nb_N}(u) \geq pE\}$. Consider the estimate
\begin{align*}
    p_N(u,E) &= \P(\mc{E}_N(+\infty)) \geq \P(\mc{E}_N(\eta)) \\
    &= \E_{N,u,\theta_{N,u}} \left[ \exp\left( - \theta_{N,u} \sum_{i=1}^N Y_i + N\nu_N(u,\theta_{N,u}) \right) \mathds{1}\{\mc{E}_N(\eta)\} \right] \\
    &\geq \exp\left\{ - N \left( \theta_{N,u} \left( \frac{pE+\eta-\sqrt{Nb_N}u}{\sqrt{a_N}}\right) - \nu_N(u,\theta_{N,u}) \right) \right\} \P_{N,u,\theta_{N,u}}(\mc{E}_N(\eta)).
\end{align*}
Thus, if we can show that, for every $\eta > \sqrt{a}\epsilon$, we have
\begin{equation}
\label{eqn:cramer-tilt}
    \liminf_{N \to \infty} \inf_{u \in [\widetilde{u}-\delta,\widetilde{u}+\delta]} \frac{1}{N} \log \P_{N,u,\theta_{N,u}}(\mc{E}_N(\eta)) \geq 0.
\end{equation}
and 
\begin{equation}
\label{eqn:cramer-rate-technical}
\begin{split}
    &\liminf_{\eta \to 0} \liminf_{\epsilon \to 0} \liminf_{\delta \to 0} \liminf_{N \to \infty} \inf_{u \in [\widetilde{u}-\delta,\widetilde{u}+\delta]} \\
    &\left( -\theta_{N,u} \left( \frac{pE+\eta-\sqrt{Nb_N}u}{\sqrt{a_N}} \right) + \nu_N(u,\theta_{N,u}) + I\left(\frac{pE-\sqrt{b}u}{\sqrt{a}}\right) \right) \geq 0,
\end{split}
\end{equation}
we can combine them to conclude. Thus it remains only to check \eqref{eqn:cramer-tilt} and \eqref{eqn:cramer-rate-technical}.

Let $Z_i = Y_i - \frac{pE-\sqrt{Nb_N}u}{\sqrt{a_N}} - \epsilon$. These are i.i.d. and centered under $\P_{N,u,\theta_{N,u}}$, and we have
\[
    \P_{N,u,\theta_{N,u}}(\mc{E}_N(\eta)) = \P_{N,u,\theta_{N,u}}\left( \frac{\eta}{\sqrt{a_N}} - \epsilon \geq \frac{1}{N} \sum_{i=1}^N Z_i \geq -\epsilon \right).
\]
Since $a_N \to a$ and $\eta > \sqrt{a}\epsilon$, there exists $\epsilon' > 0$ such that $[-\epsilon',\epsilon'] \subseteq [-\epsilon,\eta/\sqrt{a_N} - \epsilon]$ for large enough $N$. Then, as in the weak law of large numbers, we have
\[
    \P_{N,u,\theta_{N,u}}(\mc{E_N}(\eta)^c) \leq \P_{N,u,\theta_{N,u}}\left( \abs{\frac{1}{N} \sum_{i=1}^N Z_i } > \epsilon'\right) \leq \frac{\E_{N,u,\theta_{N,u}}[Z_i^2]}{N(\epsilon')^2}.
\]
Notice that
\[
    \E_{N,u,\theta_{N,u}}[Z_i^2] = \Var_{N,u,\theta_{N,u}}(Y_i) = \partial^2_{\theta} \nu_N(u,\theta_{N,u})= g'(\theta_{N,u} + s_N(u)) \leq 1,
\]
simply because $g'(\cdot) \leq 1$ on $\R$; thus we conclude \eqref{eqn:cramer-tilt} (actually in the stronger sense that the probabilities converge to one, rather than just being not exponentially small). 

For the remainder, we note from the proof of Lemma \ref{lem:e-infinity} and specifically from \eqref{eqn:first-order-in-t-specific} that 
\[
    I\left(\frac{pE-\sqrt{b}u}{\sqrt{a}}\right) = t_\ast(u) \left(\frac{pE-\sqrt{b}u}{\sqrt{a}}\right) - \frac{t_\ast(u)^2}{2} - \log\Phi\left(t_\ast(u) + s(u)\right) + \log\Phi(s(u)),
\]
where $t_\ast(u)$ solves the first-order condition 
\[
    g(t_\ast(u) + s(u)) = \frac{pE-\sqrt{b}u}{\sqrt{a}} + s(u) = \frac{pE-\gamma/2}{\sqrt{a}}.
\]
On the other hand, we know that $\theta_{N,u}$ solves
\[
    g(\theta_{N,u}+s_N(u)) = \frac{pE-\sqrt{Nb_N}u}{\sqrt{a_N}} + \epsilon + s_N(u) = \frac{pE-\gamma/2}{\sqrt{a_N}} + \epsilon.
\]
We have $\lim_{x \to +\infty} g'(x) = \lim_{x \to -\infty} g'(x) = 0$, but $g'$ is lower-bounded on finite intervals; this tells us that for every $u_\ast$ there exists a constant $C$ such that
\begin{equation}
\label{eqn:t-ast-vs-theta-n}
    \limsup_{N \to \infty} \sup_{u \in [\widetilde{u}-1,\widetilde{u}+1]} \abs{t_\ast(u) - \theta_{N,u}} \leq C\epsilon,
\end{equation}
and moreover that $\lim_{N \to \infty} \theta_{N,u}$ exists and is continuous in $u$. Thus
\begin{align*}
    &\abs{-\theta_{N,u} \left( \frac{pE+\eta-\sqrt{Nb_N}u}{\sqrt{a_N}} \right) - \nu_N(u,\theta_{N,u}) + I\left(\frac{pE-\sqrt{b}u}{\sqrt{a}}\right)} \\
    &\leq (I) + (II) + (III) + (IV) + (V), 
\end{align*}
where
\begin{align*}
    (I) &= \theta_{N,u} \eta/\sqrt{a_N}, \\
    (II) &= \abs{t_\ast(u) \left( \frac{pE-\sqrt{b}u}{\sqrt{a}} \right) - \theta_{N,u} \left( \frac{pE - \sqrt{Nb_N}u}{\sqrt{a_N}} \right)}, \\
    (III) &= \abs{\frac{t_\ast(u)^2}{2} - \frac{\theta_{N,u}^2}{2}}, \\
    (IV) &= \abs{\log\Phi(t_\ast(u) + s(u)) - \log\Phi(\theta_{N,u} + s_N(u))}, \\
    (V) &= \abs{\log\Phi(s(u)) - \log\Phi(s_N(u))}.
\end{align*}
Using that $\log\Phi$ has a lower-bounded derivative on compact sets, that $s_N(u) \to s(u)$ locally uniformly in $u$, that $Nb_N \to b$ and $a_N \to a$, the fact that $\lim_{N \to \infty} \theta_{N,u}$ exists and is continuous in $u$, as well as \eqref{eqn:t-ast-vs-theta-n}, completes the proof of \eqref{eqn:cramer-rate-technical} and thus completes the proof. 
\end{proof}

\bibliography{complexitybib-arxiv}

\end{document}